\documentclass[a4paper,12pt]{article}

\usepackage{fancyhdr}
\usepackage[latin1]{inputenc}
\usepackage[usenames,dvipsnames]{pstricks}
\usepackage[permil]{overpic}
\usepackage{subfig}
\usepackage{amsmath}
\usepackage{amssymb}
\usepackage{enumerate}
\usepackage{placeins}
\usepackage{caption}
\usepackage{psfrag}
\usepackage{bm}
\usepackage{multirow}

\begin{document}

\pagestyle{fancy}

\begin{center}
{
\Large
\textbf{Identification of Nonlinear Normal Modes of Engineering Structures under Broadband Forcing}

\vspace{1cm}
\normalsize
Authors' postprint version \\
Published in : Mechanical systems and signal processing Vol.~74 (2016)
doi:10.1016/j.ymssp.2015.04.016
http://www.sciencedirect.com/science/article/pii/S0888327015001983

\vspace{1cm}
\textbf{J.P. Noël, L. Renson, C. Grappasonni, G. Kerschen}

\vspace{0.5cm}
Space Structures and Systems Laboratory\\
Aerospace and Mechanical Engineering Department\\
University of Liège, Liège, Belgium \\

\vspace{2cm}

\rule{0.85\linewidth}{.3pt}

\vspace{-0.5cm}

\begin{abstract}
The objective of the present paper is to develop a two-step methodology integrating system identification and numerical continuation for the experimental extraction of nonlinear normal modes (NNMs) under broadband forcing. The first step processes acquired input and output data to derive an experimental state-space model of the structure. The second step converts this state-space model into a model in modal space from which NNMs are computed using shooting and pseudo-arclength continuation. The method is demonstrated using noisy synthetic data simulated on a cantilever beam with a hardening-softening nonlinearity at its free end.

\vspace{1cm}

\end{abstract}

\vspace{-0.5cm}

\rule{0.85\linewidth}{.3pt}

\vspace{1cm}
Corresponding author: Jean-Philippe Noël\\
\vspace{0.5cm}
Space Structures and Systems Laboratory\\
Aerospace and Mechanical Engineering Department\\
University of Liège\\
1, Chemin des chevreuils (B52/3), 4000 Liège, Belgium \\
Email: jp.noel@ulg.ac.be

}
\end{center}

\newpage
\rhead{\thepage}
\section{Introduction}\label{Sec:Introduction}

Experimental modal analysis of linear engineering structures is now well-established and mature~\cite{Ewins_ISMA2006}. It is routinely practiced in industry, in particular during on-ground certification of aircraft and spacecraft structures~\cite{Peeters_GVT_IMAC2008,Goge_GVT_IMAC2007,Grillenbeck_IMAC2011}, using two specific approaches, namely phase resonance and phase separation methods. Phase resonance testing, also known as force appropriation, consists in exciting the normal modes of interest one at a time using a multipoint sine forcing at the corresponding natural frequency~\cite{Wright_NNM}. Conversely, in phase separation testing, several normal modes are excited simultaneously using either broadband or swept-sine forcing, and are subsequently identified using appropriate linear system identification techniques~\cite{Peeters_PolyMAX,VODM_Book}.

The existence of nonlinear behavior in dynamic testing is today a challenge the structural engineer is more and more frequently confronted with. In this context, the development of a nonlinear counterpart to experimental modal analysis would be extremely beneficial. An interesting approach to nonlinear modal testing is the so-called nonlinear resonant decay method introduced by Wright and co-workers~\cite{Atkins_NLRDM,Platten_NLRDM}. In this approach, a burst of a sine wave is applied to the structure at the undamped natural frequency of a normal mode, and enables small groups of modes coupled by nonlinear forces to be excited. A nonlinear curve fitting in modal space is then carried out using the restoring force surface method. The identification of modes from multimodal nonlinear responses has also been attempted in the past few years. For that purpose, advanced signal processing techniques have been utilized, including the empirical mode decomposition~\cite{Poon_EMD_NNM,Vakakis_EMD,Vakakis_EMD_FrictionBolt}, time-frequency analysis tools~\cite{Pai_TF_NNM} and machine learning algorithms~\cite{Worden_NNM_IMAC2014}. Multimodal identification relying on the synthesis of frequency response functions using individual mode contributions has been proposed in Refs.~\cite{Gibert_NNM,Chong_NNM}. The difficulty with these approaches is the absence of superposition principle in nonlinear dynamics, preventing the response of a nonlinear system from being decomposed into the sum of different modal responses.

In the present study, we adopt the framework offered by the theory of nonlinear normal modes (NNMs) to perform experimental nonlinear modal analysis. The concept of normal modes was generalized to nonlinear systems by Rosenberg in the 1960s~\cite{Rosenberg_NNM1,Rosenberg_NNM2} and by Shaw and Pierre in the 1990s~\cite{ShawPierre_NNM}. NNMs possess a clear conceptual relation with the classical linear normal modes (LNMs) of vibration, while they provide a solid mathematical tool for interpreting a wide class of nonlinear dynamic phenomena, see, \textit{e.g.}, Refs.~\cite{Vakakis_NNMBook,Lacarbonara_NNM_Asymptotic,Touze_NNM_GeomNL,NNM_Part1}. There now exist effective algorithms for their computation from mathematical models~\cite{Arquier_NNM,NNM_Part2,Laxalde_NNM_Bladings,Renson_NNM_FEM}. For instance, the NNMs of full-scale aircraft and spacecraft structures and of a turbine bladed disk were computed in Refs.~\cite{Kerschen_NNM_Paris,Renson_SmallSat,Krack_NNM}, respectively. 

A nonlinear phase resonance method exploiting the NNM concept was first proposed in Ref.~\cite{PRM_Part1}, and was validated experimentally in Ref.~\cite{PRM_Part2}. Following the philosophy of force appropriation and relying on a nonlinear generalization of the phase lag quadrature criterion, this nonlinear phase resonance method excites the targeted NNMs one at a time using a multipoint, multiharmonic sine forcing. The energy-dependent frequency and modal curve of each NNM are then extracted directly from the experimental time series by virtue of the invariance principle of nonlinear oscillations. Applications of nonlinear phase resonance testing to moderately complex experimental structures were recently reported in the technical literature, in the case of a steel frame in Ref.~\cite{Zapico_NNM} and of a circular perforated plate in Ref.~\cite{Allen_NNM_IMAC2014}.

The identification of NNMs from broadband data represents a distinct challenge in view of the absence of superposition principle in nonlinear dynamics. Indeed, the measured responses cannot merely be decomposed into a sum of individual NNM contributions. To address this challenge, the present paper develops a two-step methodology integrating system identification and numerical continuation for the experimental extraction of NNMs under broadband forcing. The first step processes acquired input and output data using the frequency-domain nonlinear subspace identification (FNSI) method~\cite{Noel_FNSI} to derive an experimental state-space model of the structure. The second step converts this state-space model into a model in modal space from which NNMs are computed using shooting and pseudo-arclength continuation~\cite{NNM_Part2}. It should be noted that identification and continuation tools others than FNSI and pseudo-arclength may also qualify for the present framework. However, the two latter are adopted because of their accuracy and applicability to real-life structures.

The paper is organized as follows. The fundamental properties of NNMs defined as periodic solutions of the underlying undamped system are briefly reviewed in Section~\ref{Sec:NNMs}. The existing nonlinear phase resonance method introduced in Ref.~\cite{PRM_Part1} is also described. In Section~\ref{Sec:NPSM}, the two building blocks of the proposed NNM identification methodology, namely the FNSI method and the pseudo-arclength continuation algorithm, are presented. The methodology is demonstrated in Section~\ref{Sec:Demo} using noisy synthetic data simulated on a cantilever beam with a hardening-softening nonlinearity at its free end. Since it can be viewed as a nonlinear generalization of linear phase separation techniques, the proposed methodology is also compared in Section~\ref{Sec:Comparison} with the previously-developed nonlinear phase resonance method. The conclusions of the study are finally summarized in Section~\ref{Sec:Conclusion}.

\section{Brief review of nonlinear normal modes (NNMs) and identification using phase resonance}\label{Sec:NNMs}

In this work, an extension of Rosenberg's definition of a NNM is considered~\cite{NNM_Part1}. Specifically, a NNM is defined as a nonnecessarily synchronous, periodic motion of the undamped, unforced, $n_{p}$-degree-of-freedom (DOF) system
\begin{equation}
\mathbf{M} \: \mathbf{\ddot{q}}(t) + \mathbf{K} \: \mathbf{q}(t) + \mathbf{f}(\mathbf{q}(t)) = 0 ,
\label{Eq:TD_Model_Undamped}
\end{equation}
where $\mathbf{M}$ and $\mathbf{K} \in \mathbb{R}^{\: n_{p} \times n_{p}}$ are the mass and linear stiffness matrices, respectively; $\mathbf{q} \in \mathbb{R}^{\: n_{p}}$ is the generalized displacement vector; $\mathbf{f}(\mathbf{q}(t)) \in \mathbb{R}^{\: n_{p}}$ is the nonlinear restoring force vector encompassing elastic terms only. This definition of a NNM may appear to be restrictive in the case of nonconservative systems. However, as shown in Refs.~\cite{NNM_Part1,Noel_SmallSat_JSV}, the topology of the underlying conservative NNMs of a system yields considerable insight into its damped dynamics.

\begin{figure}[ht]
\begin{center}
\begin{overpic}[width=140mm]
		{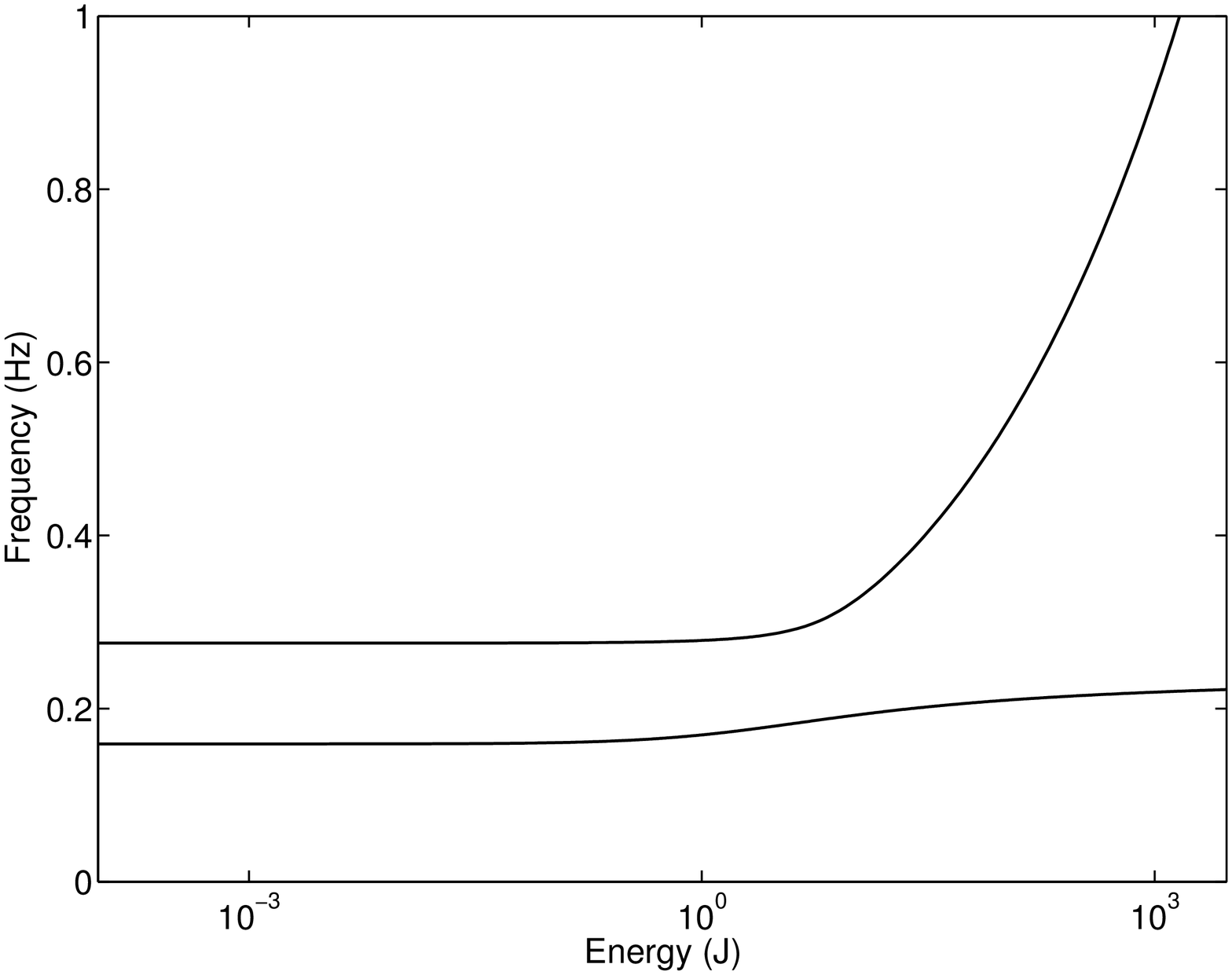}
		\put(130,292){\includegraphics[width=11.2mm]{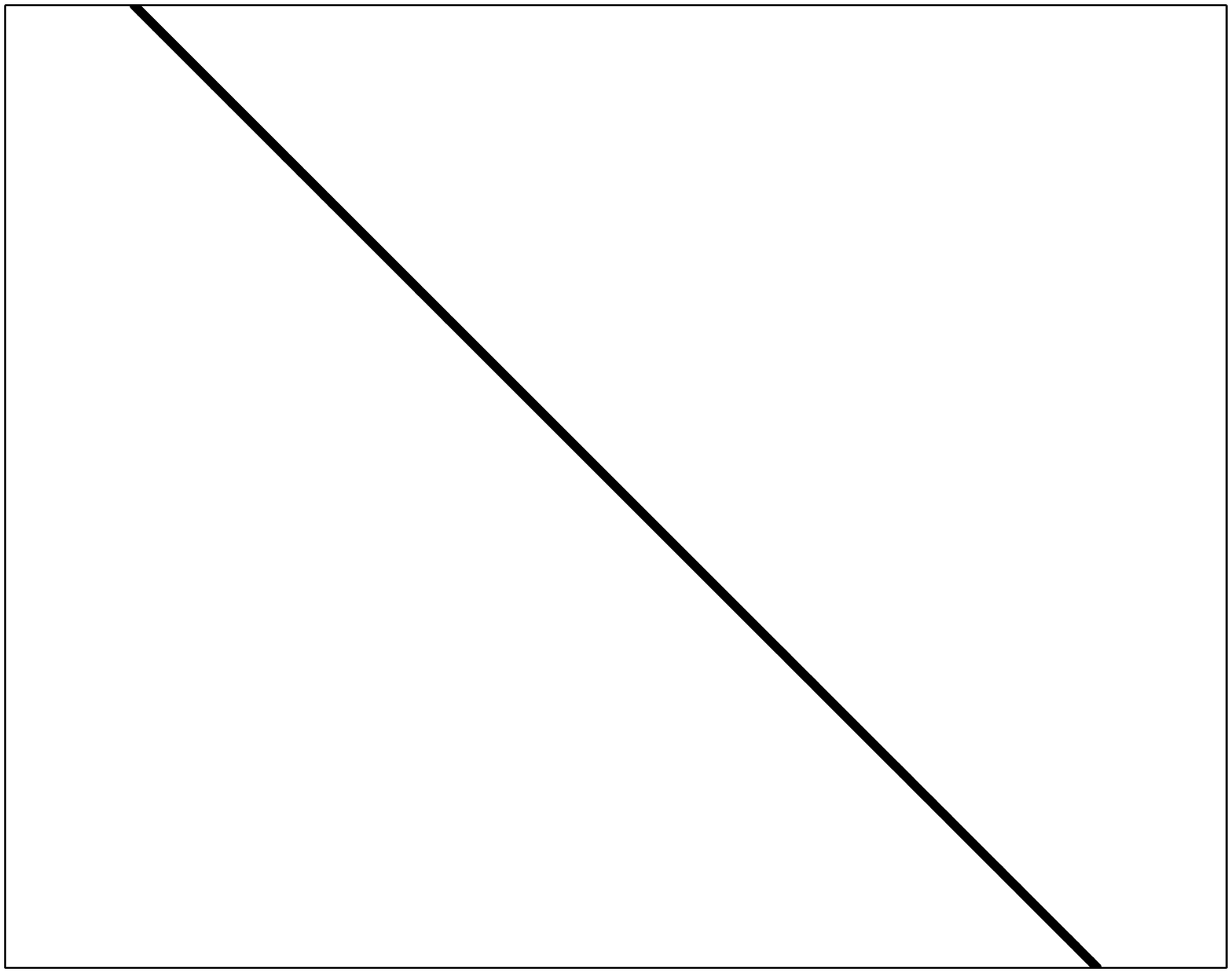}}
		\put(300,292){\includegraphics[width=11.2mm]{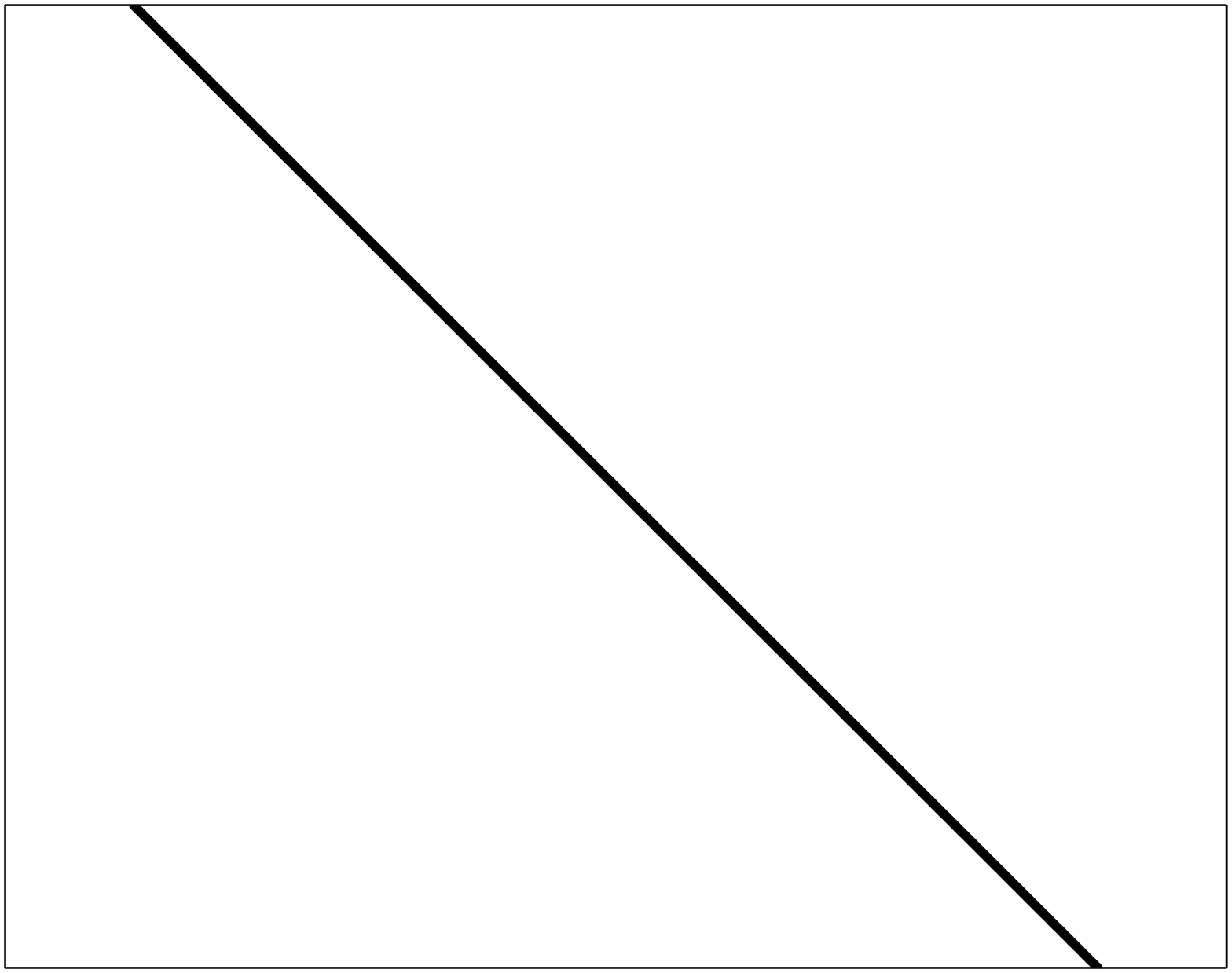}}
		\put(470,292){\includegraphics[width=11.2mm]{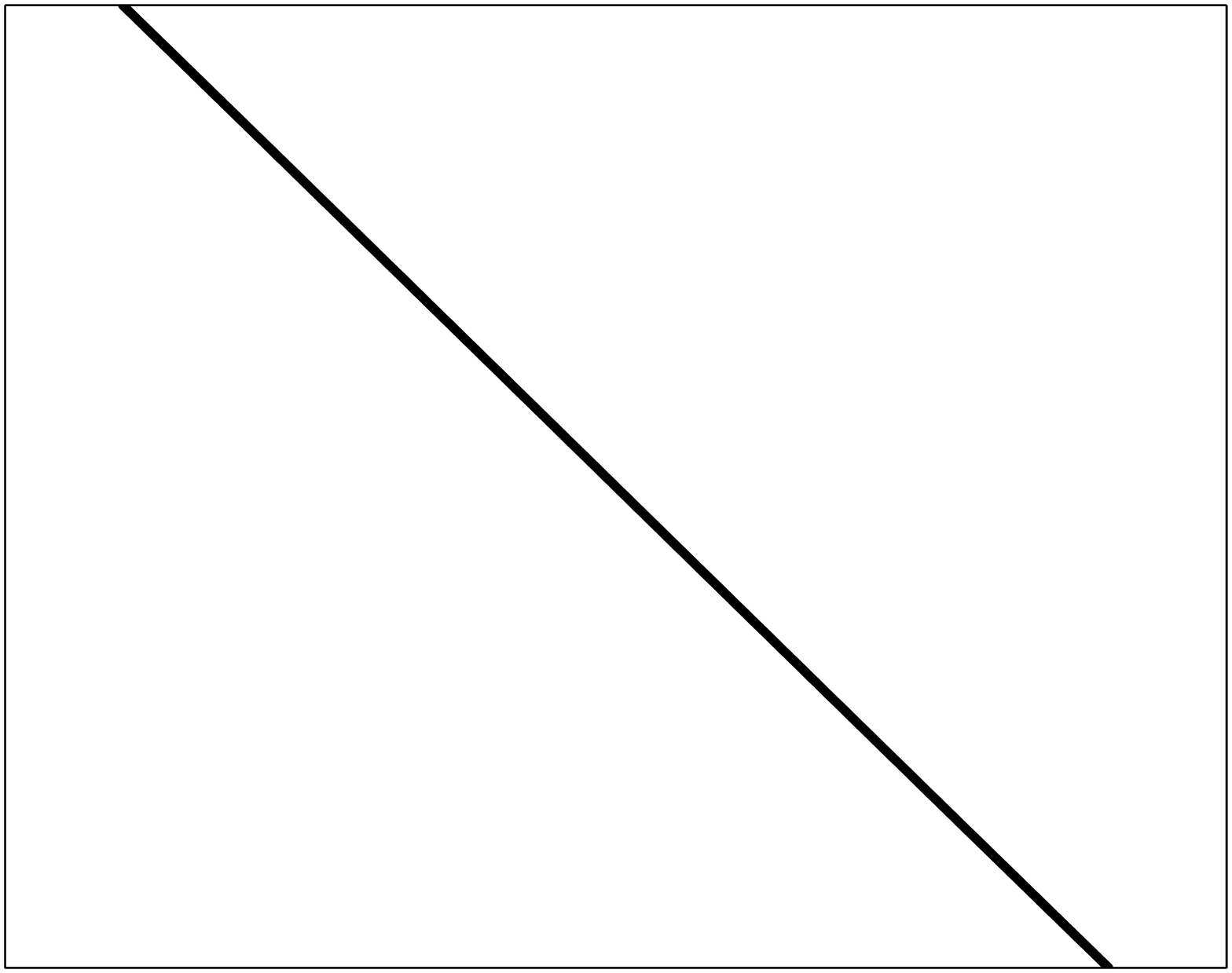}}
		\put(620,324){\includegraphics[width=11.2mm]{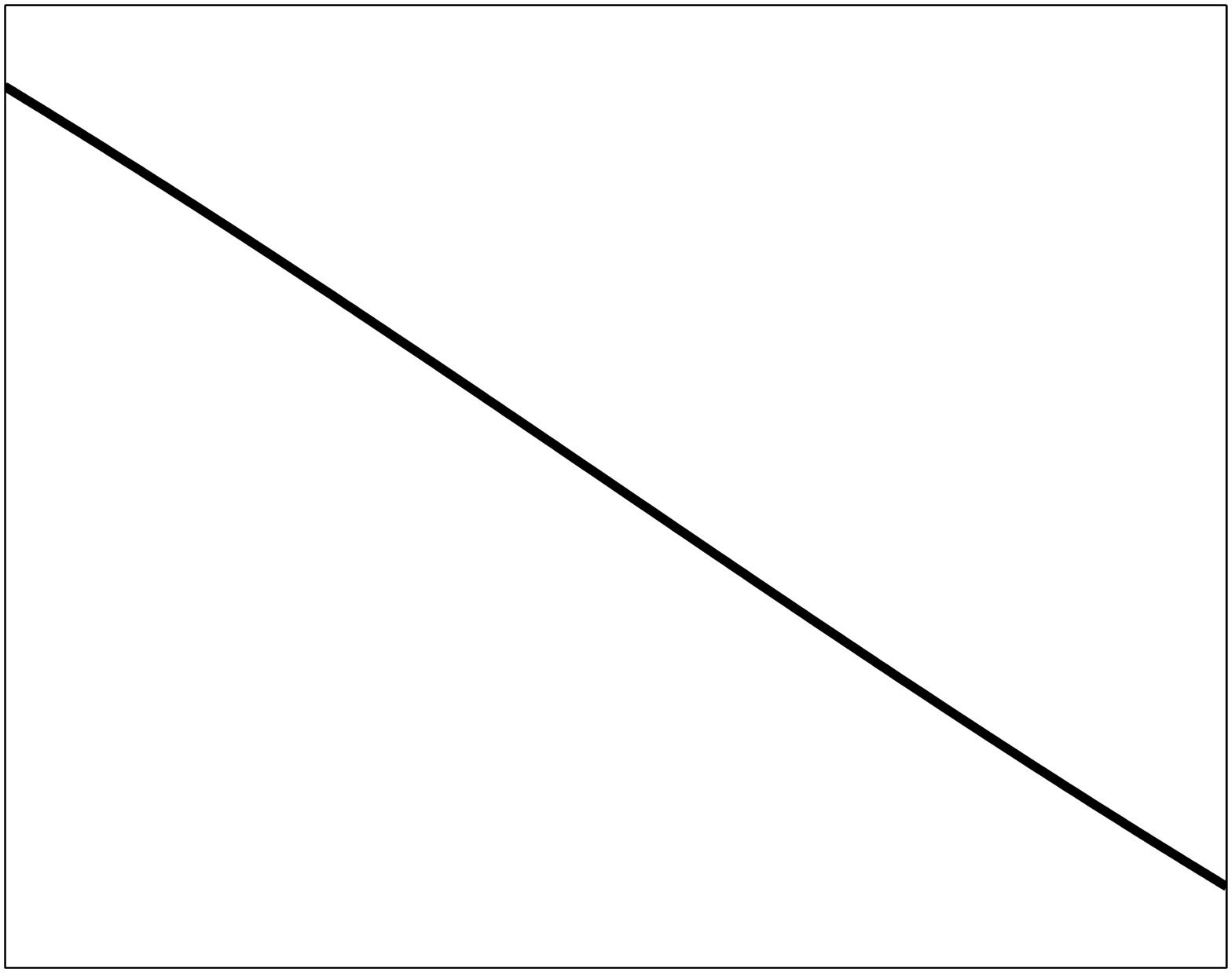}}
		\put(710,451){\includegraphics[width=11.2mm]{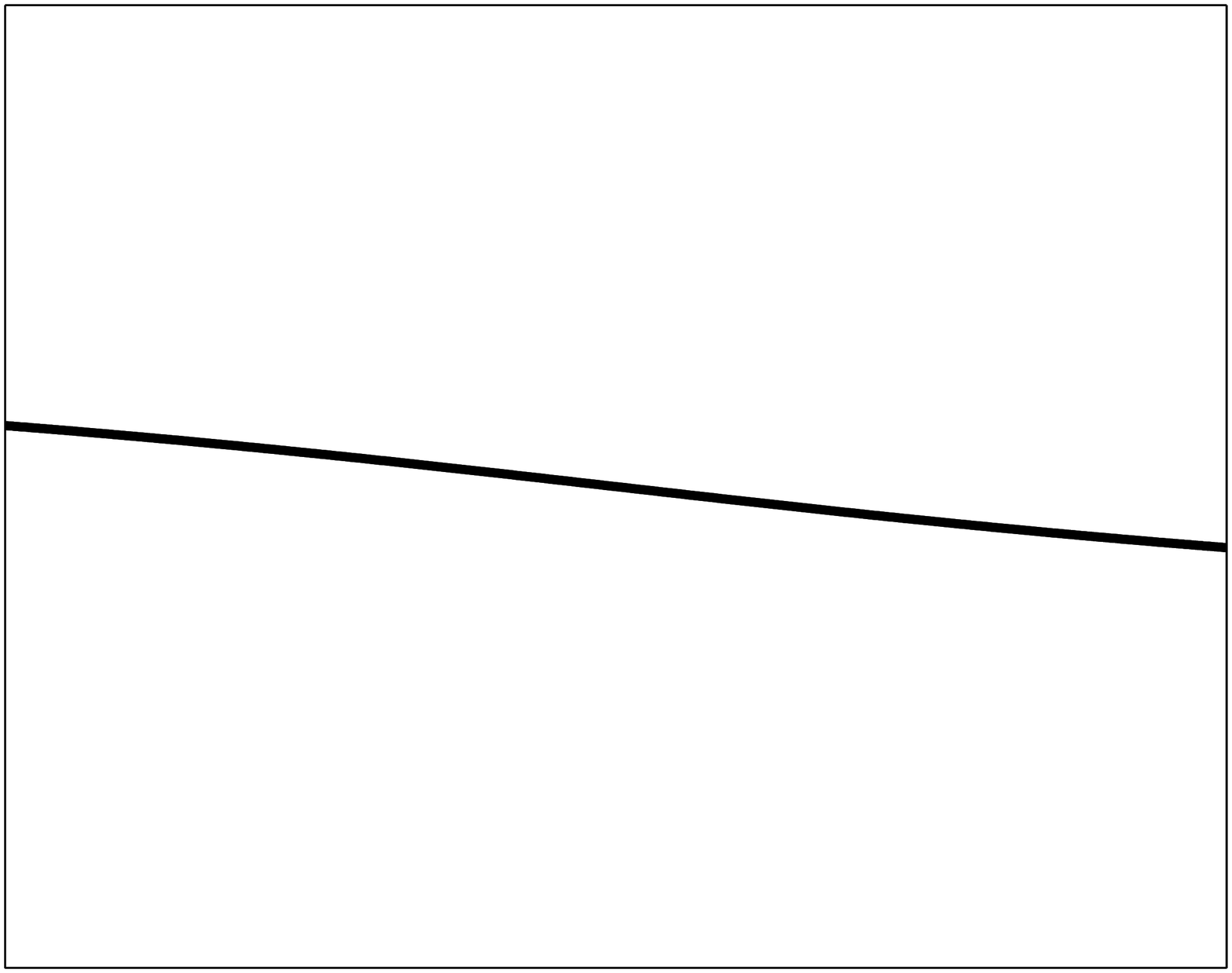}}
		\put(822,656.5){\includegraphics[width=11.2mm]{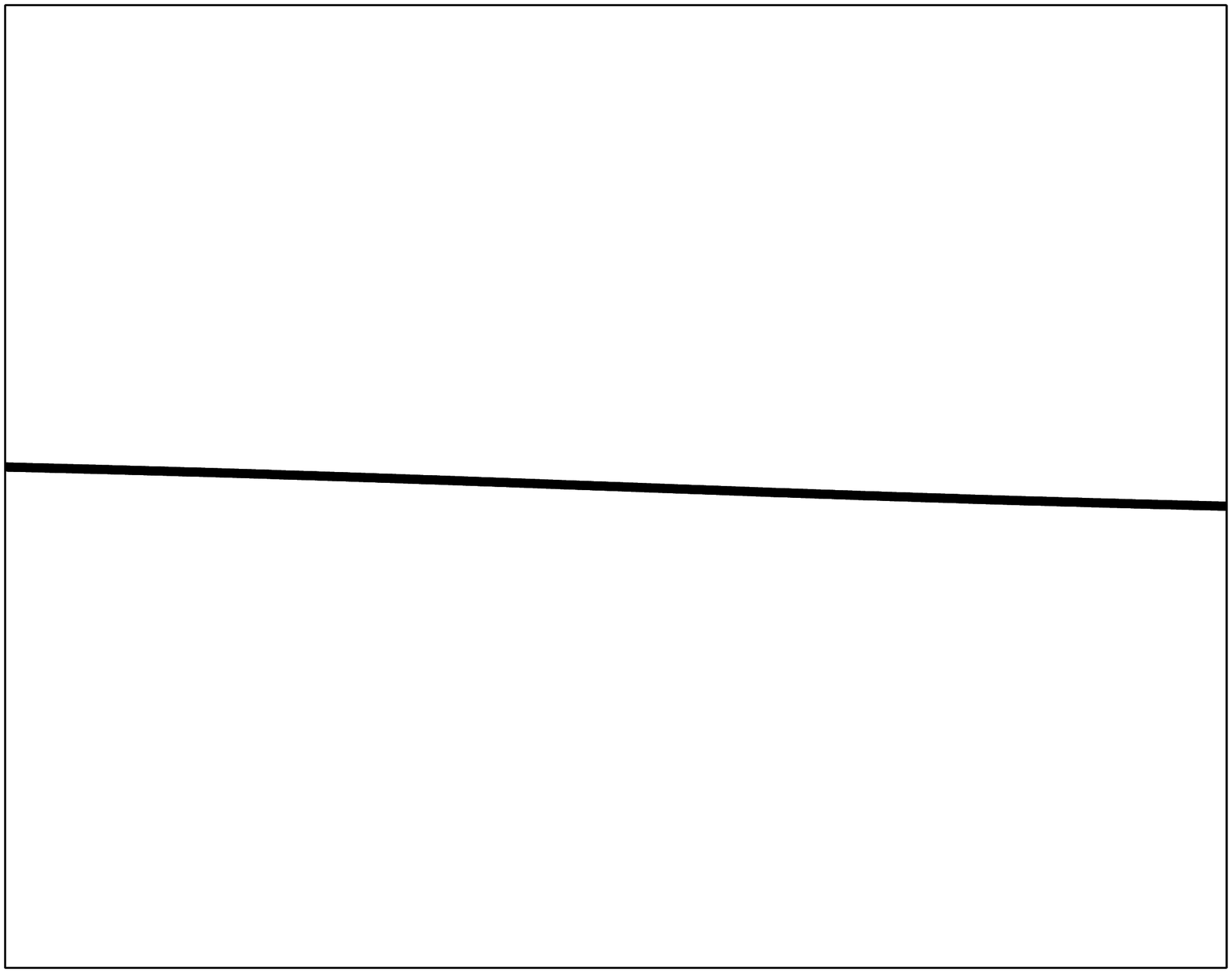}}
		\put(170,269){\pscircle(0,0){0.1}}
		\put(340,269){\pscircle(0,0){0.1}}
		\put(510,269.5){\pscircle(0,0){0.1}}
		\put(660,285){\pscircle(0,0){0.1}}
		\put(821,454){\pscircle(0,0){0.1}}
		\put(939,720.5){\pscircle(0,0){0.1}}
		\put(130,99.5){\includegraphics[width=11.2mm]{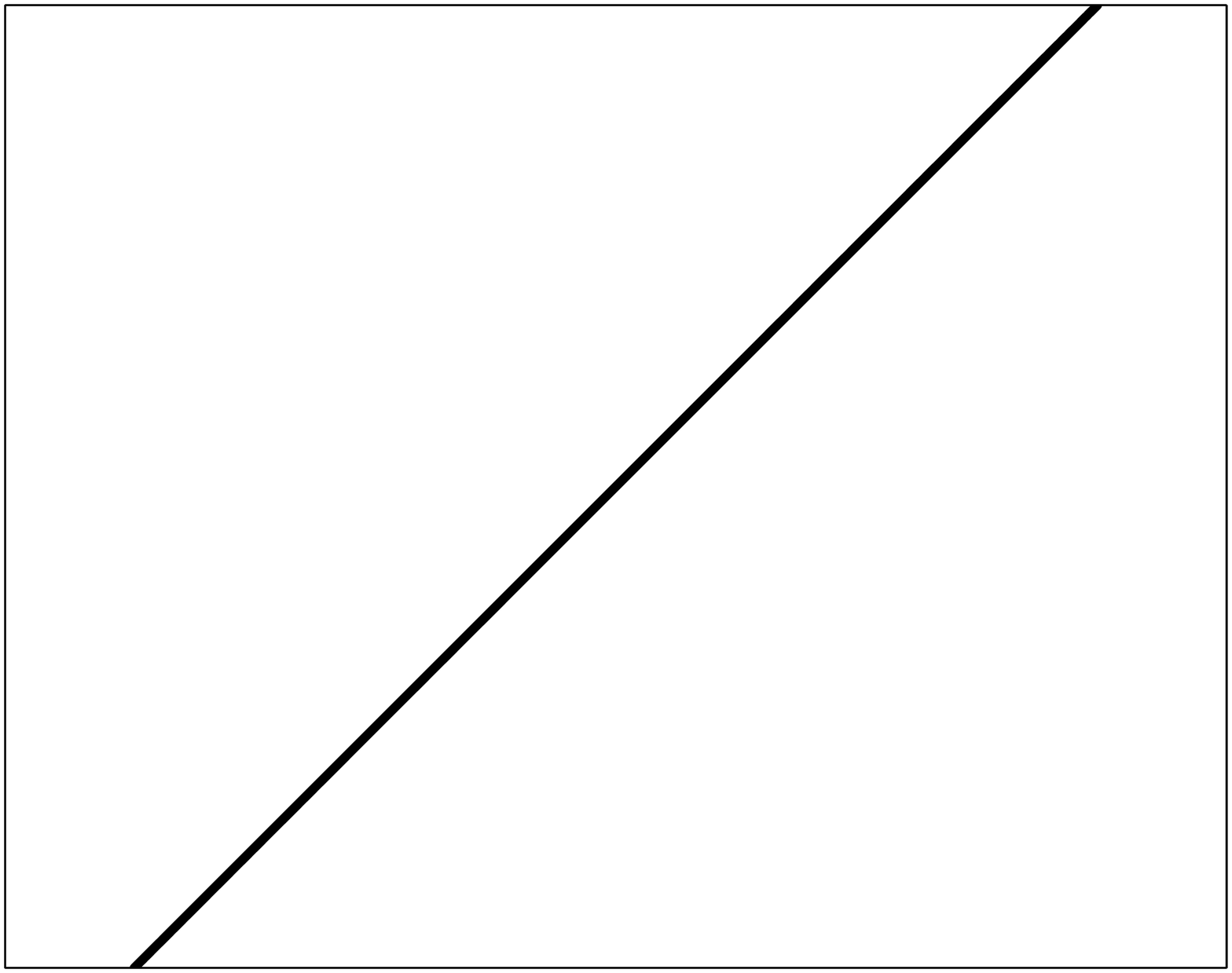}}
		\put(300,99.5){\includegraphics[width=11.2mm]{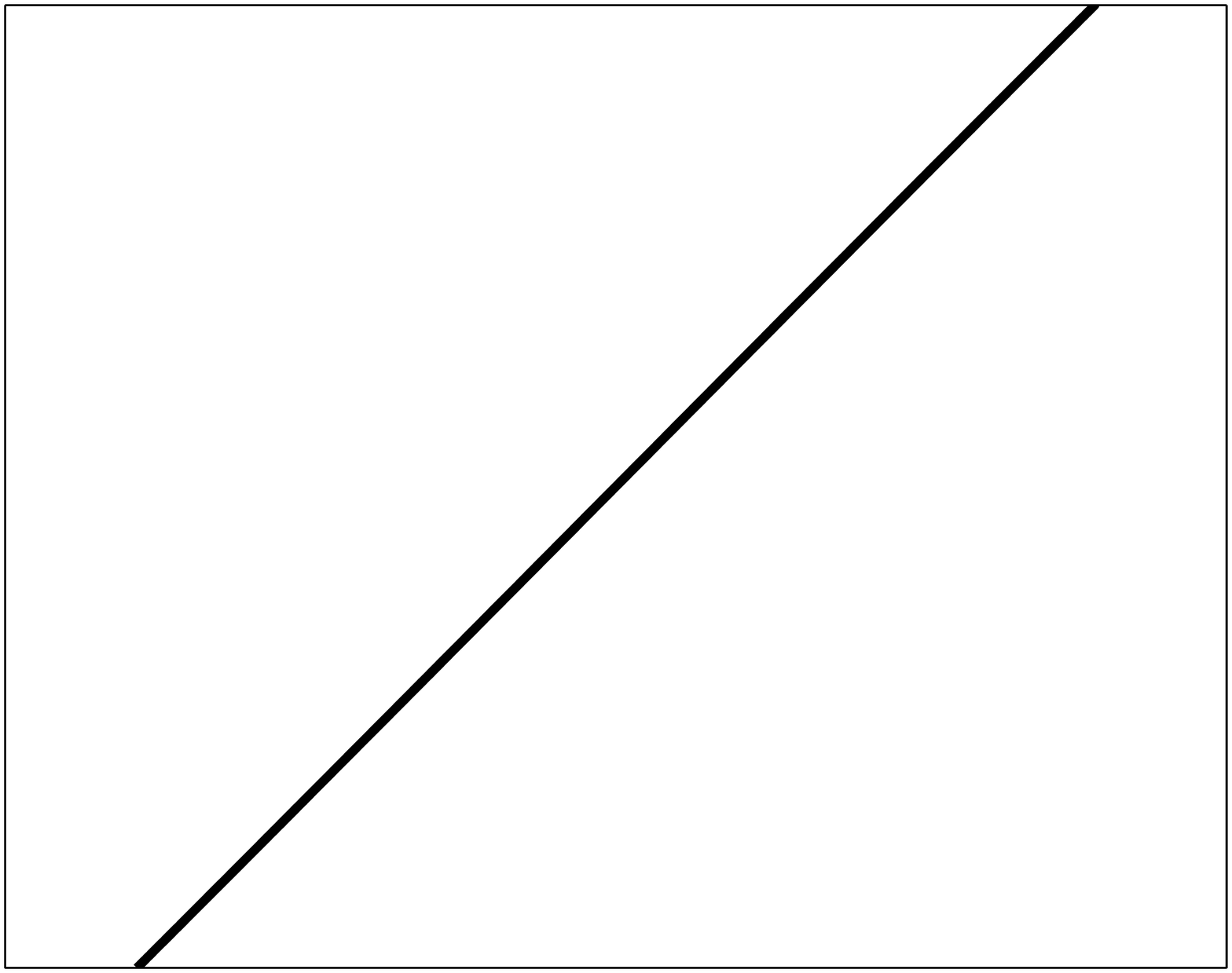}}
		\put(470,99.5){\includegraphics[width=11.2mm]{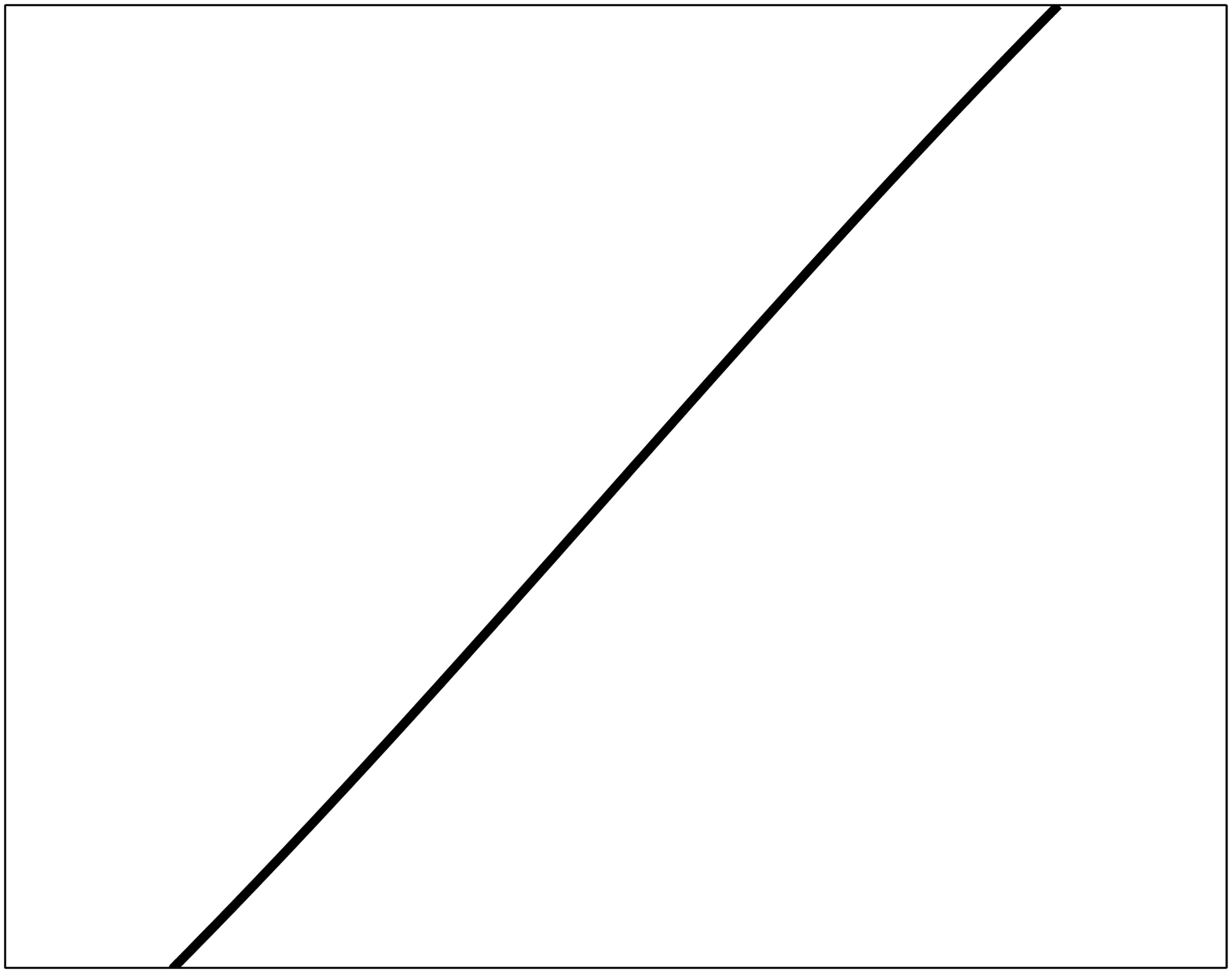}}
		\put(620,109){\includegraphics[width=11.2mm]{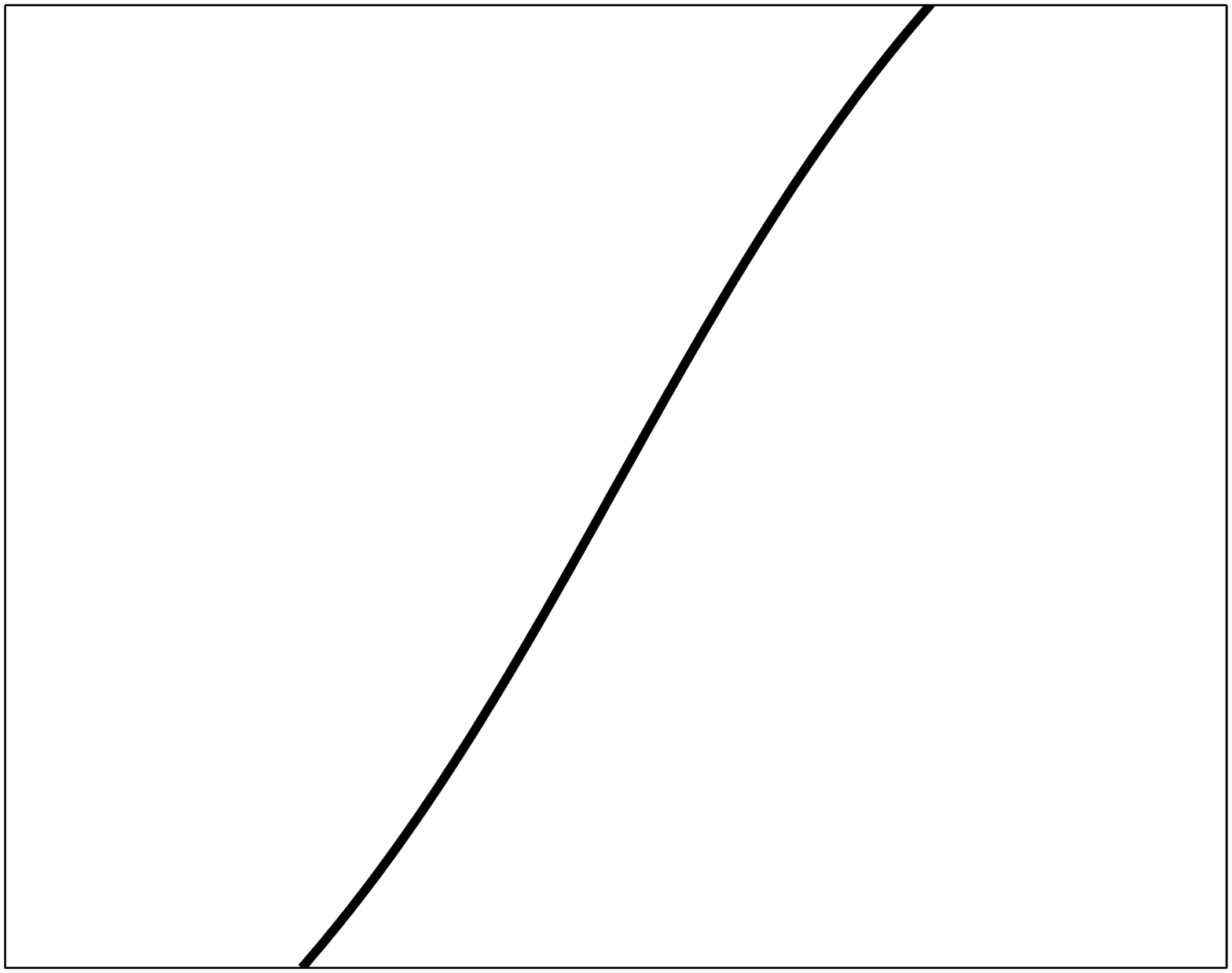}}
		\put(755,119.5){\includegraphics[width=11.2mm]{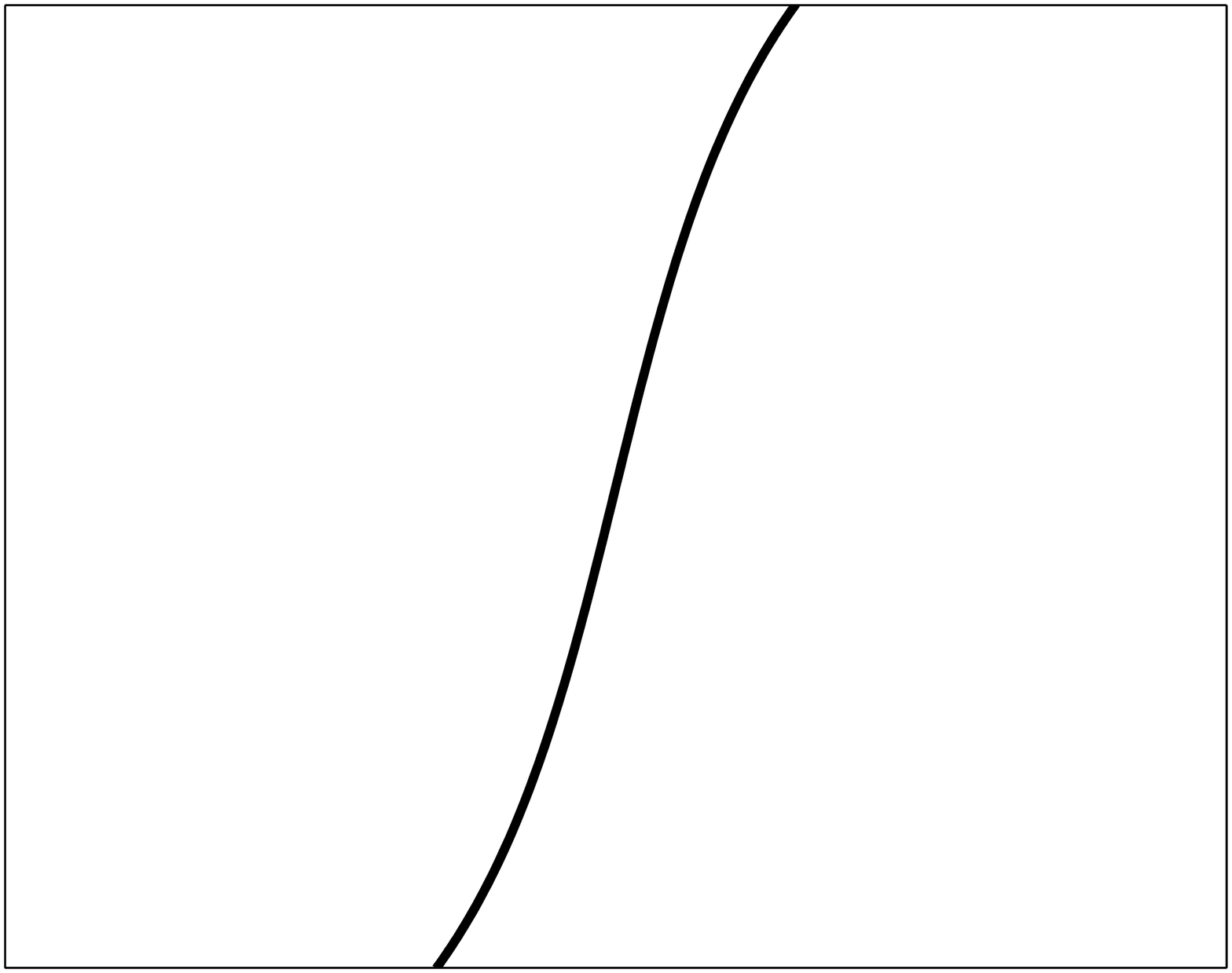}}
		\put(880,130){\includegraphics[width=11.2mm]{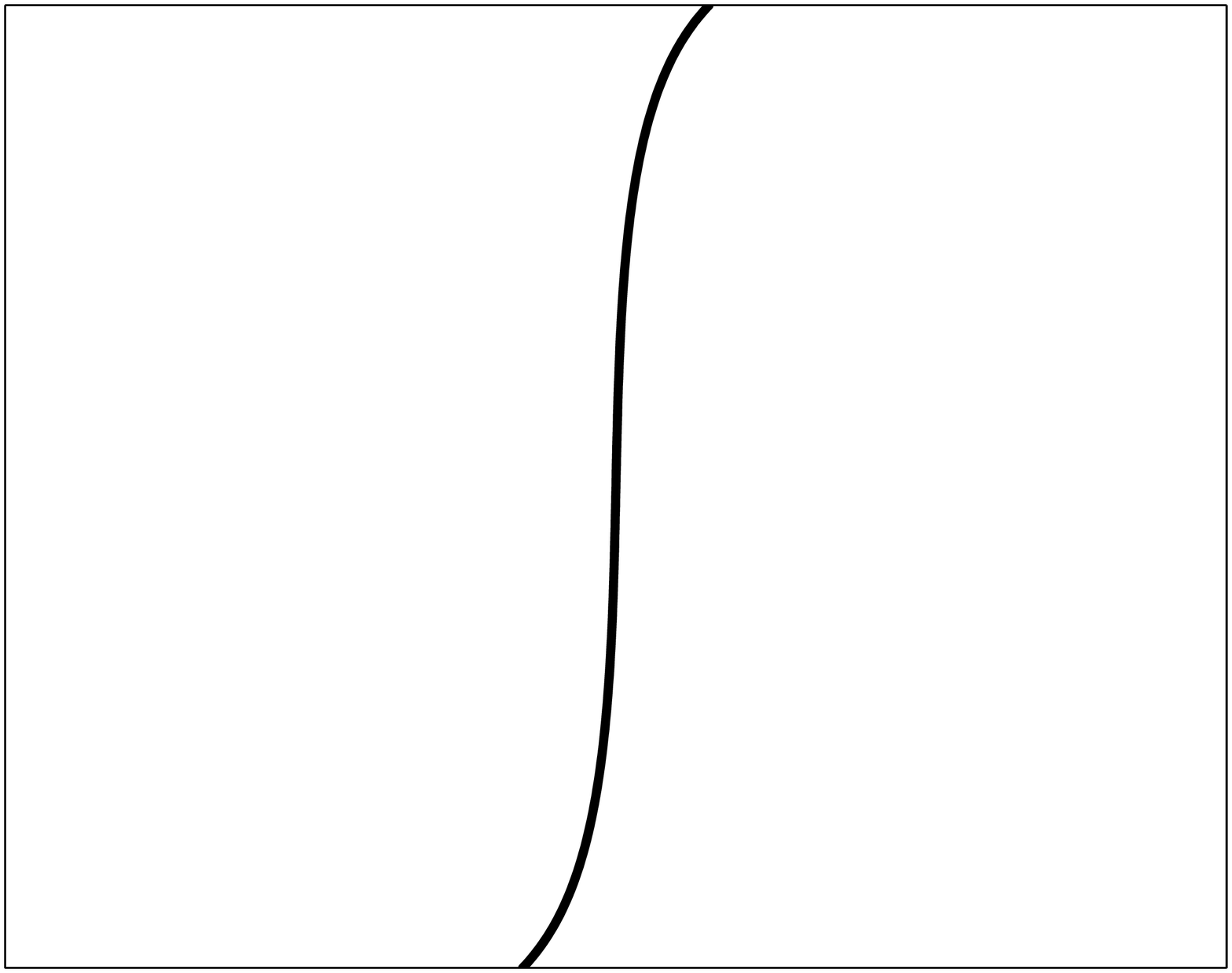}}
		\put(170,187.5){\pscircle(0,0){0.1}}
		\put(340,187.5){\pscircle(0,0){0.1}}
		\put(510,190){\pscircle(0,0){0.1}}
		\put(660,206){\pscircle(0,0){0.1}}
		\put(795,221){\pscircle(0,0){0.1}}
		\put(920,228.5){\pscircle(0,0){0.1}}
\end{overpic}
\caption{FEP of the two-DOF system described by Eqs.~(\ref{Eq:2DOF_Example}). NNM motions depicted in displacement space are inset. The horizontal and vertical axes in these plots are the displacements of the first and second DOF of the system, respectively.}
\label{Fig:FEP_2DOFsystem}
\end{center}
\end{figure}

Because a salient property of nonlinear systems is the frequency-energy dependence of their oscillations, the depiction of NNMs is conveniently realized in a frequency-energy plot (FEP). A NNM motion in a FEP is represented by a point associated with the fundamental frequency of the periodic motion, and with the total conserved energy accompanying the motion. A branch in a FEP details the complete frequency-energy dependence of the considered mode. Fig.~\ref{Fig:FEP_2DOFsystem} illustrates the FEP  of the two-DOF system described by the equations
\begin{equation}
\begin{array}{r c l}
    \ddot{q}_{1} + (2 \: q_{1} - q_{2}) + 0.5 \: q_{1}^{3} & = & 0 \\
    \ddot{q}_{2} + (2 \: q_{2} - q_{1}) & = & 0 .
\end{array}
\label{Eq:2DOF_Example}
\end{equation}
The plot features two branches corresponding to the in-phase and out-of-phase synchronous NNMs of the system. These fundamental NNMs are the direct nonlinear extension of the corresponding LNMs. The nonlinear modal parameters, \textit{i.e.} the frequencies of oscillation and the modal curves, are found to depend markedly on the energy. In particular, the frequency of the two fundamental NNMs increases with the energy level, revealing the hardening characteristic of the cubic stiffness nonlinearity in the system.

\begin{figure}[p]
\vspace*{-0.5cm}
\begin{center}
\scalebox{1} 
{
\begin{pspicture}(0,0)(15.0,23.7)
%
\definecolor{lightgrey}{rgb}{0.8,0.8,0.8}
%
%
\put(13.5,13.6){(a)}
%
\rput[bl](3.8,18.8){\psbezier[linecolor=black,linewidth=0.04,fillstyle=solid,fillcolor=lightgrey](1.90,3.80)(2.64,4.35)(5.74,3.76)(6.43,3.08)(7.12,2.40)(7.21,1.45)(5.90,0.97)(4.59,0.49)(3.78,0.12)(2.80,0.06)(1.82,0)(0.95,0.26)(0.48,1.08)(0,1.90)(0.37,2.64)(1.11,3.20)(1.85,3.75)(1.16,3.24)(1.90,3.80)}
\rput[bl](5,20.7){Structure}
\rput[bl](4.9,20.2){under test}
%
\rput[bl]{45}(8.3,20.675){\psline[linecolor=black,linewidth=0.04](0,0)(0.3,0)(0.42,-0.25)(0.54,0.25)(0.66,-0.25)(0.78,0.25)(0.90,-0.25)(1.02,0.25)(1.14,-0.25)(1.26,0)(1.56,0)}
\rput[bl]{45}(8.3,20.675){\psdots[dotsize=0.1](0,0)}
\rput[bl]{45}(8.3,20.675){\psdots[dotsize=0.1](1.56,0)}
\rput[bl]{45}(8.3,20.675){\psline[arrows=->,linewidth=0.02,arrowinset=0,arrowlength=1.2,arrowsize=3pt 2](0.5,-0.4)(1.05,0.6)}
%
\rput[bl](1.05,23.0){$\mathbf{p}(t)$}
\put(0.4,18.8){\includegraphics[width=2cm]{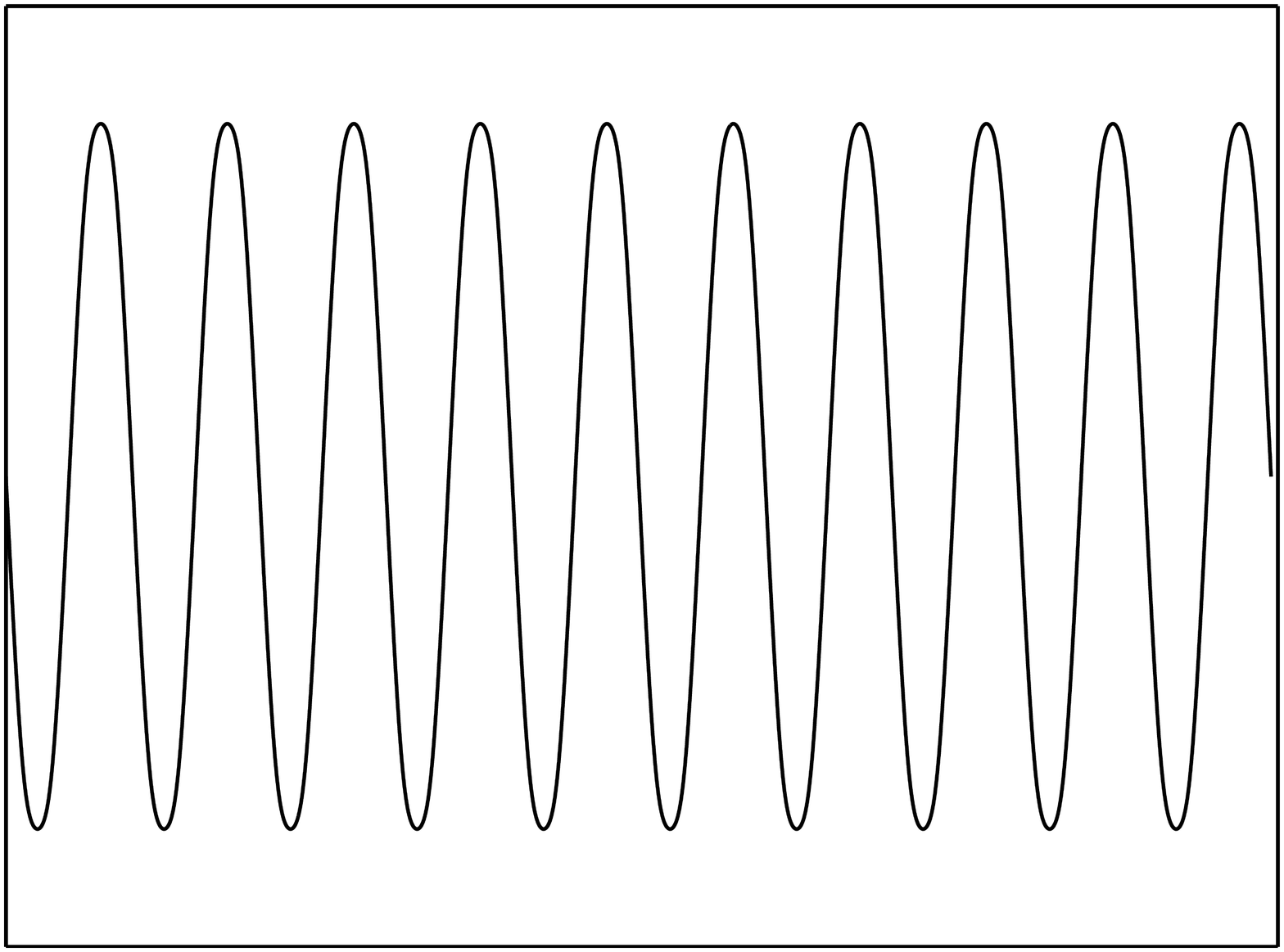}}
\psline[arrows=->,linewidth=0.02,arrowinset=0,arrowlength=1.2,arrowsize=3pt 2](2.6,19.575)(4.1,19.575)
\psdots[dotsize=0.1](3.1,19.575)
\psdots[dotsize=0.1](3.6,19.575)
\rput[bl](3.05,18.95){On}
\put(0.4,20.95){\includegraphics[width=2cm]{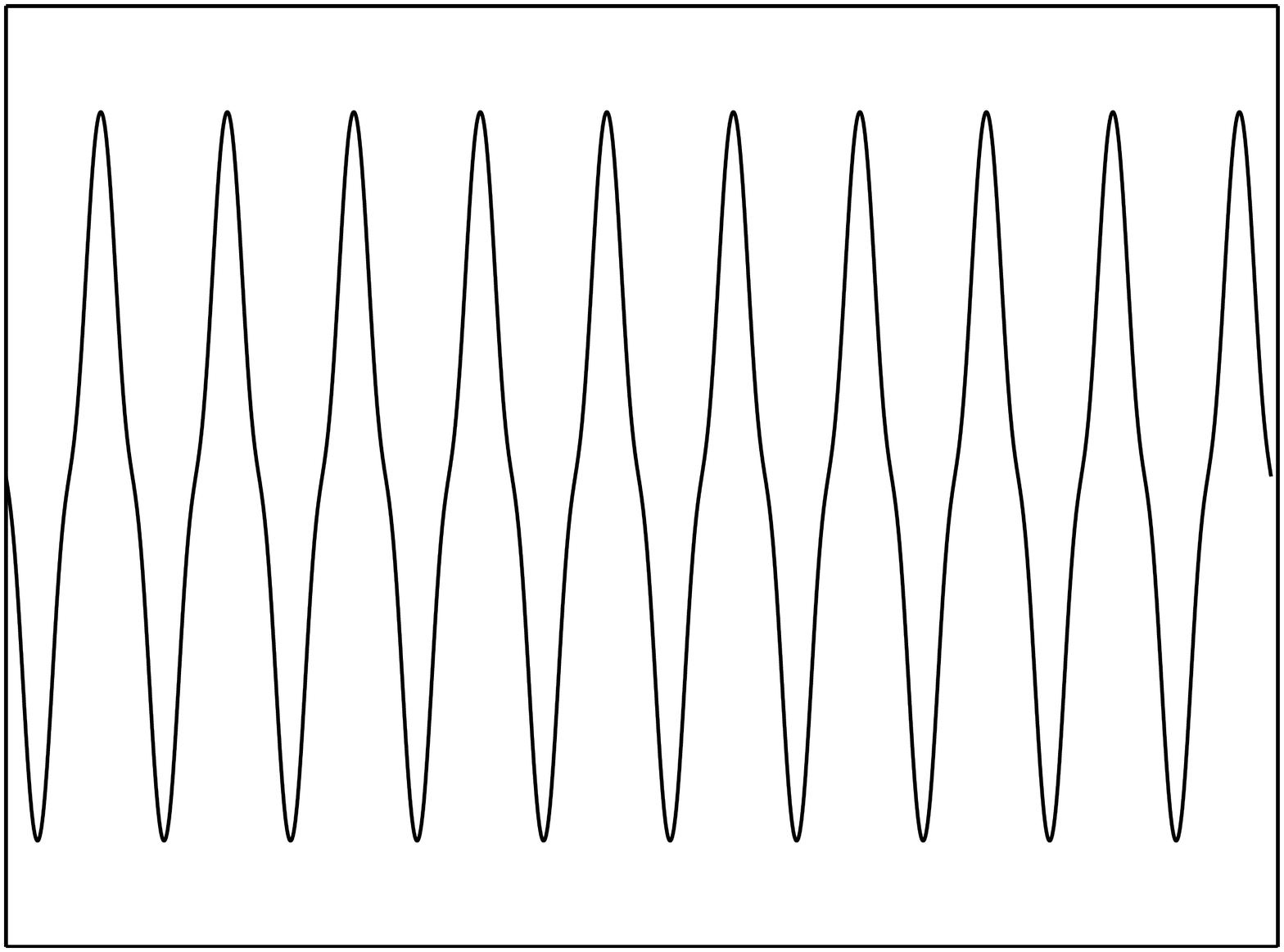}}
\psline[arrows=->,linewidth=0.02,arrowinset=0,arrowlength=1.2,arrowsize=3pt 2](2.6,21.725)(4.1,21.725)
\psdots[dotsize=0.1](3.1,21.725)
\psdots[dotsize=0.1](3.6,21.725)
\rput[bl](3.05,21.1){On}
\rput[bl](0.25,17.0){Stepped-sine}
\rput[bl](0.50,16.6){excitation}
%
\rput[bl](13.25,23.35){$\mathbf{q}(t)$}
\put(12.6,18.2){\includegraphics[width=2cm]{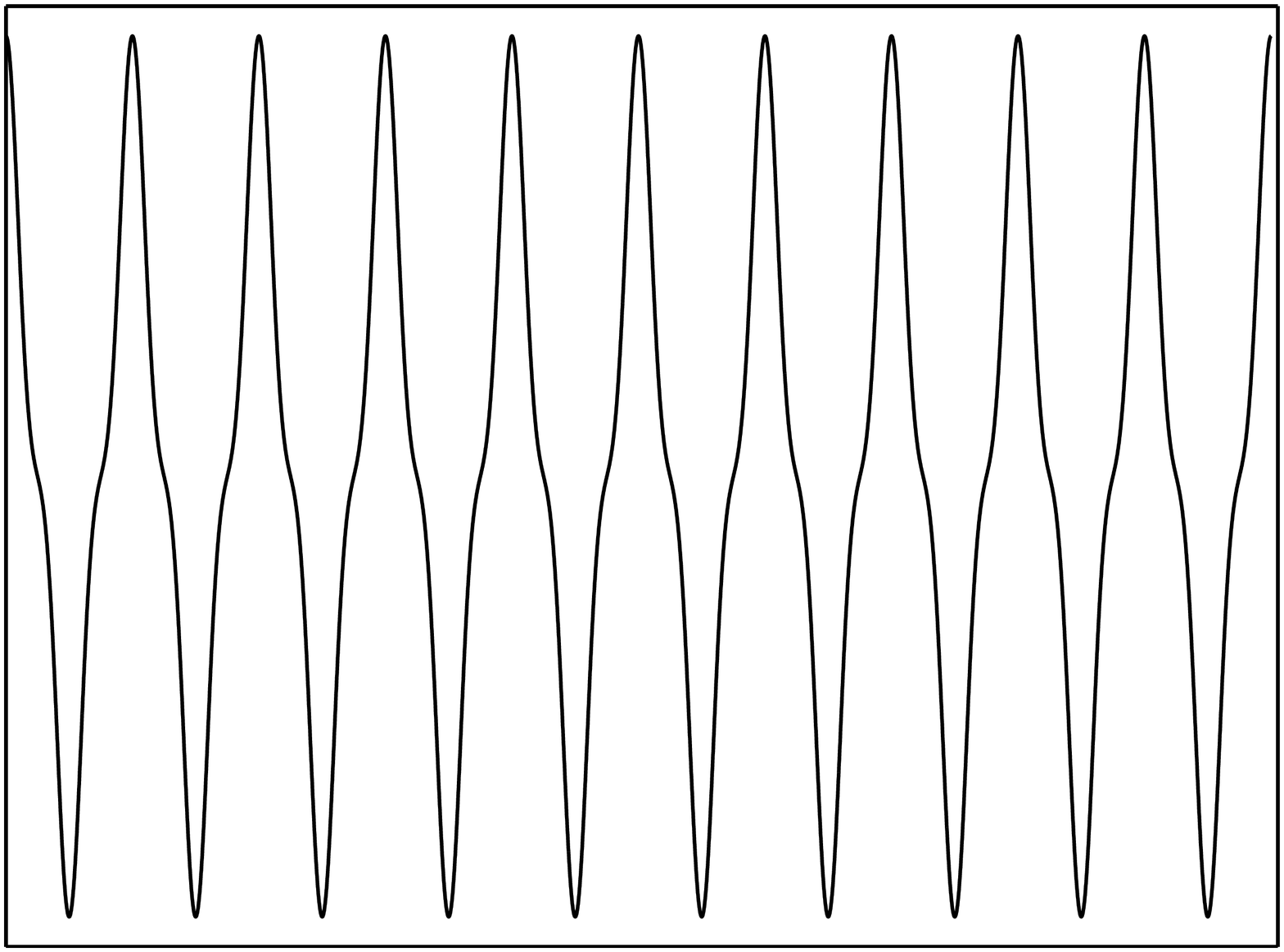}}
\psline[arrows=->,linewidth=0.02,arrowinset=0,arrowlength=1.2,arrowsize=3pt 2](10.9,19.55)(12.4,19.55)
\put(12.6,19.9){\includegraphics[width=2cm]{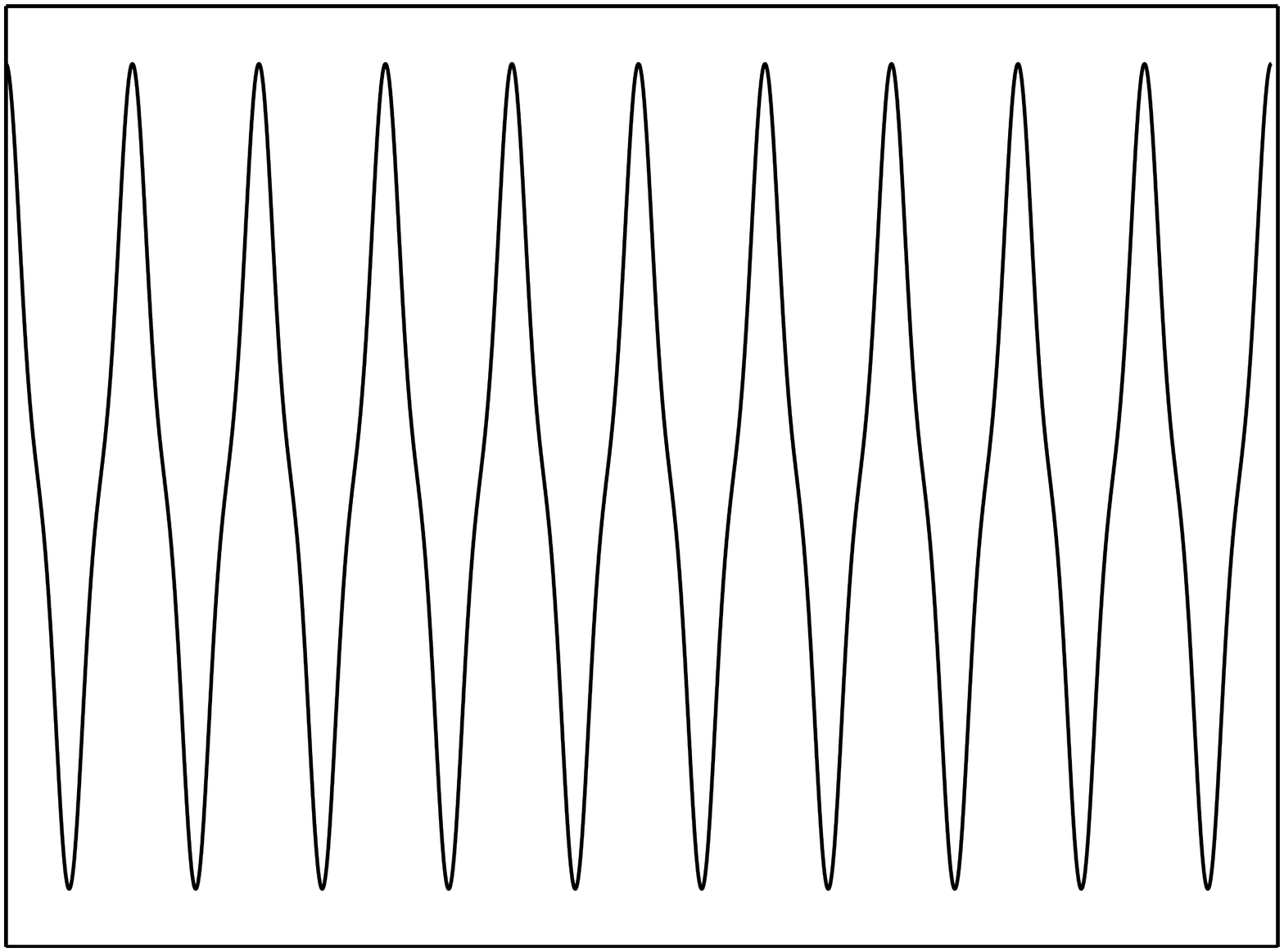}}
\psline[arrows=->,linewidth=0.02,arrowinset=0,arrowlength=1.2,arrowsize=3pt 2](10.9,20.675)(12.4,20.675)
\put(12.6,21.6){\includegraphics[width=2cm]{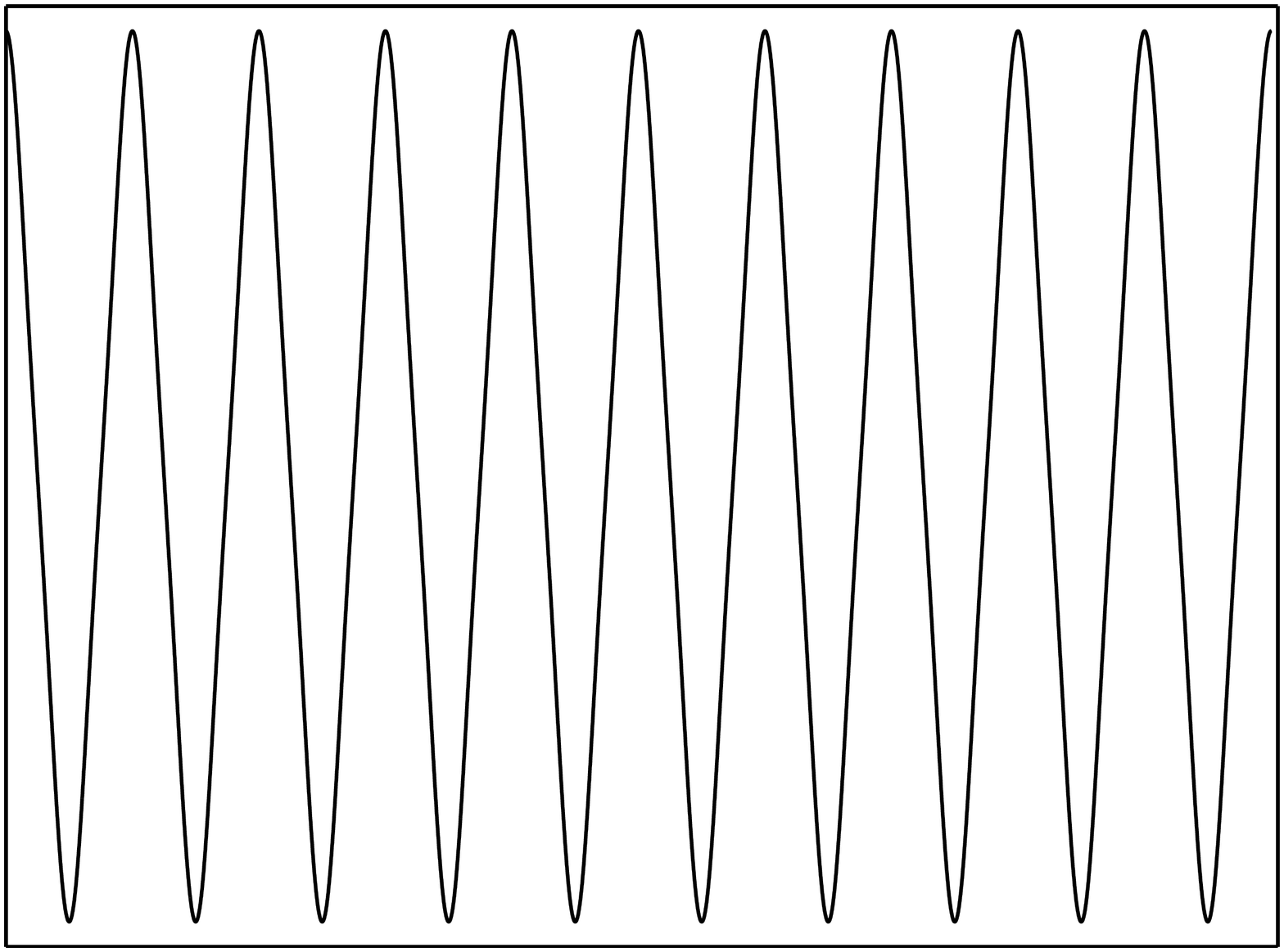}}
\psline[arrows=->,linewidth=0.02,arrowinset=0,arrowlength=1.2,arrowsize=3pt 2](10.9,21.8)(12.4,21.8)
\rput[bl](12.5,17.0){Steady-state}
\rput[bl](12.85,16.5){response}
%
\psline[arrows=->,linewidth=0.02,arrowinset=0,arrowlength=1.2,arrowsize=5pt 2](3,17.0)(5.3,17.0)
\psline[arrows=->,linewidth=0.02,arrowinset=0,arrowlength=1.2,arrowsize=5pt 2](12,17.0)(9.7,17.0)
\psframe[linewidth=0.02,dimen=outer](5.5,16.2)(9.5,17.8)
\rput[bl](6.65,17.0){Phase lag}
\rput[bl](6.6,16.6){estimation}
\psline[arrows=->,linewidth=0.02,arrowinset=0,arrowlength=1.2,arrowsize=5pt 2](7.5,16.0)(7.5,15.0)
%
\psline[linewidth=0.02](7.5,14.8)(7.5,12.65)
\psline[linewidth=0.02](7.5,14.8)(5.5,13.725)
\psline[linewidth=0.02](7.5,12.65)(5.5,13.725)
\rput[bl](6.3,13.6){90 $^{\circ}$ ?}
%
\psline[arrows=->,linewidth=0.02,arrowinset=0,arrowlength=1.2,arrowsize=5pt 2](5.3,13.725)(3,13.725)
\rput[bl](4,14.0){No}
%
\psframe[linewidth=0.02,dimen=outer](0,12.65)(2.8,14.8)
\rput[bl](0.5,14.1){Increment}
\rput[bl](0.5,13.6){excitation}
\rput[bl](0.55,13.0){frequency}
\psline[arrows=->,linewidth=0.02,arrowinset=0,arrowlength=1.2,arrowsize=5pt 2](1.4,15.0)(1.4,16.0)
%
\psline[arrows=->,linewidth=0.02,arrowinset=0,arrowlength=1.2,arrowsize=5pt 2](7.5,12.25)(7.5,10.5)
\rput[bl](8,11.5){Yes, see step 2 in (b)}
%
%
\psline[linewidth=0.02](0,11.2)(15,11.2)
\put(13.5,0.7){(b)}
%
\rput[bl](3.8,5.8){\psbezier[linecolor=black,linewidth=0.04,fillstyle=solid,fillcolor=lightgrey](1.90,3.80)(2.64,4.35)(5.74,3.76)(6.43,3.08)(7.12,2.40)(7.21,1.45)(5.90,0.97)(4.59,0.49)(3.78,0.12)(2.80,0.06)(1.82,0)(0.95,0.26)(0.48,1.08)(0,1.90)(0.37,2.64)(1.11,3.20)(1.85,3.75)(1.16,3.24)(1.90,3.80)}
\rput[bl](5,7.7){Structure}
\rput[bl](4.9,7.2){under test}
%
\rput[bl]{45}(8.3,7.675){\psline[linecolor=black,linewidth=0.04](0,0)(0.3,0)(0.42,-0.25)(0.54,0.25)(0.66,-0.25)(0.78,0.25)(0.90,-0.25)(1.02,0.25)(1.14,-0.25)(1.26,0)(1.56,0)}
\rput[bl]{45}(8.3,7.675){\psdots[dotsize=0.1](0,0)}
\rput[bl]{45}(8.3,7.675){\psdots[dotsize=0.1](1.56,0)}
\rput[bl]{45}(8.3,7.675){\psline[arrows=->,linewidth=0.02,arrowinset=0,arrowlength=1.2,arrowsize=3pt 2](0.5,-0.4)(1.05,0.6)}
%
\rput[bl](1.05,10.0){$\mathbf{p}(t)$}
\put(0.4,5.8){\includegraphics[width=2cm]{Illust_NPR_FA_Input1.eps}}
\psline[linewidth=0.02](2.6,6.575)(3.1,6.575)
\psdots[dotsize=0.1](3.1,6.575)
\psline[linewidth=0.02](3.2,7.05)(3.6,6.575)
\psdots[dotsize=0.1](3.2,7.05)
\psline[arrows=->,linewidth=0.02,arrowinset=0,arrowlength=1.2,arrowsize=3pt 2](3.6,6.575)(4.1,6.575)
\psdots[dotsize=0.1](3.6,6.575)
\rput[bl](3.05,5.95){Off}
\put(0.4,7.95){\includegraphics[width=2cm]{Illust_NPR_FA_Input2.eps}}
\psline[linewidth=0.02](2.6,8.725)(3.1,8.725)
\psdots[dotsize=0.1](3.1,8.725)
\psline[linewidth=0.02](3.2,9.2)(3.6,8.725)
\psdots[dotsize=0.1](3.2,9.2)
\psline[arrows=->,linewidth=0.02,arrowinset=0,arrowlength=1.2,arrowsize=3pt 2](3.6,8.725)(4.1,8.725)
\psdots[dotsize=0.1](3.6,8.725)
\rput[bl](3.05,8.1){Off}
\rput[bl](1.15,4.1){No}
\rput[bl](0.50,3.6){excitation}
%
\rput[bl](13.25,10.35){$\mathbf{q}(t)$}
\put(12.6,5.2){\includegraphics[width=2cm]{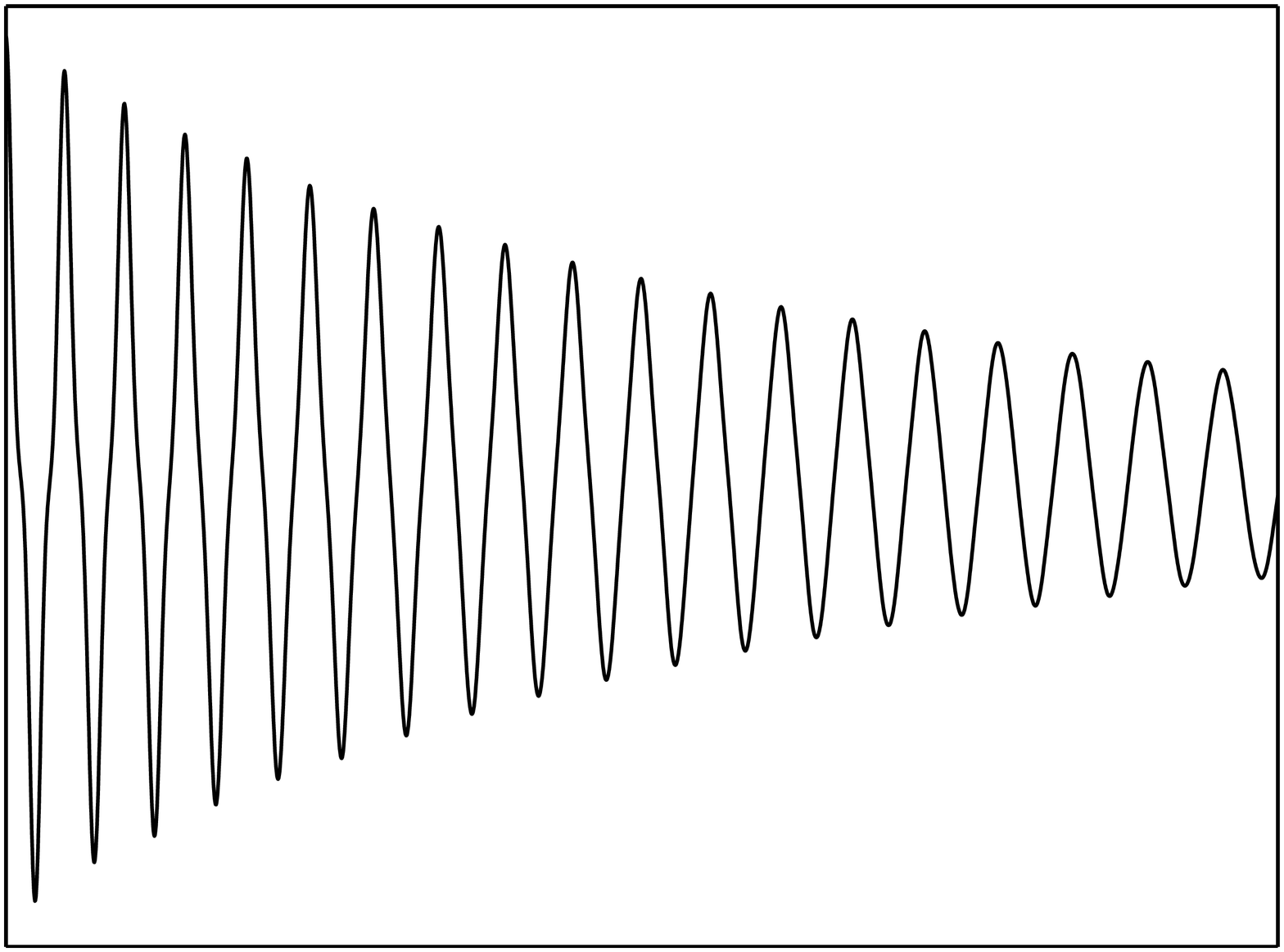}}
\psline[arrows=->,linewidth=0.02,arrowinset=0,arrowlength=1.2,arrowsize=3pt 2](10.9,6.55)(12.4,6.55)
\put(12.6,6.9){\includegraphics[width=2cm]{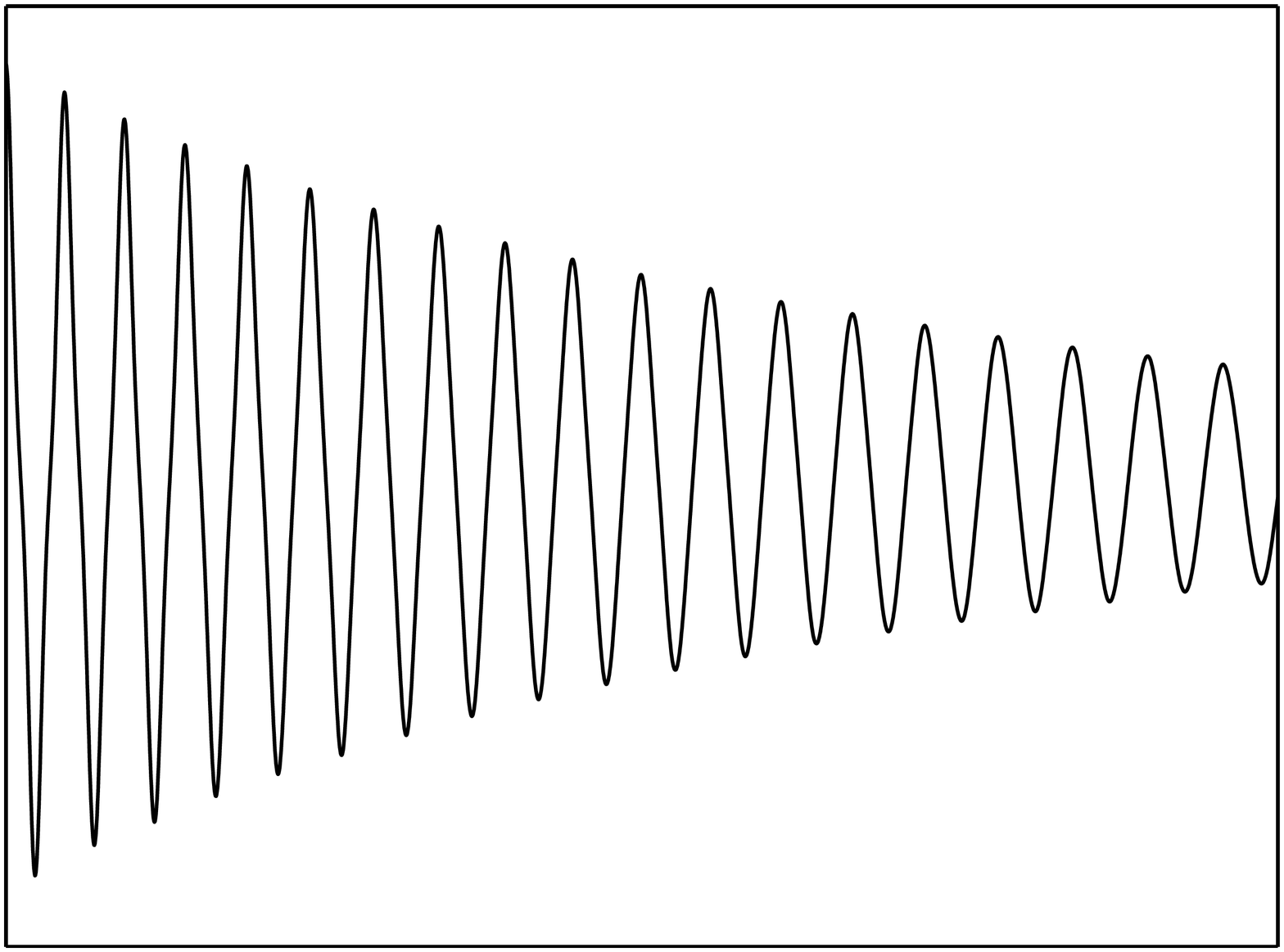}}
\psline[arrows=->,linewidth=0.02,arrowinset=0,arrowlength=1.2,arrowsize=3pt 2](10.9,7.675)(12.4,7.675)
\put(12.6,8.6){\includegraphics[width=2cm]{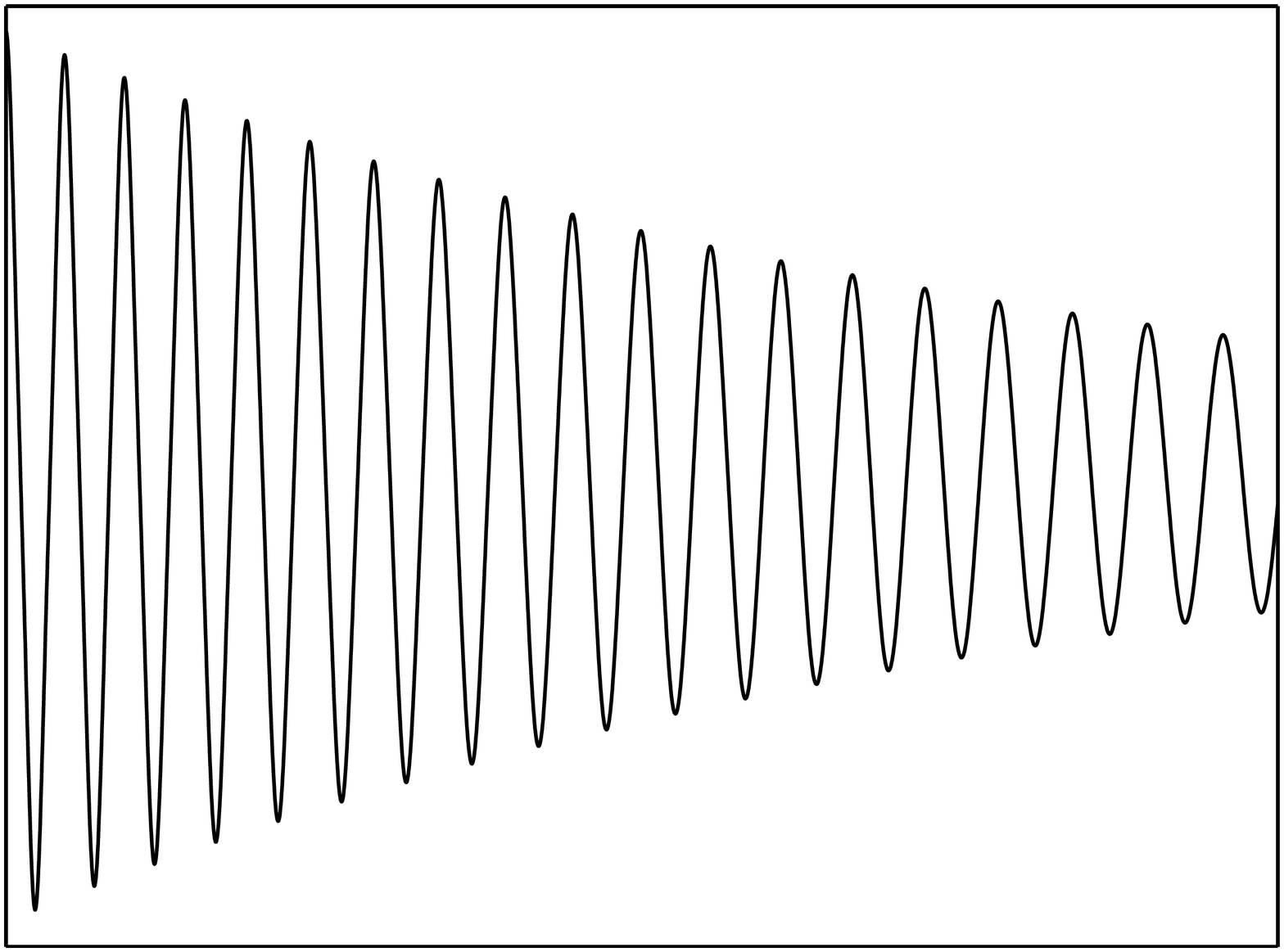}}
\psline[arrows=->,linewidth=0.02,arrowinset=0,arrowlength=1.2,arrowsize=3pt 2](10.9,8.8)(12.4,8.8)
\rput[bl](12.675,4.0){Free-decay}
\rput[bl](12.85,3.5){response}
%
\psline[arrows=->,linewidth=0.02,arrowinset=0,arrowlength=1.2,arrowsize=5pt 2](12,4.0)(9.7,4.0)
\psframe[linewidth=0.02,dimen=outer](5.5,3.2)(9.5,4.8)
\rput[bl](6.1,4.0){Time-frequency}
\rput[bl](6.8,3.5){analysis}
%
\psline[arrows=->,linewidth=0.02,arrowinset=0,arrowlength=1.2,arrowsize=5pt 2](7.5,3.0)(7.5,2.0)
\psframe[linewidth=0.02,dimen=outer](5.5,0.2)(9.5,1.8)
\rput[bl](6.625,1.1){Nonlinear}
\rput[bl](5.875,0.5){modal parameters}
\end{pspicture}
}
\caption{Experimental modal analysis of nonlinear systems using phase resonance~\cite{PRM_Part1}. (a) Step 1: NNM force appropriation; (b) step 2: NNM free-decay identification.}
\label{Fig:NPR_Methodology}
\end{center}
\end{figure}

Two essential properties of linear systems are preserved in the presence of nonlinearity. First, forced resonances of nonlinear systems occur in the neighborhood of NNMs~\cite{Vakakis_NNMBook}. Second, NNMs obey the invariance principle, which states that if the motion is initiated on one specific NNM, the remaining NNMs are quiescent for all time~\cite{ShawPierre_NNM}. These two properties were exploited in Ref.~\cite{PRM_Part1} to develop a nonlinear phase resonance method. The procedure comprises two steps, as illustrated in Fig.~\ref{Fig:NPR_Methodology}. During the first step, termed NNM force appropriation, the system is excited using a stepped-sine signal to induce a single-NNM motion at a prescribed energy level. This step is facilitated by a generalized phase lag quadrature criterion applicable to nonlinear systems~\cite{PRM_Part1}. This criterion asserts that a structure vibrates according to an underlying conservative NNM if the measured displacements possess, for all harmonics, a phase difference of ninety degrees with respect to the excitation. The second step of the procedure, termed NNM free-decay identification, turns off the excitation to track the energy dependence of the appropriated NNM. The associated modal parameters are extracted directly from the free damped system response through time-frequency analysis. This nonlinear phase resonance method was found to be highly accurate but, as in linear testing, very time-consuming. In addition, to reach the neighborhood of the resonance where a specific NNM lives may require a trial-and-error approach to deal with the shrinking basins of attraction along forced resonance peaks. The methodology described in the next section precisely addresses these two issues.

\newpage
\section{A two-step methodology for NNM identification under broadband forcing}\label{Sec:NPSM}

\begin{figure}[p]
\begin{center}
\scalebox{1} 
{
\begin{pspicture}(0,0)(15.0,19.2)
%
\definecolor{lightgrey}{rgb}{0.8,0.8,0.8}
%
\rput[bl](3.8,14.0){\psbezier[linecolor=black,linewidth=0.04,fillstyle=solid,fillcolor=lightgrey](1.90,3.80)(2.64,4.35)(5.74,3.76)(6.43,3.08)(7.12,2.40)(7.21,1.45)(5.90,0.97)(4.59,0.49)(3.78,0.12)(2.80,0.06)(1.82,0)(0.95,0.26)(0.48,1.08)(0,1.90)(0.37,2.64)(1.11,3.20)(1.85,3.75)(1.16,3.24)(1.90,3.80)}
\rput[bl](5,15.9){Structure}
\rput[bl](4.9,15.4){under test}
%
\rput[bl]{45}(8.3,15.875){\psline[linecolor=black,linewidth=0.04](0,0)(0.3,0)(0.42,-0.25)(0.54,0.25)(0.66,-0.25)(0.78,0.25)(0.90,-0.25)(1.02,0.25)(1.14,-0.25)(1.26,0)(1.56,0)}
\rput[bl]{45}(8.3,15.875){\psdots[dotsize=0.1](0,0)}
\rput[bl]{45}(8.3,15.875){\psdots[dotsize=0.1](1.56,0)}
\rput[bl]{45}(8.3,15.875){\psline[arrows=->,linewidth=0.02,arrowinset=0,arrowlength=1.2,arrowsize=3pt 2](0.5,-0.4)(1.05,0.6)}
%
\rput[bl](1.05,18.2){$\mathbf{p}(t)$}
\put(0.4,14.0){\includegraphics[width=2cm]{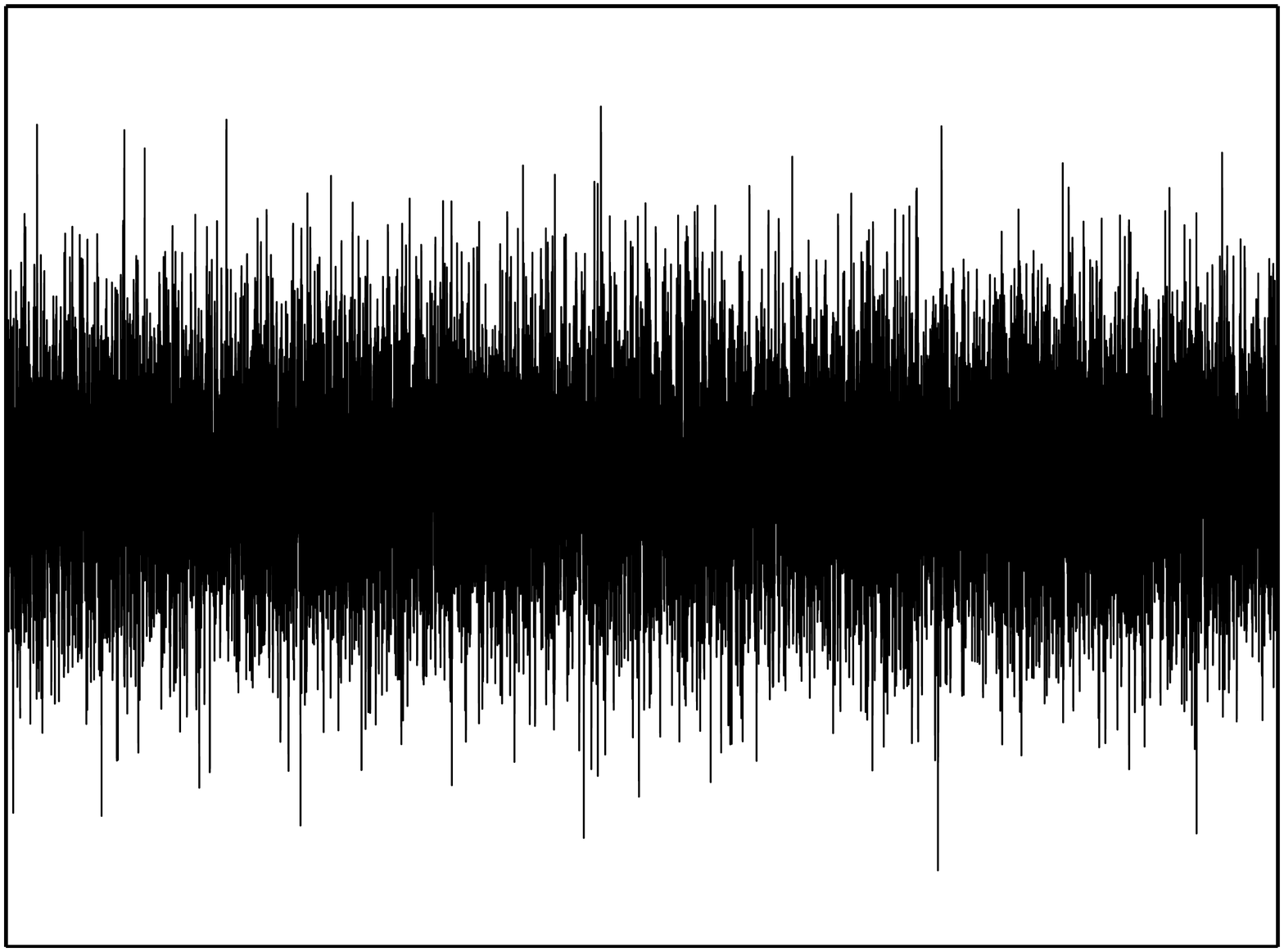}}
\psline[arrows=->,linewidth=0.02,arrowinset=0,arrowlength=1.2,arrowsize=3pt 2](2.6,14.775)(4.1,14.775)
\put(0.4,16.15){\includegraphics[width=2cm]{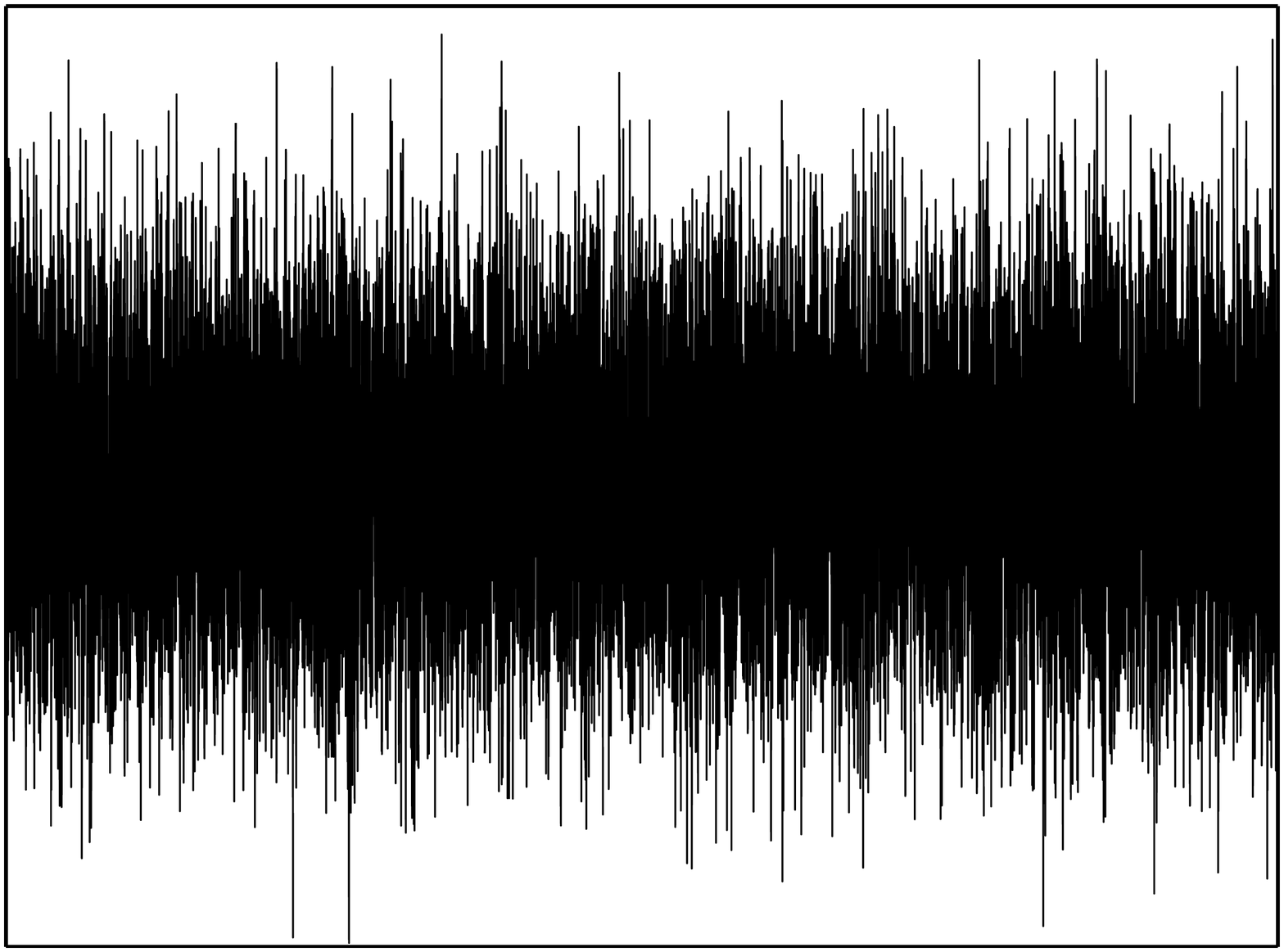}}
\psline[arrows=->,linewidth=0.02,arrowinset=0,arrowlength=1.2,arrowsize=3pt 2](2.6,16.925)(4.1,16.925)
\rput[bl](0.40,12.3){Broadband}
\rput[bl](0.50,11.8){excitation}
%
\rput[bl](13.25,18.55){$\mathbf{q}(t)$}
\put(12.6,13.4){\includegraphics[width=2cm]{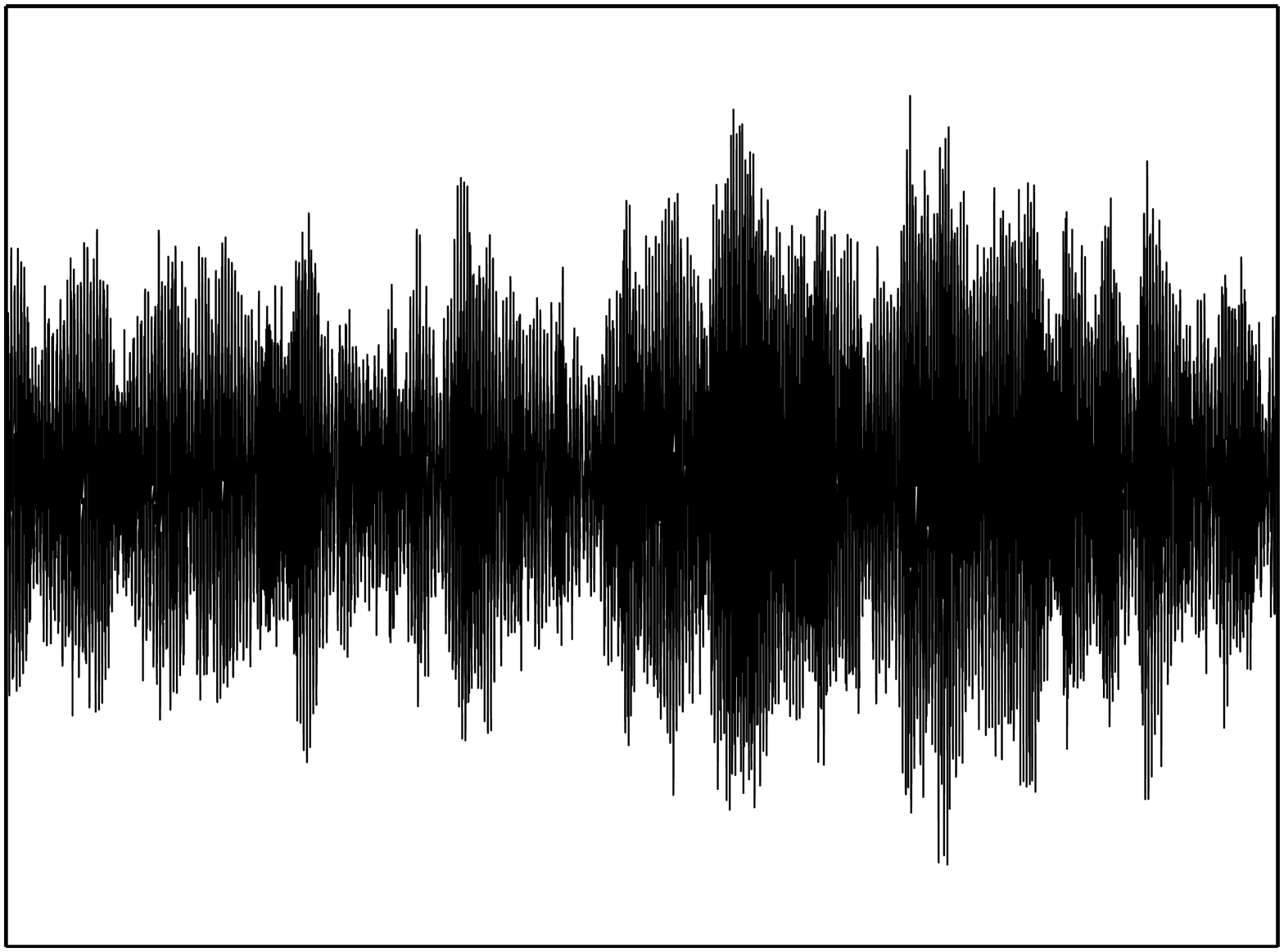}}
\psline[arrows=->,linewidth=0.02,arrowinset=0,arrowlength=1.2,arrowsize=3pt 2](10.9,14.75)(12.4,14.75)
\put(12.6,15.1){\includegraphics[width=2cm]{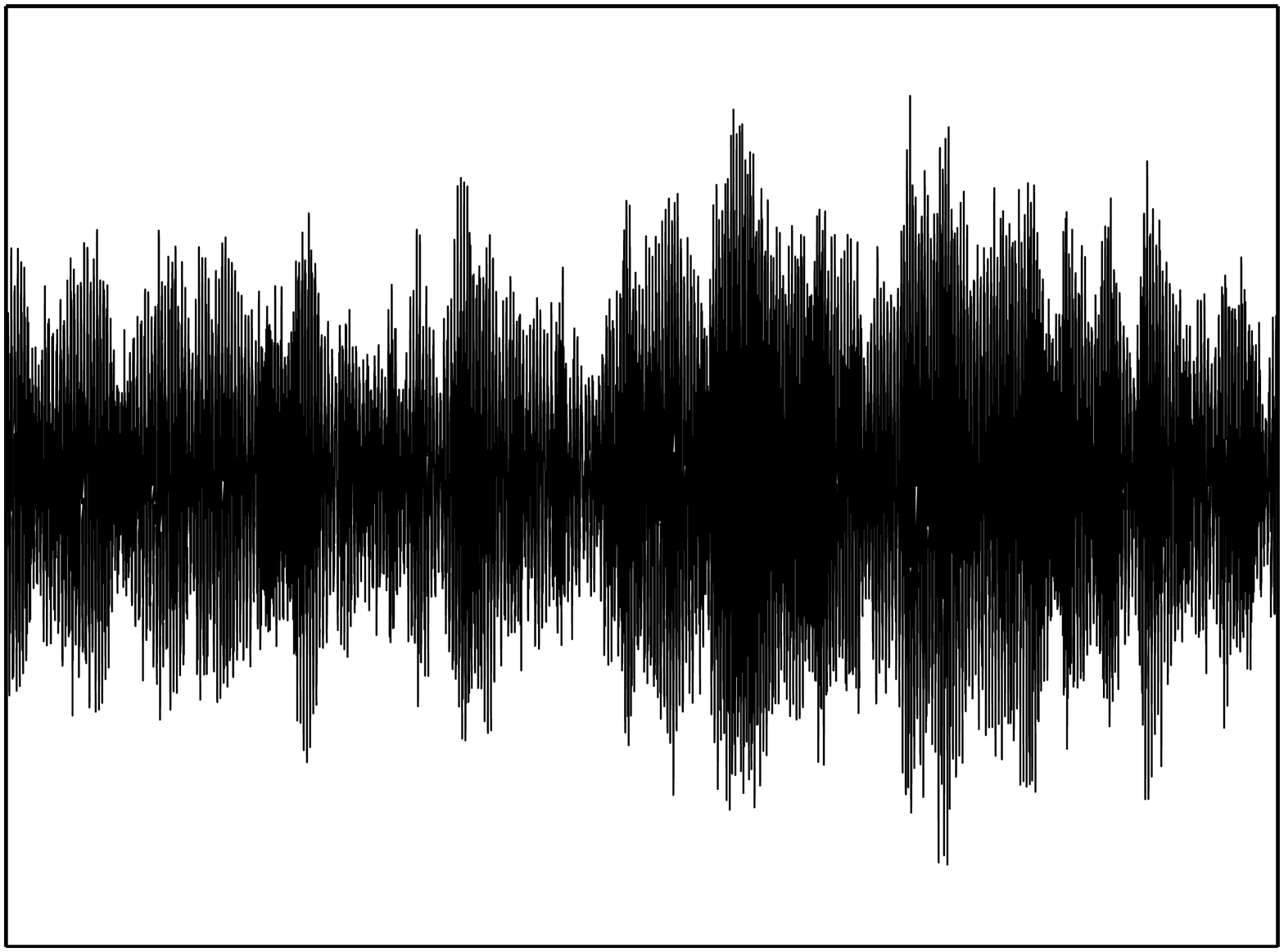}}
\psline[arrows=->,linewidth=0.02,arrowinset=0,arrowlength=1.2,arrowsize=3pt 2](8.3,15.875)(12.4,15.875)
\put(12.6,16.8){\includegraphics[width=2cm]{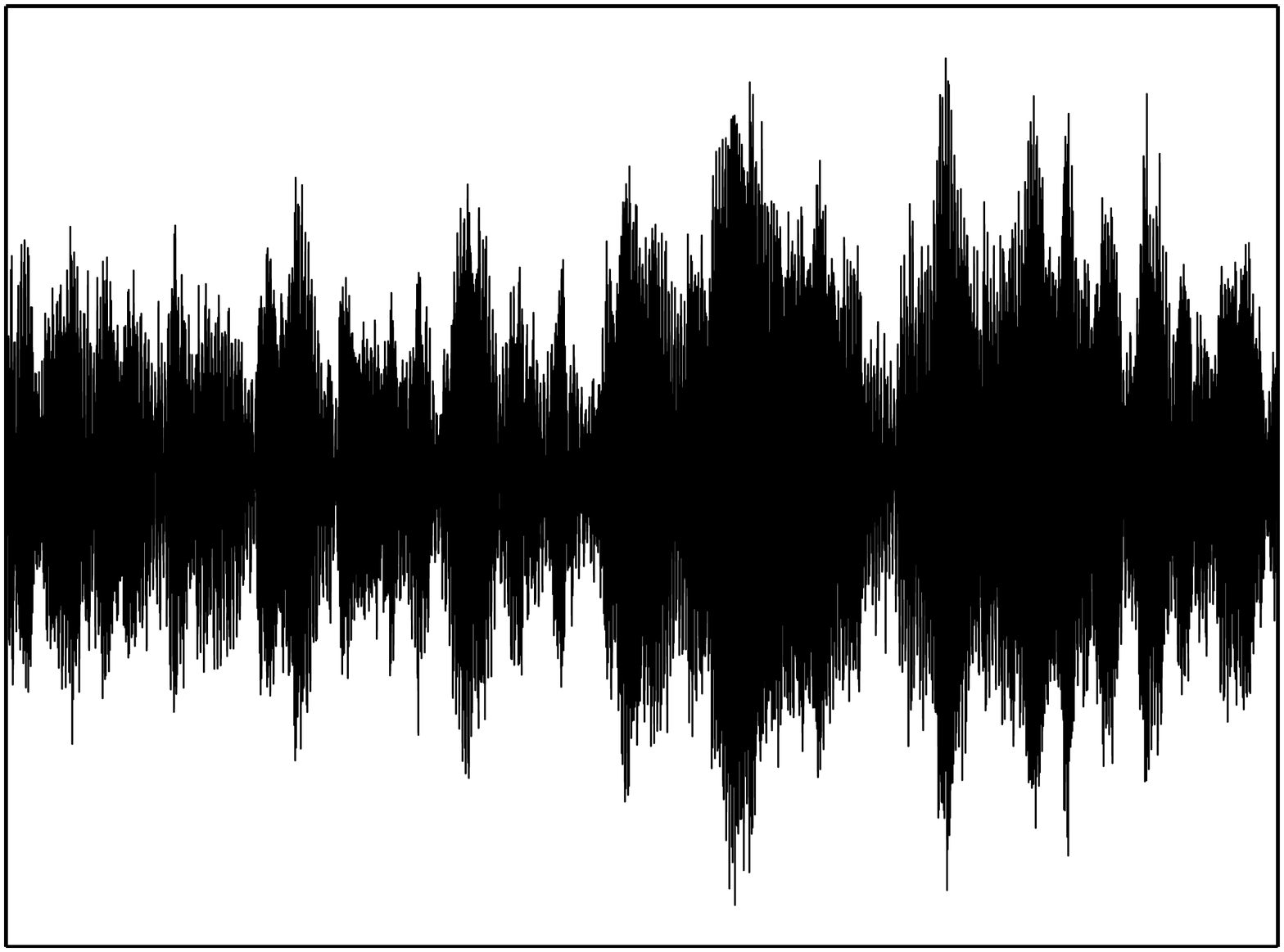}}
\psline[arrows=->,linewidth=0.02,arrowinset=0,arrowlength=1.2,arrowsize=3pt 2](9.4,17.0)(12.4,17.0)
\rput[bl](12.6,12.3){Broadband}
\rput[bl](12.85,11.7){response}
%
\psline[arrows=->,linewidth=0.02,arrowinset=0,arrowlength=1.2,arrowsize=5pt 2](3,12.2)(5.3,12.2)
\psline[arrows=->,linewidth=0.02,arrowinset=0,arrowlength=1.2,arrowsize=5pt 2](12,12.2)(9.7,12.2)
\psframe[linewidth=0.02,dimen=outer](5.5,11.4)(9.5,13.0)
\rput[bl](5.7,12.2){1. Nonlinear system}
\rput[bl](6.35,11.8){identification}
\psline[arrows=->,linewidth=0.02,arrowinset=0,arrowlength=1.2,arrowsize=5pt 2](7.5,11.2)(7.5,10.2)
%
\psframe[linewidth=0.02,dimen=outer](4.75,6.2)(10.25,10.0)
\rput[bl](5.25,9.1){Experimental undamped}
\rput[bl](6.3,8.7){modal model:}
\rput[bl](6.15,7.7){$\mathbf{\ddot{\overline{q}}}(t) + \mathbf{\overline{K}} \: \mathbf{\overline{q}}(t) + $}
\rput[bl](5.5,6.4){$\displaystyle \sum^{s}_{a=1} {c_{a} \: \mathbf{\Phi}^{T} \: \mathbf{h}_{a}(\mathbf{q}(t))} = 0$}
\psline[arrows=->,linewidth=0.02,arrowinset=0,arrowlength=1.2,arrowsize=5pt 2](7.5,6.0)(7.5,5.0)
%
\psframe[linewidth=0.02,dimen=outer](5.5,3.2)(9.5,4.8)
\rput[bl](6.325,4.1){2. Numerical}
\rput[bl](6.4,3.6){continuation}
\psline[arrows=->,linewidth=0.02,arrowinset=0,arrowlength=1.2,arrowsize=5pt 2](7.5,3.0)(7.5,2.0)
%
\psframe[linewidth=0.02,dimen=outer](5.5,0.2)(9.5,1.8)
\rput[bl](6.625,1.1){Nonlinear}
\rput[bl](5.875,0.5){modal parameters}
\end{pspicture}
}
\caption{Proposed methodology for the identification of NNMs based on broadband measurements. It comprises two major steps, namely nonlinear system identification and numerical continuation.}
\label{Fig:NPS_Methodology}
\end{center}
\end{figure}

The proposed methodology, presented in Fig.~\ref{Fig:NPS_Methodology}, comprises two major steps. The first step, described in Section~\ref{Sec:FNSI}, processes acquired input and output data using the FNSI method to derive an experimental state-space model of the structure. The second step, described in Section~\ref{Sec:NNMCont}, converts this state-space model into a model in modal space from which the energy-dependent frequencies and modal curves of the excited NNMs are computed individually using shooting and pseudo-arclength continuation.

\subsection{Identification of a nonlinear state-space model}\label{Sec:FNSI}

The FNSI method is capable of deriving models of nonlinear vibrating systems directly from measured data, and without resorting to a preexisting numerical model, \textit{e.g.}, a finite element model~\cite{Noel_FNSI}. It is applicable to multi-input, multi-output structures with high damping and high modal density, and makes no assumption as to the importance of nonlinearity in the measured dynamics~\cite{Noel_FNSI_SmallSat,Noel_SolarPanels}. 

\subsubsection{Feedback interpretation and state-space model identification}\label{Sec:Feedback}

The vibrations of damped nonlinear systems obey Newton's second law of dynamics
\begin{equation}
\mathbf{M} \: \mathbf{\ddot{q}}(t) + \mathbf{C}_{v} \: \mathbf{\dot{q}}(t) + \mathbf{K} \: \mathbf{q}(t) + \mathbf{f}(\mathbf{q}(t)) = \mathbf{p}(t)
\label{Eq:TD_Model}
\end{equation}
where $\mathbf{C}_{v} \in \mathbb{R}^{\: n_{p} \times n_{p}}$ is the linear viscous damping matrix; $\mathbf{p}(t) \in \mathbb{R}^{\: n_{p}}$ is the generalized external force vector; $\mathbf{f}(\mathbf{q}(t)) \in \mathbb{R}^{\: n_{p}}$ is the nonlinear restoring force vector, encompassing elastic terms only as this paper concentrates on stiffness nonlinearities. Note that Eq.~(\ref{Eq:TD_Model}) represents the damped and forced generalization of Eq.~(\ref{Eq:TD_Model_Undamped}). The nonlinear restoring force term in Eq.~(\ref{Eq:TD_Model}) is expressed by means of a linear combination of basis functions $\mathbf{h}_{a}(\mathbf{q}(t))$ as
\begin{equation}
\mathbf{f}(\mathbf{q}(t)) = \displaystyle \sum^{s}_{a=1} {c_{a} \: \mathbf{h}_{a}(\mathbf{q}(t))} .
\label{Eq:NL_Def}
\end{equation}
Given measurements of $\mathbf{p}(t)$ and $\mathbf{q}(t)$ or its derivatives, and an appropriate selection of the functionals $\mathbf{h}_{a}(\mathbf{q}(t))$, the objective of the FNSI method is to identify a state-space model from which the nonlinear coefficients $c_{a}$ can be estimated. The nonlinear components in the structure must therefore be instrumented on both sides in order to measure the relative displacement required in the formulation of $\mathbf{h}_{a}(\mathbf{q}(t))$, as illustrated in Fig.~\ref{Fig:NPS_Methodology}.

The FNSI approach builds on a block-oriented interpretation of nonlinear structural dynamics, which sees nonlinearities as a feedback into the linear system in the forward loop~\cite{Adams_NIFO2}. This interpretation boils down to moving the nonlinear internal forces in Eq.~(\ref{Eq:TD_Model}) to the right-hand side, and viewing them as additional external forces applied to the underlying linear structure, that is,
\begin{equation}
\mathbf{M} \: \mathbf{\ddot{q}}(t) + \mathbf{C}_{v} \: \mathbf{\dot{q}}(t) + \mathbf{K} \: \mathbf{q}(t) = \mathbf{p}(t) - \displaystyle \sum^{s}_{a=1} {c_{a} \: \mathbf{h}_{a}(\mathbf{q}(t))} .
\label{Eq:Feedback}
\end{equation}

Without loss of generality, it is assumed in this section that the structural response is measured in terms of generalized displacements. Defining the state vector $\mathbf{x} = \left( \mathbf{q}^{T} \ \ \mathbf{\dot{q}}^{T} \right)^{T}~\in~\mathbb{R}^{\: n_{s}}$, where $n_{s} = 2 \: n_{p}$ is the dimension of the state space and $T$ the matrix transpose, Eq.~(\ref{Eq:Feedback}) is recast into the set of first-order equations
\begin{equation}
\left\lbrace
\begin{array}{r c l}
    \mathbf{\dot{x}}(t) & = & \mathbf{A} \: \mathbf{x}(t) + \mathbf{B} \: \mathbf{e}(t) \\
    \mathbf{q}(t) & = & \mathbf{C} \: \mathbf{x}(t) + \mathbf{D} \: \mathbf{e}(t) ,
\end{array} \right.
\label{Eq:StateSpace}
\end{equation}
where the vector $\mathbf{e} \in \mathbb{R}^{\: (s+1) \: n_{p}}$, termed \textit{extended input vector}, concatenates the external forces $\mathbf{p}(t)$ and the nonlinear basis functions $\mathbf{h}_{a}(t)$. The matrices $\mathbf{A} \in \mathbb{R}^{\: n_{s} \times n_{s}}$, $\mathbf{B} \in \mathbb{R}^{\: n_{s} \times (s+1) \: n_{p}}$, $\mathbf{C} \in \mathbb{R}^{\: n_{p} \times n_{s}}$ and $\mathbf{D} \in \mathbb{R}^{\: n_{p} \times (s+1) \: n_{p}}$ are the state, extended input, output and direct feedthrough matrices, respectively. State-space and physical-space matrices correspond through the relations
$$\begin{array}{c c}
\mathbf{A} = \left( \begin{array}{c c}
      \mathbf{0}^{\: n_{p} \times n_{p}} & \mathbf{I}^{\: n_{p} \times n_{p}} \\
      -\mathbf{M}^{-1} \: \mathbf{K} & -\mathbf{M}^{-1} \: \mathbf{C}_{v} \\
        \end{array} \right) \: ; & \mathbf{B} = \left( \begin{array}{c c c c c}
                                                    \mathbf{0}^{\: n_{p} \times n_{p}} & \mathbf{0}^{\: n_{p} \times n_{p}} & \mathbf{0}^{\: n_{p} \times n_{p}} & \ldots & \mathbf{0}^{\: n_{p} \times n_{p}} \\
                                                    \mathbf{M}^{-1} & -c_{1} \: \mathbf{M}^{-1} & -c_{2} \: \mathbf{M}^{-1} & \ldots & -c_{s} \: \mathbf{M}^{-1} \\
                                                    \end{array} \right) \\
\end{array}$$
\begin{equation}
\begin{array}{c c}
\mathbf{C} = \left( \begin{array}{c c}
             \mathbf{I}^{\: n_{p} \times n_{p}} & \mathbf{0}^{\: n_{p} \times n_{p}}\\
             \end{array} \right) \: ; & \mathbf{D} = \mathbf{0}^{\: n_{p} \times (s+1) \: n_{p}} , \\
\end{array}
\label{Eq:PSpace2SSpace}
\end{equation}
where $\mathbf{0}$ and $\mathbf{I}$ are the zero and identity matrices, respectively. 

The FNSI estimation of the four system matrices $\mathbf{A}$, $\mathbf{B}$, $\mathbf{C}$ and $\mathbf{D}$ is achieved in the frequency domain using a classical subspace resolution scheme. This resolution essentially involves the reformulation of Eqs.~(\ref{Eq:StateSpace}) in matrix form, and the computation of estimates of $\mathbf{A}$, $\mathbf{B}$, $\mathbf{C}$ and $\mathbf{D}$ through geometrical manipulations of input and output data. The interested reader is referred to Ref.~\cite{Noel_FNSI} for a complete introduction to the theoretical and practical aspects of the FNSI method.

\subsubsection{Estimation of the nonlinear coefficients}\label{Sec:SS2PMS}

It is well-known that the matrices $\mathbf{A}$, $\mathbf{B}$, $\mathbf{C}$ and $\mathbf{D}$ defined in Eqs.~(\ref{Eq:PSpace2SSpace}) can be retrieved using subspace identification only up to an unknown similarity transformation of the state-space basis~\cite{VODM_Book}. This implies that estimates of the nonlinear coefficients $c_{a}$ cannot be obtained from a direct inspection of matrix $\mathbf{B}$. This issue was resolved in Ref.~\cite{Marchesiello_TNSI} by forming the transfer function matrix of the state-space model
\begin{equation}
\mathbf{G}_{s}(\omega) = \mathbf{C} \left( j \omega \: \mathbf{I}^{\: n_{s} \times n_{s}} - \mathbf{A} \right)^{-1} \mathbf{B} + \mathbf{D} ,
\label{Eq:G_StateSpace_CT}
\end{equation}
where $\omega$ is the pulsation and $j$ the imaginary unit. 

Matrix $\mathbf{G}_{s}(\omega)$ is invariant with respect to any similarity transformation, and relates the extended input vector to the measured response. Indeed, substituting Eq.~(\ref{Eq:NL_Def}) into Eq.~(\ref{Eq:TD_Model}) and moving to the frequency domain yields
\begin{equation}
\mathbf{G}^{-1}(\omega) \mathbf{Q}(\omega) + \displaystyle \sum^{s}_{a=1} {c_{a} \: \mathbf{H}_{a}(\omega)} = \mathbf{P}(\omega) ,
\label{Eq:FD_Model}
\end{equation}	
where $\mathbf{G}(\omega) = \left( -\omega^{2}\:\mathbf{M} + j\:\omega\:\mathbf{C}_{v} + \mathbf{K} \right)^{-1}$ is the transfer function matrix of the underlying linear system, and where $\mathbf{Q}(\omega)$, $\mathbf{H}_{a}(\omega)$ and $\mathbf{P}(\omega)$ are the continuous Fourier transforms of $\mathbf{q}(t)$, $\mathbf{h}_{a}(t)$ and $\mathbf{p}(t)$, respectively. The concatenation of $\mathbf{P}(\omega)$ and $\mathbf{H}_{a}(\omega)$ into the extended input spectrum $\mathbf{E}(\omega)$ finally gives
\begin{equation}
\mathbf{Q}(\omega) = \mathbf{G}(\omega) \left[ \begin{array}{c c c c}
\mathbf{I}^{\: n_{p} \times n_{p}} & -c_{1} \: \mathbf{I}^{\: n_{p} \times n_{p}} & \ldots & -c_{s} \: \mathbf{I}^{\: n_{p} \times n_{p}} \\
																\end{array} \right] \: \mathbf{E}(\omega) = \mathbf{G}_{s}(\omega) \: \mathbf{E}(\omega) .
\label{Eq:ExtendedFRF}
\end{equation}
The nonlinear coefficients $c_{a}$, together with the frequency response functions (FRFs) in $\mathbf{G}(\omega)$, can be directly extracted from Eq.~(\ref{Eq:ExtendedFRF}), given the transfer function matrix $\mathbf{G}_{s}(\omega)$ estimated from Eq.~(\ref{Eq:G_StateSpace_CT}).

\subsection{Computation of NNMs in modal space}\label{Sec:NNMCont}

The calculation of NNMs is not realized in this work in state space but in modal space in order to be compatible with the computational framework of Ref~\cite{NNM_Part2}. For that purpose, Eq.~(\ref{Eq:TD_Model_Undamped}) is recast into
\begin{equation}
\mathbf{\ddot{\overline{q}}}(t) + \mathbf{\overline{K}} \: \mathbf{\overline{q}}(t) + \displaystyle \sum^{s}_{a=1} {c_{a} \: \mathbf{\Phi}^{T} \: \mathbf{h}_{a}(\mathbf{q}(t))} = 0 ,
\label{Eq:TD_ModalModel_Undamped}
\end{equation}
where $\mathbf{\overline{q}}(t)=\mathbf{\Phi}^{-1}\mathbf{q}(t)$ is the vector of linear modal coordinates. Matrix $\mathbf{\Phi}$ contains the mode shapes $\boldsymbol\phi^{(i)}$ scaled to unit modal mass, \textit{i.e.} $\mathbf{\overline{M}} = \mathbf{\Phi}^{T} \: \mathbf{M} \: \mathbf{\Phi} = \mathbf{I}^{\: n_{p} \times n_{p}}$ and $ \mathbf{\overline{K}} = \mathbf{\Phi}^{T} \: \mathbf{K} \: \mathbf{\Phi} = \text{diag}\left( \omega^{2}_{i,0} \right)$. 

To simulate Eq.~(\ref{Eq:TD_ModalModel_Undamped}), the knowledge of the undamped frequencies $\omega_{i,0}$, nonlinear coefficients $c_{a}$, scaled mode shapes $\boldsymbol\phi^{(i)}$ and basis functions $\mathbf{h}_{a}$ is required. The nonlinear coefficients and basis functions are known from Section~\ref{Sec:FNSI}. The undamped frequencies $\omega_{i,0}$ are the absolute values of the complex eigenvalues $\lambda_{i}$ of matrix $\mathbf{A}$ estimated using FNSI in the previous section, that is, 
\begin{equation}
\mathbf{A} \: \boldsymbol\psi^{(i)} = \lambda_{i} \: \boldsymbol\psi^{(i)}, \ \ i = 1,...,n_{s}.
\label{Eq:EigA}
\end{equation}
The state-space mode shapes $\boldsymbol\psi^{(i)}$ are converted into the corresponding modes $\widetilde{\boldsymbol\phi}^{(i)}$ in physical space utilizing the output matrix $\mathbf{C}$ as
\begin{equation}
\widetilde{\boldsymbol\phi}^{(i)} = \mathbf{C} \: \boldsymbol\psi^{(i)}.
\label{Eq:SS2PS_Mode}
\end{equation}
Each mode shape vector $\widetilde{\boldsymbol\phi}^{(i)}$ is then scaled using the residue $R_{kk}^{(i)}$ of the driving point FRF $G_{kk}(\omega)$ of the underlying linear system formulated as
\begin{equation}
G_{kk}(\omega) = \sum_{i=1}^{n_{p}} \frac{R_{kk}^{(i)}}{j \: \omega - \lambda_{i}} + \frac{R_{kk}^{(i)^{\text{\scriptsize $\ast$}}}}{j \: \omega - \lambda^{\ast}_{i}} ,
\label{Eq:Residue_k}
\end{equation}
where $k$ is the location of the excited DOF, and where a star denotes the complex conjugate operation. Eq.~(\ref{Eq:Residue_k}) is an overdetermined algebraic system of equations with $n_{p}$ unknowns and as many equations as the number of processed frequency lines. The $i$-th scaled mode shape vector at the driving point $\phi_{k}^{(i)}$ is finally obtained by enforcing a unit modal mass, \textit{i.e.}
\begin{equation}
R_{kk}^{(i)} =  \frac{\phi_{k}^{(i)} \: \phi_{k}^{(i)}}{2 \: j \: \omega_{i,0}} ,
\label{Eq:ScaledMode}
\end{equation}
while the other components of the mode shape vector are scaled accordingly. This mode scaling is rigorously valid in the case of linear proportional damping~\cite{Geradin_Book}, which implies a real-valued mode shape at the driving point $\phi_{k}^{(i)}$. In general, experimental mode shapes are however complex-valued. They can be enforced to be real by rotating each mode in the complex plane by an angle equal to the mean of the phase angles of the mode components, and subsequently neglecting the imaginary parts of the rotated components.

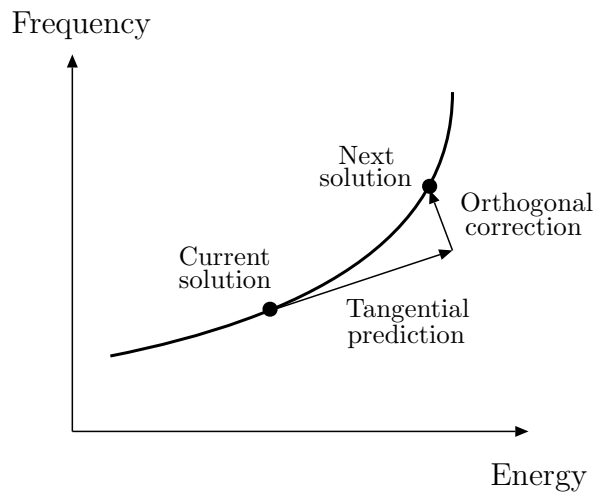
\begin{figure}[ht]
\begin{center}
\scalebox{1} 
{
\begin{pspicture}(0,0)(8,7)
%
\psline[arrows=->,linewidth=0.02,arrowinset=0,arrowlength=1.2,arrowsize=3pt 2](1,1)(7,1)
\psline[arrows=->,linewidth=0.02,arrowinset=0,arrowlength=1.2,arrowsize=3pt 2](1,1)(1,6)
\rput[bl](0.2,6.2){Frequency}
\rput[bl](6.5,0.2){Energy}
\psbezier[linecolor=black, linewidth=0.04](1.5,2.0)(4,2.5)(6,3.5)(6,5.5)
\rput[bl]{45}(3.6,2.62){\psdots[dotsize=0.20](0,0)}
\psline[arrows=->,linewidth=0.02,arrowinset=0,arrowlength=1.2,arrowsize=3pt 2](3.6,2.6)(6.0,3.4)
\rput[bl](2.4,3.2){\footnotesize Current}
\rput[bl](2.4,2.9){\footnotesize solution}
\rput[bl](4.6,2.45){\footnotesize Tangential}
\rput[bl](4.62,2.1){\footnotesize prediction}
\psline[arrows=->,linewidth=0.02,arrowinset=0,arrowlength=1.2,arrowsize=3pt 2](6.0,3.4)(5.7,4.2)
\rput[bl]{45}(5.7,4.25){\psdots[dotsize=0.20](0,0)}
\rput[bl](4.5,4.55){\footnotesize Next}
\rput[bl](4.25,4.25){\footnotesize solution}
\rput[bl](6.1,3.85){\footnotesize Orthogonal}
\rput[bl](6.17,3.6){\footnotesize correction}
\end{pspicture} 
}
\caption{Computation of a family of periodic solutions using a pseudo-arclength continuation scheme including prediction and correction steps.}
\label{Fig:Continuation}
\end{center}
\end{figure}

The algorithm described in Ref.~\cite{NNM_Part2} can now be applied to seek NNMs, given that all quantities in Eq.~(\ref{Eq:TD_ModalModel_Undamped}) are known. To obtain the family of periodic solutions that describe the considered NNM, shooting is combined with a pseudo-arclength continuation technique. Starting from a known periodic solution, continuation proceeds in two steps, namely a prediction and a correction, as illustrated in Fig.~\ref{Fig:Continuation}. In the prediction step, a guess of the next periodic solution along the NNM branch is generated in the direction of the tangent vector to the branch at the current solution. Next, the prediction is corrected using a shooting procedure, forcing the variations of the period and the initial conditions to be orthogonal to the prediction direction.

\newpage
\section{Numerical demonstration using a cantilever beam possessing a hardening-softening nonlinearity}\label{Sec:Demo}

In this section, the methodology is demonstrated based on numerical experiments carried out on the nonlinear beam structure proposed as a benchmark during the European COST Action F3~\cite{Thouverez_NLBeam}. This structure consists of a main cantilever beam whose free end is connected to a thin beam clamped on the other side. The thin beam can exhibit geometrically nonlinear behavior for sufficiently large displacements.  

\begin{figure}[ht]
\begin{center}
\scalebox{1} 
{
\begin{pspicture}(0,0)(12.5,3)
%
\psline[linewidth=0.05](0.5,1.5)(0.5,2.75)
\psline[linewidth=0.02](0.5,2.5)(11,2.5)
\psline[linewidth=0.02](11,2.5)(11,1.75)
\psline[linewidth=0.02](11,1.75)(0.5,1.75)
%
\psline[linewidth=0.02](1.25,1.75)(1.25,2.5)
\psline[linewidth=0.02](2.00,1.75)(2.00,2.5)
\psline[linewidth=0.02](2.75,1.75)(2.75,2.5)
\psline[linewidth=0.02](3.50,1.75)(3.50,2.5)
\psline[linewidth=0.02](4.25,1.75)(4.25,2.5)
\psline[linewidth=0.02](5.00,1.75)(5.00,2.5)
\psline[linewidth=0.02](5.75,1.75)(5.75,2.5)
\psline[linewidth=0.02](6.50,1.75)(6.50,2.5)
\psline[linewidth=0.02](7.25,1.75)(7.25,2.5)
\psline[linewidth=0.02](8.00,1.75)(8.00,2.5)
\psline[linewidth=0.02](8.75,1.75)(8.75,2.5)
\psline[linewidth=0.02](9.50,1.75)(9.50,2.5)
\psline[linewidth=0.02](10.25,1.75)(10.25,2.5)
\psline[linewidth=0.02](11.00,1.75)(11.00,2.5)
\rput[bl](0.8,2){1}
\rput[bl](1.55,2){2}
\rput[bl](2.3,2){3}
\rput[bl](3.05,2){4}
\rput[bl](3.8,2){5}
\rput[bl](4.55,2){6}
\rput[bl](5.3,2){7}
\rput[bl](6.05,2){8}
\rput[bl](6.8,2){9}
\rput[bl](7.4,2){10}
\rput[bl](8.15,2){11}
\rput[bl](8.9,2){12}
\rput[bl](9.65,2){13}
\rput[bl](10.4,2){14}
%
\psline[linewidth=0.02](11,2.20)(11.9,2.20)
\psline[linewidth=0.02](11.9,2.20)(11.9,2.05)
\psline[linewidth=0.02](11.9,2.05)(11,2.05)
\psline[linewidth=0.05](11.9,1.5)(11.9,2.75)
%
\psline[linewidth=0.02](11.3,2.05)(11.3,2.20)
\psline[linewidth=0.02](11.6,2.05)(11.6,2.20)
%

\psline[linewidth=0.05](10.375,0.2)(11.625,0.2)
\rput[bl]{90}(11,0.20){\psline[linecolor=black,linewidth=0.02](0,0)(0.3,0)(0.42,-0.25)(0.54,0.25)(0.66,-0.25)(0.78,0.25)(0.90,-0.25)(1.02,0.25)(1.14,-0.25)(1.26,0)(1.56,0)}
\rput[bl]{-15}(10.3,1.1){\psline[arrows=->,linewidth=0.02,arrowinset=0,arrowlength=1.2,arrowsize=3pt 2](0.5,-0.4)(1.05,0.6)}
\rput[bl](9.5,0.8){$c_{1}$, $c_{2}$}
\end{pspicture} 
}
\caption{Finite element model of the nonlinear beam.}
\label{Fig:FEM_NLBeam}
\end{center}
\end{figure}
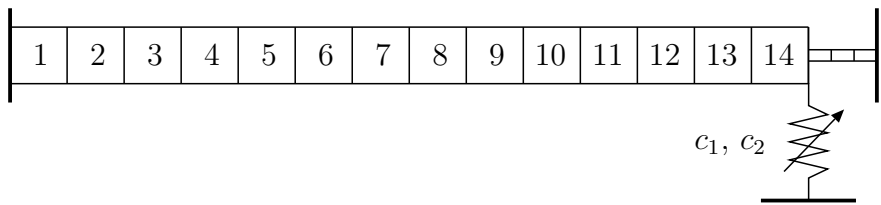

The linear finite element model of the structure, shown in Fig.~\ref{Fig:FEM_NLBeam}, is identical to the linear model experimentally updated in Ref.~\cite{PRM_Part2}. It comprises 14 two-dimensional beam elements for the main beam and 3 elements for the thin beam. The geometrical and mechanical properties of the structure are listed in Tables~\ref{Table:GeoProperties_NLBeam} and~\ref{Table:MechProperties_NLBeam}, respectively. As discussed in Refs.~\cite{Noel_FNSI,Kerschen_CRP_NLBeam}, the nonlinear dynamics induced by the thin beam can be represented through a grounded cubic spring associated with a coefficient $c_{1}$, and positioned at the connection between the two beams. Moreover, it was observed in Ref.~\cite{Grappasonni_NLBeam_IMAC2014} that a slight asymmetry in the clamping conditions of the thin beam may be responsible for significant softening distortions in the system response. This physics is modeled herein by adding a quadratic grounded spring with a negative coefficient $c_{2}$ in parallel to the cubic spring element.

\begin{table}[ht]
\begin{center}
\begin{tabular*}{1.00\textwidth}{@{\extracolsep{\fill}} c c c c}
\hline
					& Length ($m$) & Width ($mm$) & Thickness ($mm$) \\
Main beam & 0.7 & 14 & 14 \\
Thin beam & 0.04 & 14 & 0.5 \\
\hline
\end{tabular*}
\caption{Geometrical properties of the nonlinear beam.}
\label{Table:GeoProperties_NLBeam}
\end{center}
\end{table}

\begin{table}[ht]
\begin{center}
\begin{tabular*}{1.00\textwidth}{@{\extracolsep{\fill}} c c c c}
\hline
Young's modulus ($N/m^{2}$) & Density ($kg/m^{3}$) & $c_{1}$ ($N/m^{3}$) & $c_{2}$ ($N/m^{2}$) \\
$2.05 \: 10^{11}$ & 7800 & $8 \: 10^{9}$ & $-1.05 \: 10^{7}$ \\
\hline
\end{tabular*}
\caption{Mechanical properties of the nonlinear beam.}
\label{Table:MechProperties_NLBeam}
\end{center}
\end{table}

According to Fig.~\ref{Fig:NPS_Methodology}, the testing procedure applies a broadband excitation signal to the nonlinear beam structure. The FNSI method can address classical random excitations, including Gaussian noise, periodic random, burst random and pseudo random signals. Impulsive excitations also fall within the scope of the method, but they involve windowing to avoid leakage and generally lead to low signal-to-noise ratios (SNRs). Swept-sine excitations are not applicable because of the inability of the FNSI method to handle nonstationary signals, \textit{i.e.} signals with time-varying frequency content~\cite{Noel_FNSI}. One opts herein for pseudo random signals, also known as random phase multisine signals. A random phase multisine is a periodic random signal with a user-controlled amplitude spectrum, and a random phase spectrum drawn from a uniform distribution. If an integer number of periods is measured, the amplitude spectrum is perfectly realized, unlike Gaussian noise. One of the other main advantages of a multisine is that its periodic nature can be utilized to separate transient from steady-state oscillations in response time histories. This, in turn, eliminates the systematic errors due to leakage in the identification. Periodicity also allows the estimation of the covariance matrix of the noise perturbations affecting the system outputs.

\begin{figure}[ht]
\begin{center}
\begin{tabular}{c c}
\subfloat[]{\label{NLBeam_Force_1P_Amp}\includegraphics[width=75mm]{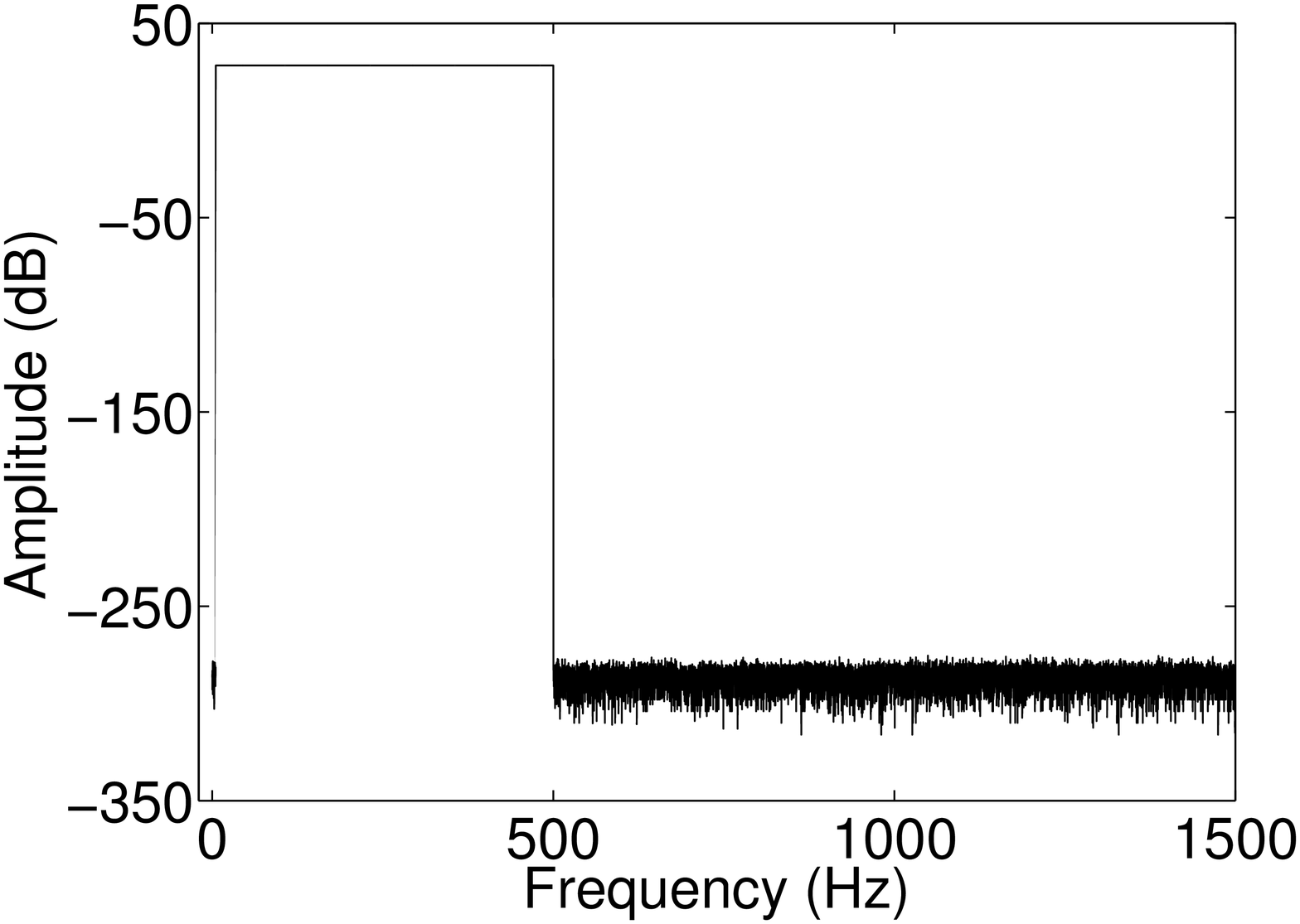}} &
\subfloat[]{\label{NLBeam_Force_1P_Phase}\includegraphics[width=72mm]{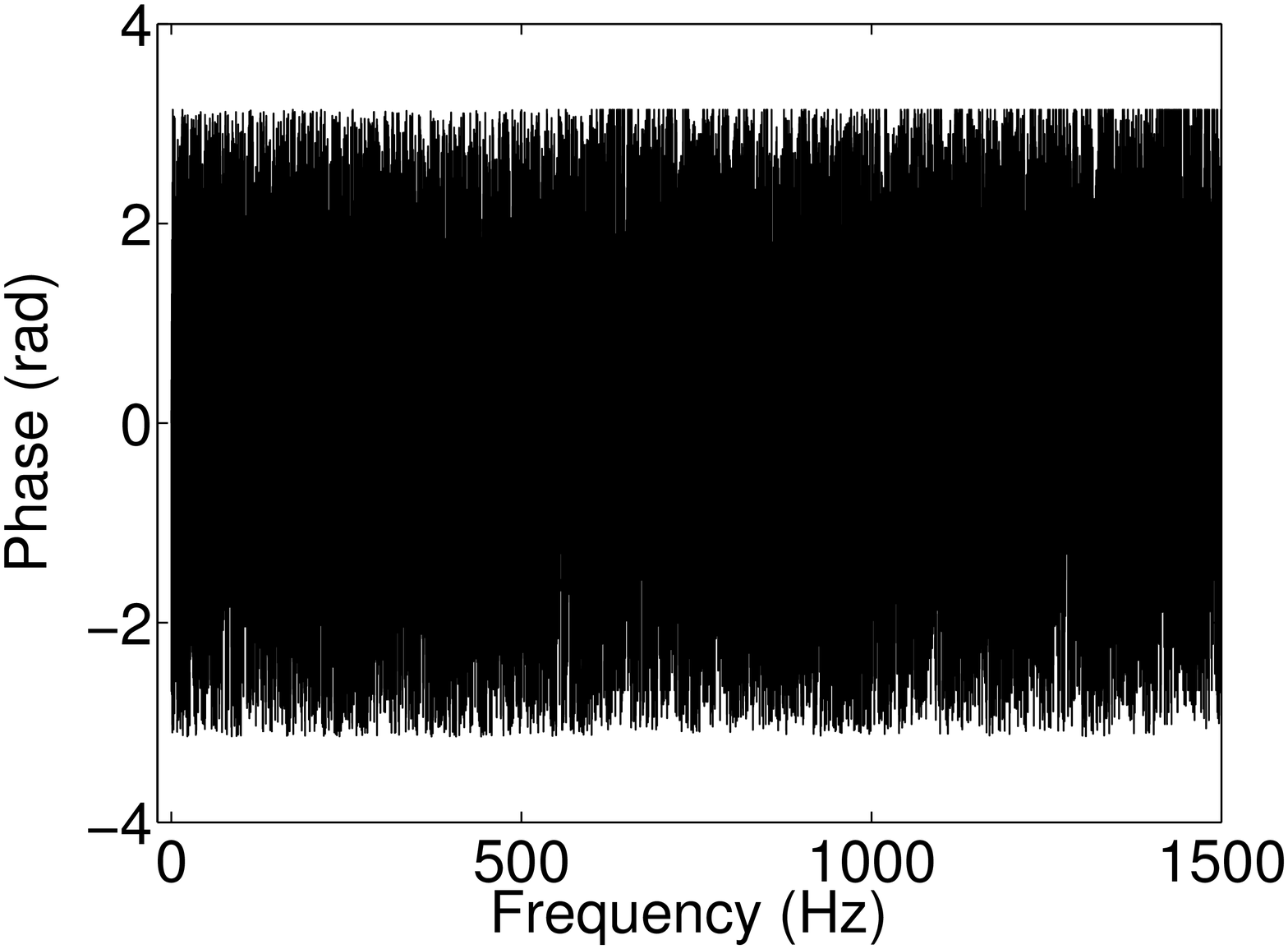}} \\
\subfloat[]{\label{NLBeam_Force_AllP}\includegraphics[width=72mm]{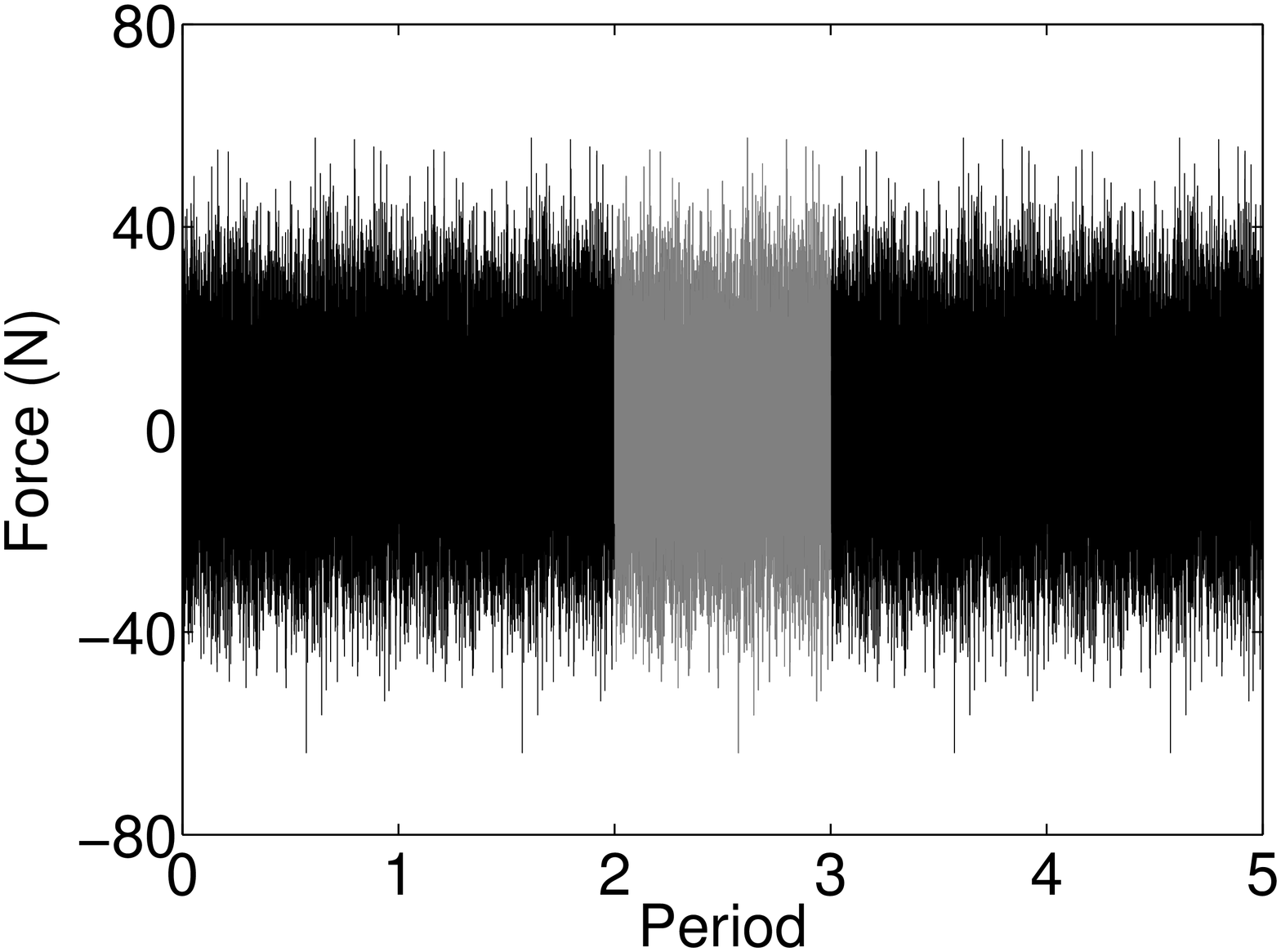}} &
\subfloat[]{\label{NLBeam_Force_TransientdB}\includegraphics[width=73mm]{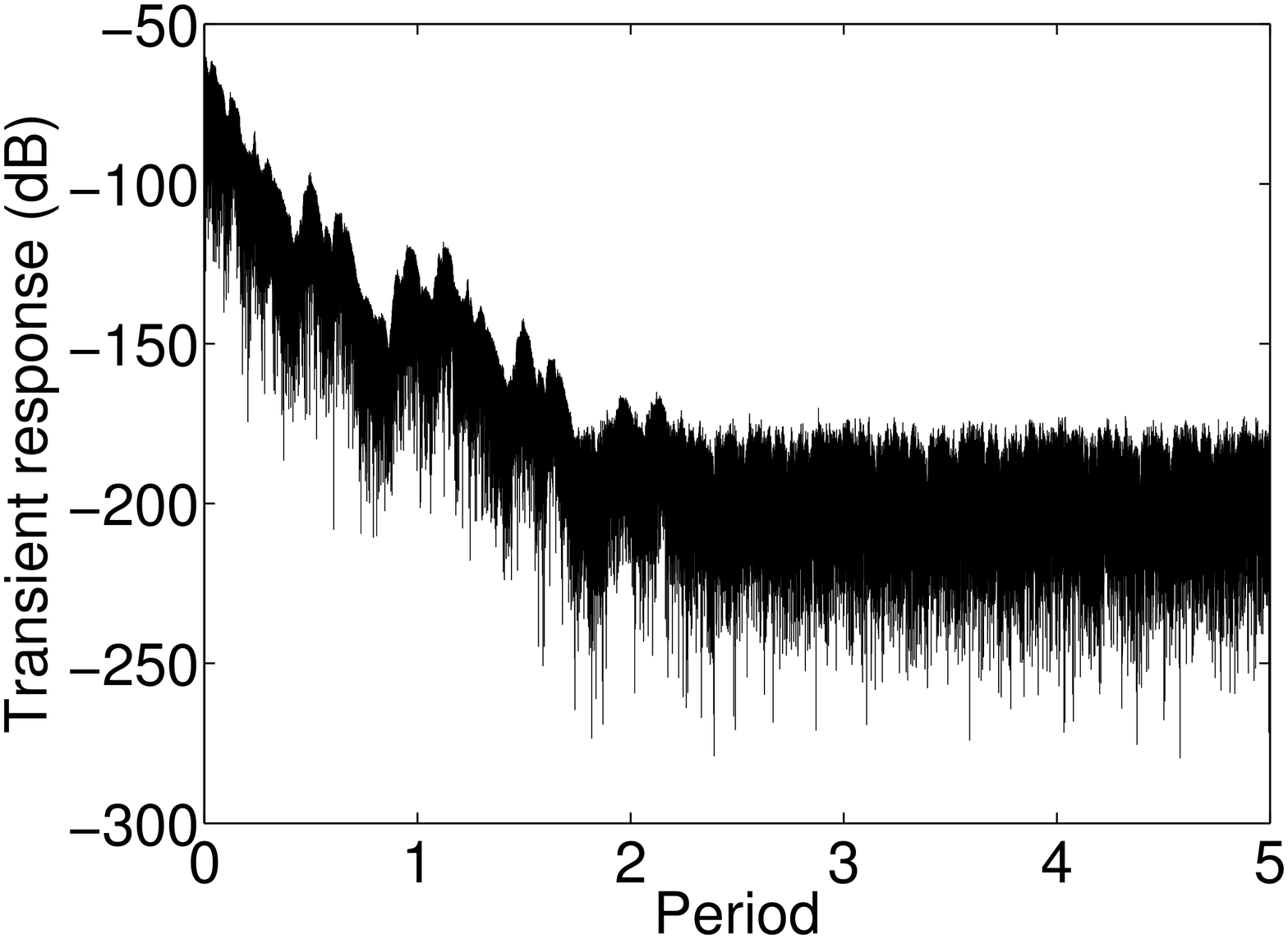}} \\
\end{tabular}
\caption{Random phase multisine excitation signal. (a -- b) Amplitude and phase spectrum of a single period; (c) first 5 periods in the time domain with one specific period highlighted in gray; (d) decay of the transient system response illustrated using a logarithmic scaling at the main beam tip over 5 periods.}
\label{Fig:NLBeam_Force}
\end{center}
\end{figure}

A multisine excitation with a flat amplitude spectrum and a root-mean-squared (RMS) amplitude of 15 $N$ was applied vertically to node 4 of the structure (see Fig.~\ref{Fig:FEM_NLBeam}). The excited band spans the 5 -- 500 $Hz$ interval to encompass the three linear modes of interest. The response of the nonlinear beam to this excitation was simulated over 20 periods of $2^{15}=32768$ samples each. Fig.~\ref{Fig:NLBeam_Force}~(a -- b) shows the amplitude and the phase spectrum of one period of the multisine input. The first 5 periods of the signal in the time domain are also depicted in Fig.~\ref{Fig:NLBeam_Force}~(c), where one specific period is highlighted in gray.

Numerical experiments were conducted by direct time integration using a nonlinear Newmark scheme. To this end, a linear proportional damping matrix $\mathbf{C}_{v} = \alpha \: \mathbf{K} + \beta \: \mathbf{M}$, with $\alpha = 3\:10^{-7}$ and $\beta = 5$, was introduced in the model. The resulting linear natural frequencies and damping ratios of the first three bending modes of the beam structure are given in Table~\ref{Table:NLBeam_LinProp}. The sampling frequency during time simulation was set to 60000 $Hz$ to ensure the accuracy of the integration. Synthetic time series were then decimated down to 3000 $Hz$ for practical use, considering low-pass filtering to avoid aliasing.  The decay of the transient system response is illustrated using a logarithmic scaling in Fig.~\ref{Fig:NLBeam_Force}~(d). It is seen to die out after 3 periods. This latter plot was generated at the main beam tip by subtracting from the entire measured signal its last period, \textit{i.e.} its twentieth period, assumed to be in steady state.

\begin{table}[ht]
\begin{center}
\begin{tabular*}{1.00\textwidth}{@{\extracolsep{\fill}} c c c}
\hline
Mode & Natural frequency $\omega_{0}$ ($Hz$) & Damping ratio $\zeta$ ($\%$) \\
		 &  &  \\
1 & 31.28  & 1.28 \\
2 & 143.64 & 0.29 \\
3 & 397.87 & 0.14 \\
\hline
\end{tabular*}
\caption{Linear natural frequencies $\omega_{0}$ and damping ratios $\zeta$ of the first three bending modes of the nonlinear beam.}
\label{Table:NLBeam_LinProp}
\end{center}
\end{table}

\begin{figure}[ht]
\begin{center}
\includegraphics[width=100mm]{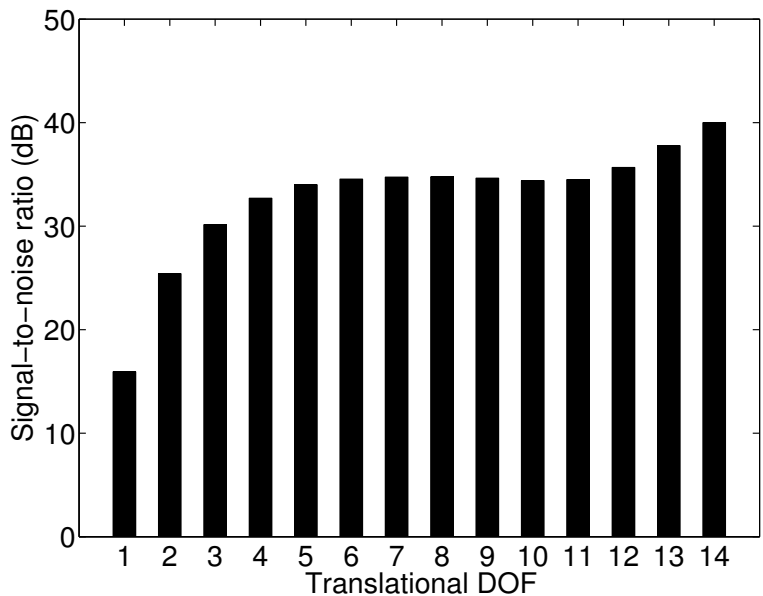} \\
\caption{Signal-to-noise ratio calculated at the translational DOFs along the main beam. The noise level is set to 1 $\%$ (= 40 $dB$) of the RMS displacement amplitude at beam tip, \textit{i.e.} at translational DOF 14.}
\label{Fig:NLBeam_NSR}
\end{center}
\end{figure}

Simulated time series were finally corrupted by adding white noise, recreating the mechanical and electrical disturbances observed in a typical measurement setup. The noise level was set to 1 $\%$ of the RMS displacement amplitude at the main beam tip. The resulting SNR calculated at the translational DOFs along the beam is plotted in Fig.~\ref{Fig:NLBeam_NSR}. It is found that imposing 1 $\%$ noise at the tip, \textit{i.e.} a SNR of 40 $dB$ at DOF 14 in Fig.~\ref{Fig:NLBeam_NSR}, results in more severe noise conditions at all other sensors. In particular, at mid-span, the SNR is around 34 $dB$, while it is lower than 30 $dB$ close to the left clamping. 

\subsection{Identification using the FNSI method}\label{Sec:Demo_FNSI}

\subsubsection{Selection of the nonlinear basis functions}\label{Sec:Demo_Characterization}

The application of the FNSI method to measured data requires the selection of appropriate basis functions $\mathbf{h}_{a}(\mathbf{q}(t))$ to describe the nonlinearity in the system. This task, referred to as the characterization of nonlinearity, is in general challenging because of the various sources of nonlinear behavior that may exist in engineering structures, and the plethora of dynamic phenomena they may cause~\cite{Kerschen_Review}. In this work, a gray-box methodology is adopted by exploiting cubic splines, \textit{i.e.} piecewise third-order polynomials~\cite{DeBoor_SplinesBook}. Cubic splines have the advantage of being as simple as ordinary polynomials, while overcoming some of their drawbacks~\cite{Noel_SolarPanels}. The use of splines as nonlinear basis functions is a key aspect of the present study, as it allows NNMs to be identified without assuming the functional forms $\mathbf{h}_{a}(\mathbf{q}(t))$ of the nonlinearities in the system under test. Specifically, the hardening-softening nonlinearity in the beam dynamics is modeled herein using a spline function of the beam tip displacement, considering a division of the measured displacement range into 10, equally-wide segments.

\subsubsection{Selection of the model order}\label{Sec:Demo_ModelOrder}

\begin{figure}[t]
\begin{center}
\includegraphics[width=140mm]{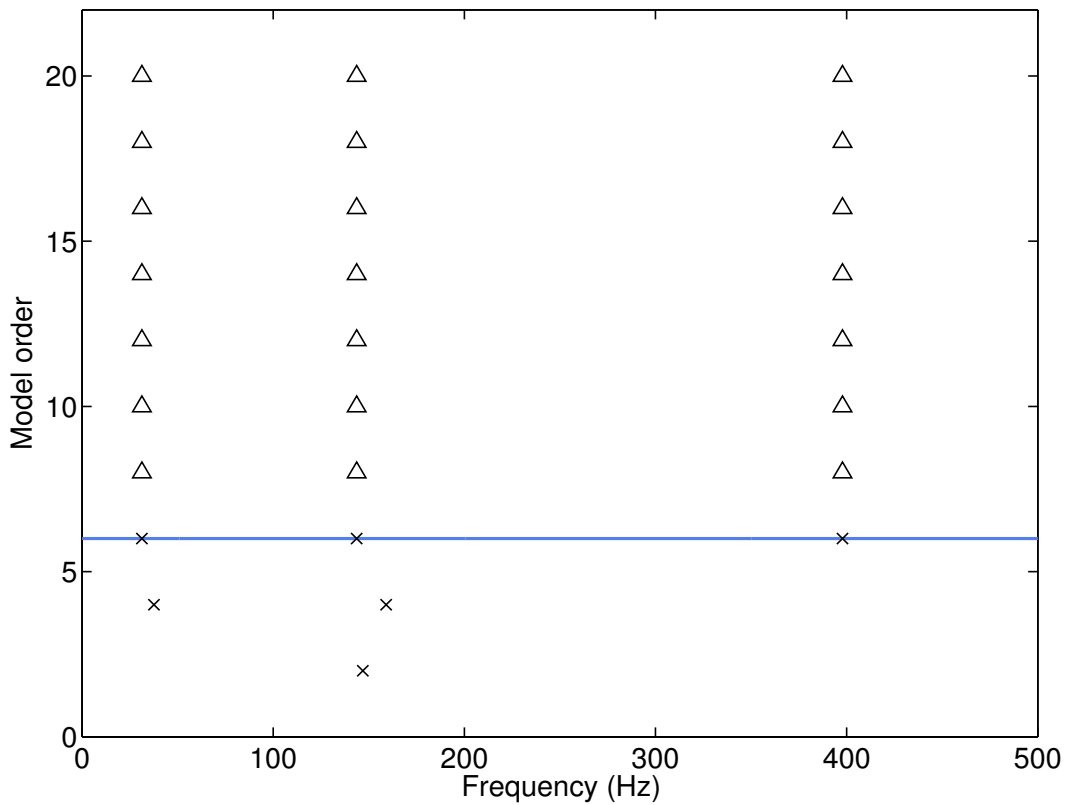} \\
\caption{Stabilization diagram. Cross: stabilization in natural frequency; circle: extra stabilization in damping ratio; triangle: full stabilization. Stabilization thresholds in natural frequency, damping ratio and MAC are 1 $\%$, 5 $\%$ and 0.98, respectively. The blue line indicates the selected order.}
\label{Fig:NLBeam_FNSI_SD}
\end{center}
\end{figure}

The order of the state-space model derived using the FNSI method is equal to twice the number of linear modes activated in the measured data~\cite{Noel_FNSI}. This order is conveniently estimated using a stabilization diagram, similarly to the current practice in linear system identification. Fig.~\ref{Fig:NLBeam_FNSI_SD} charts the stabilization of the natural frequencies, damping ratios and mode shapes of the structure computed at 15 $N$ RMS for model orders up to 20. In this diagram, the modal assurance criterion (MAC) is utilized to quantify the correspondence between mode shapes at different orders. The knowledge of the output noise covariance matrix gained via the periodicity of the excitation was incorporated in the diagram following the discussion in Ref.~\cite{McKelvey_Subspace}. Fig.~\ref{Fig:NLBeam_FNSI_SD} shows full stabilization of three modes in the input band, which leads to the selection of the order 6. Opting for higher model orders would here result in overfitting issues, and would hence obviously deteriorate the identification accuracy. Note that the usual practice in linear system identification is to pick up poles at different model orders. In the nonlinear case, this is no longer possible as a single model order should be selected to further estimate nonlinear coefficients~\cite{Noel_FNSI}.

\subsubsection{Estimation of the nonlinear coefficients}\label{Sec:Demo_EstNLLin}

\begin{figure}[p]
\begin{center}
\includegraphics[width=140mm,tics=50,trim={3.2cm 0 0 0},clip=true]{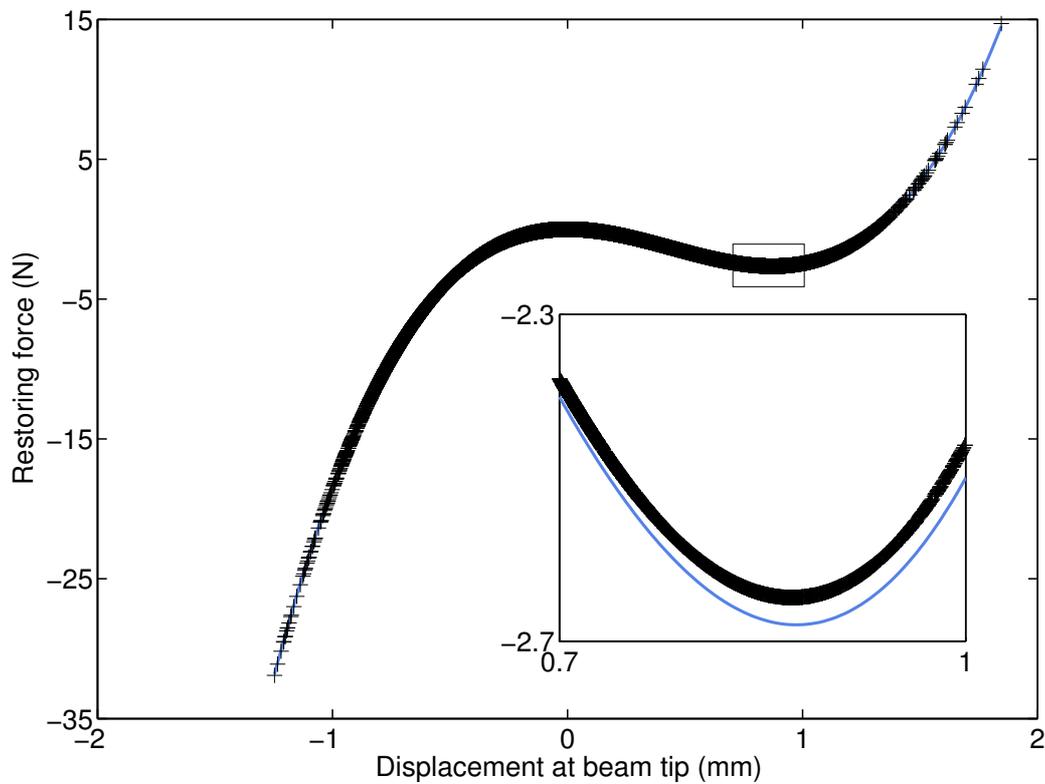} \\
\caption{Exact (in blue) and estimated (black markers) nonlinear restoring force curves. The estimated curve is a cubic spline function of the displacement measured at the main beam tip.}
\label{Fig:NLBeam_FNSI_NLForce}
\end{center}
\end{figure}

\begin{figure}[p]
\begin{center}
\begin{tabular}{c c}
\subfloat[]{\label{NLBeam_FNSI_AbsError}\includegraphics[width=75mm]{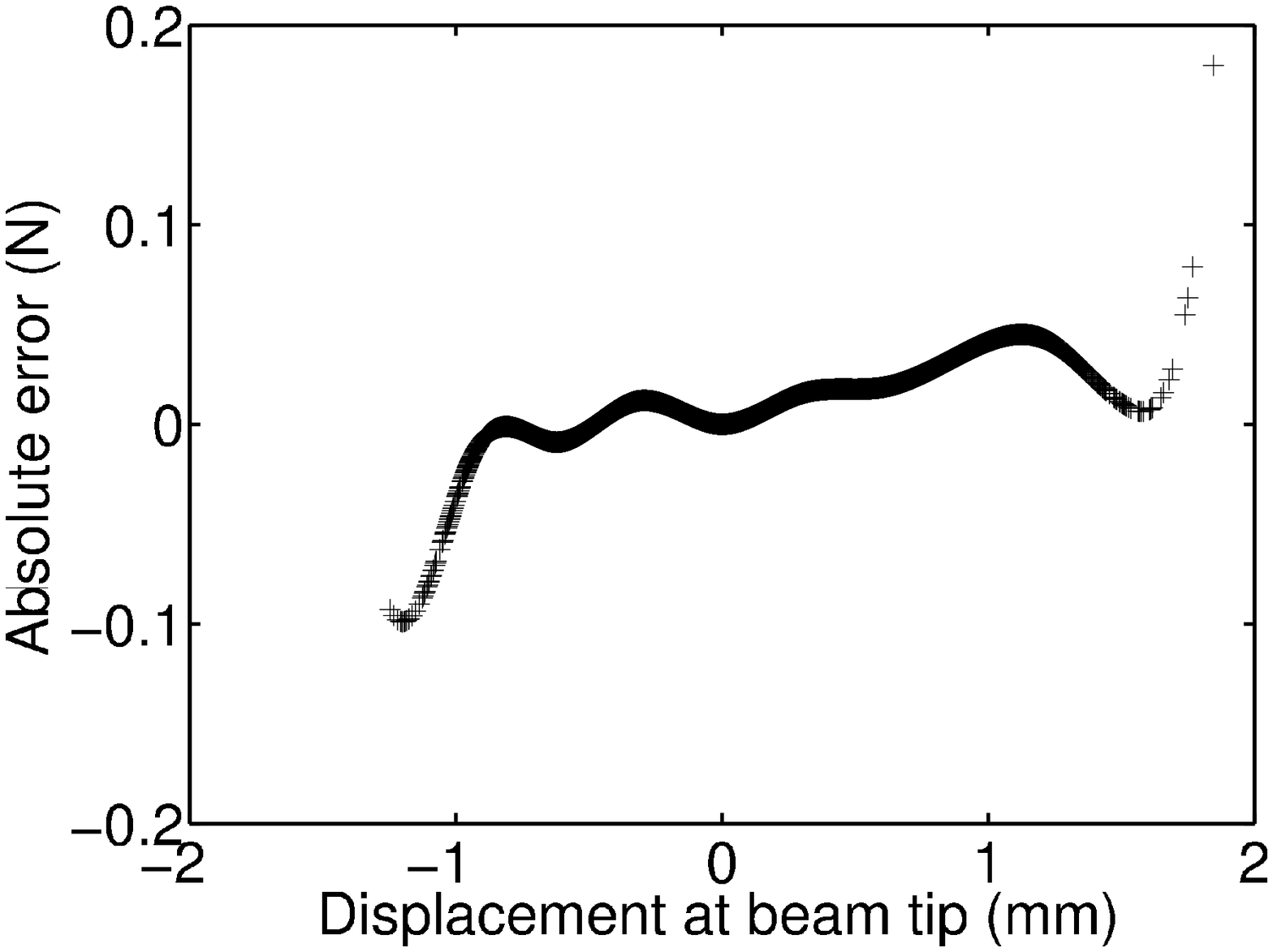}} &
\subfloat[]{\label{NLBeam_FNSI_RatioRI}\includegraphics[width=72mm]{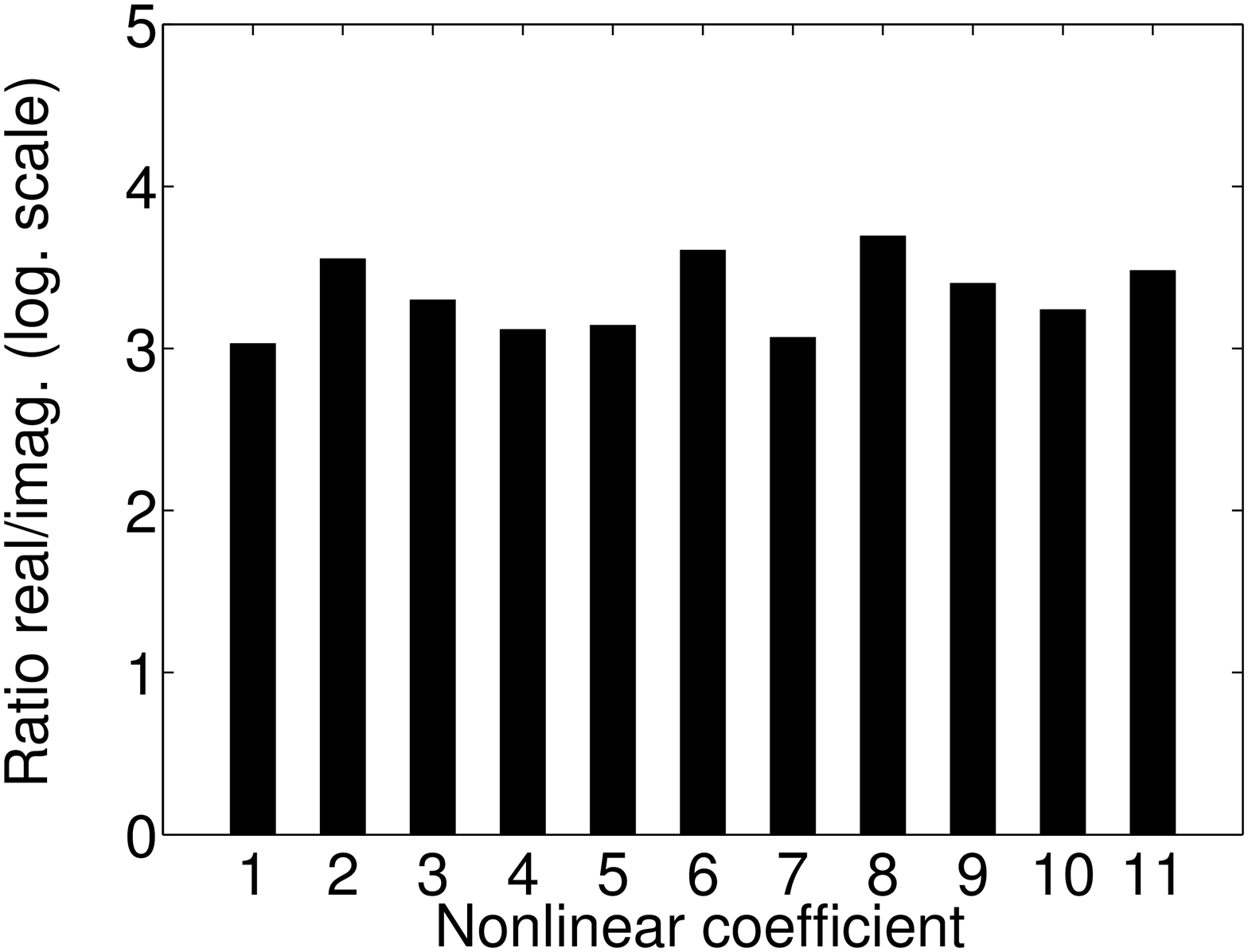}} \\
\end{tabular}
\caption{(a) Absolute error between exact and estimated nonlinear restoring force curves; (b) ratios in logarithmic scaling between real and imaginary parts of the 11 coefficients of the identified spline nonlinearity}
\label{Fig:NLBeam_Error_RI}
\end{center}
\end{figure}

The nonlinear restoring force identified at beam tip is displayed in Fig.~\ref{Fig:NLBeam_FNSI_NLForce}. The agreement between the exact (in blue) and estimated (black markers) force curves is excellent. The inset close-up in Fig.~\ref{Fig:NLBeam_FNSI_NLForce}, together with the absolute error plot in Fig.~\ref{Fig:NLBeam_Error_RI}~(a), reveal that discrepancies are of the order of 0.1 $N$. The 11 coefficients of the identified spline nonlinearity are spectral quantities, \textit{i.e.} complex-valued and frequency-dependent, as a result of Eq.~(\ref{Eq:ExtendedFRF}). Fig.~\ref{Fig:NLBeam_Error_RI}~(b) depicts in logarithmic scaling the ratios between their real and imaginary parts. They are all found to be greater than 3, which confirms the high quality of the parameter estimates.

\subsection{Computation of the first two NNMs using continuation}\label{Sec:Demo_Continuation}

In this section, the first two NNMs of the beam structure are computed by applying the algorithm of Section~\ref{Sec:NNMCont} to Eqs.~(\ref{Eq:TD_ModalModel_Undamped}) populated with the nonlinear and linear parameters estimated using FNSI. Nonlinear parameter estimates were discussed in the previous section. Table~\ref{Table:NLBeam_FNSI_Lin} lists the relative errors on the linear natural frequencies and damping ratios together with the diagonal MAC values. The results in this table demonstrate the ability of the FNSI method to recover accurately the modal properties of the underlying linear structure from nonlinear data. Note that the third mode will not be further analyzed herein as it involves virtually no nonlinear distortions.

\begin{table}[ht]
\begin{center}
\begin{tabular*}{1.00\textwidth}{@{\extracolsep{\fill}} c c c c}
\hline
Mode & Error on $\omega_{0}$ ($\%$) & Error on $\zeta$ ($\%$) & MAC \\
 & & & \\
1  &  0.0008 & -0.0758 & 1.00 \\
2  & -0.0015 &  0.0709 & 1.00 \\
3  & -0.0144 & -0.1015 & 1.00 \\
\hline
\end{tabular*}
\caption{Relative errors on the estimated natural frequencies and damping ratios (in $\%$) and diagonal MAC values of the first three modes of the beam computed at order 6.} 
\label{Table:NLBeam_FNSI_Lin}
\end{center}
\end{table}

\begin{figure}[p]
\begin{center}
\begin{overpic}[width=140mm]
		{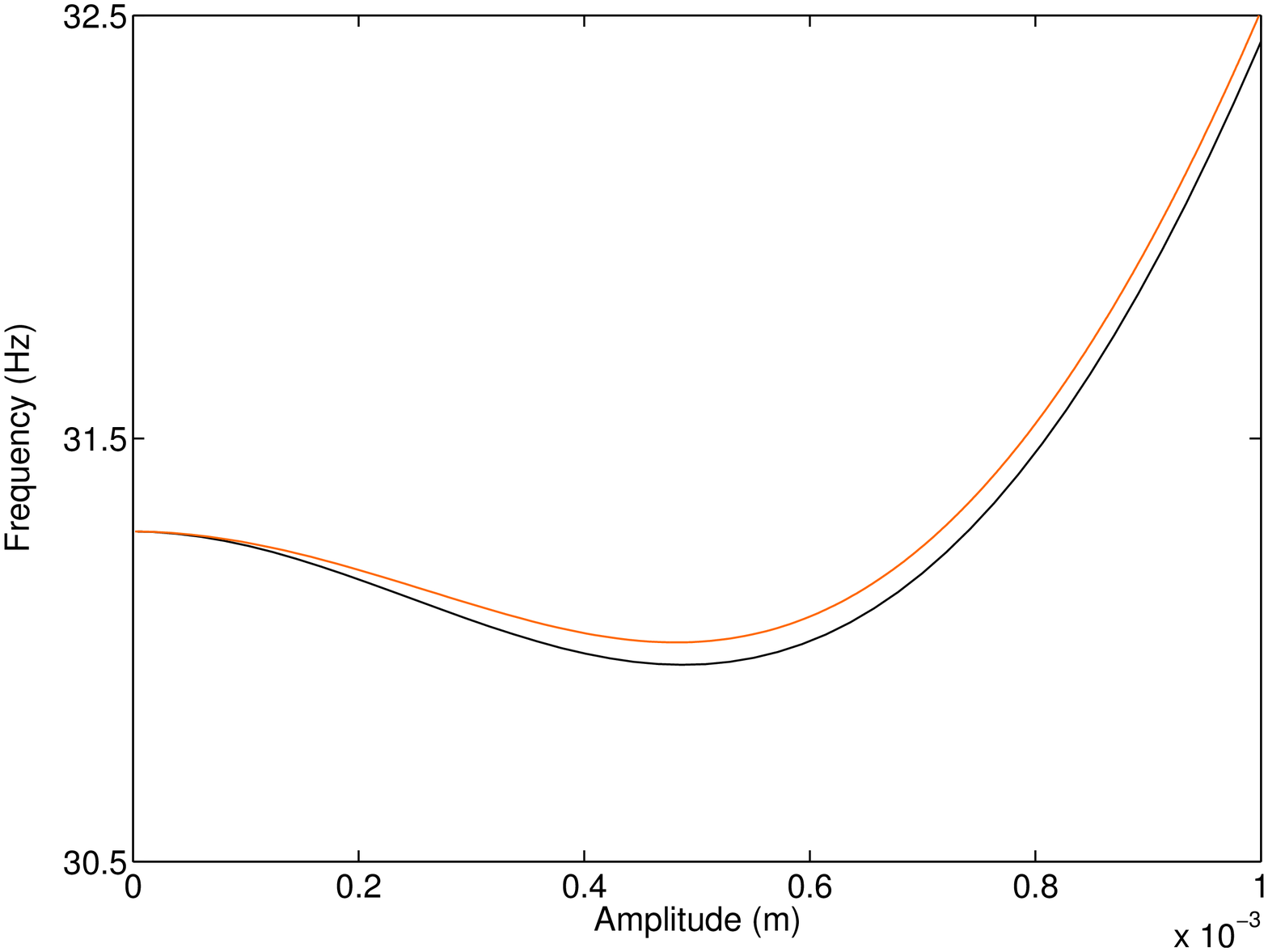}
		\put(140,560){\includegraphics[width=25mm]{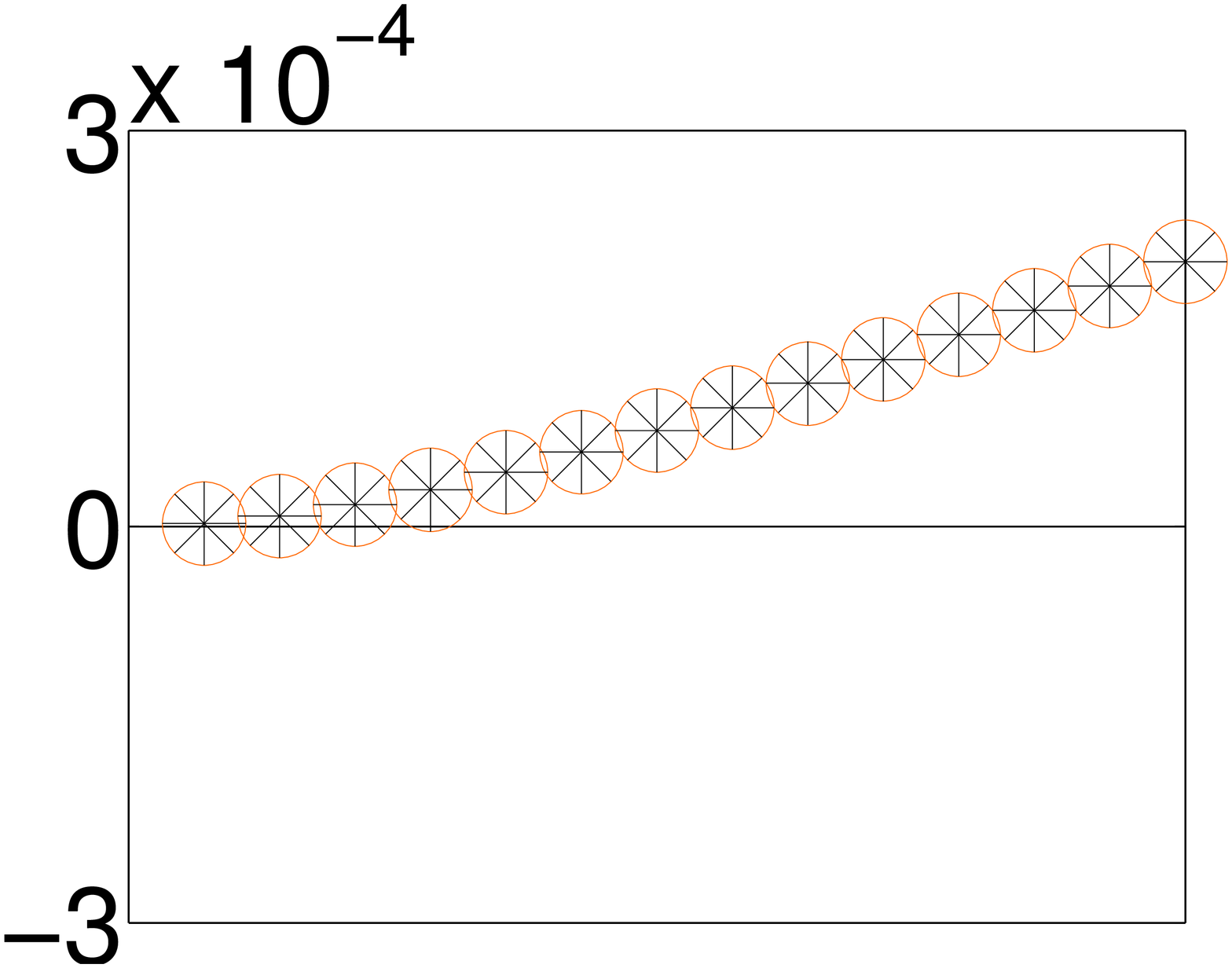}}
		\put(337,560){\includegraphics[width=25mm]{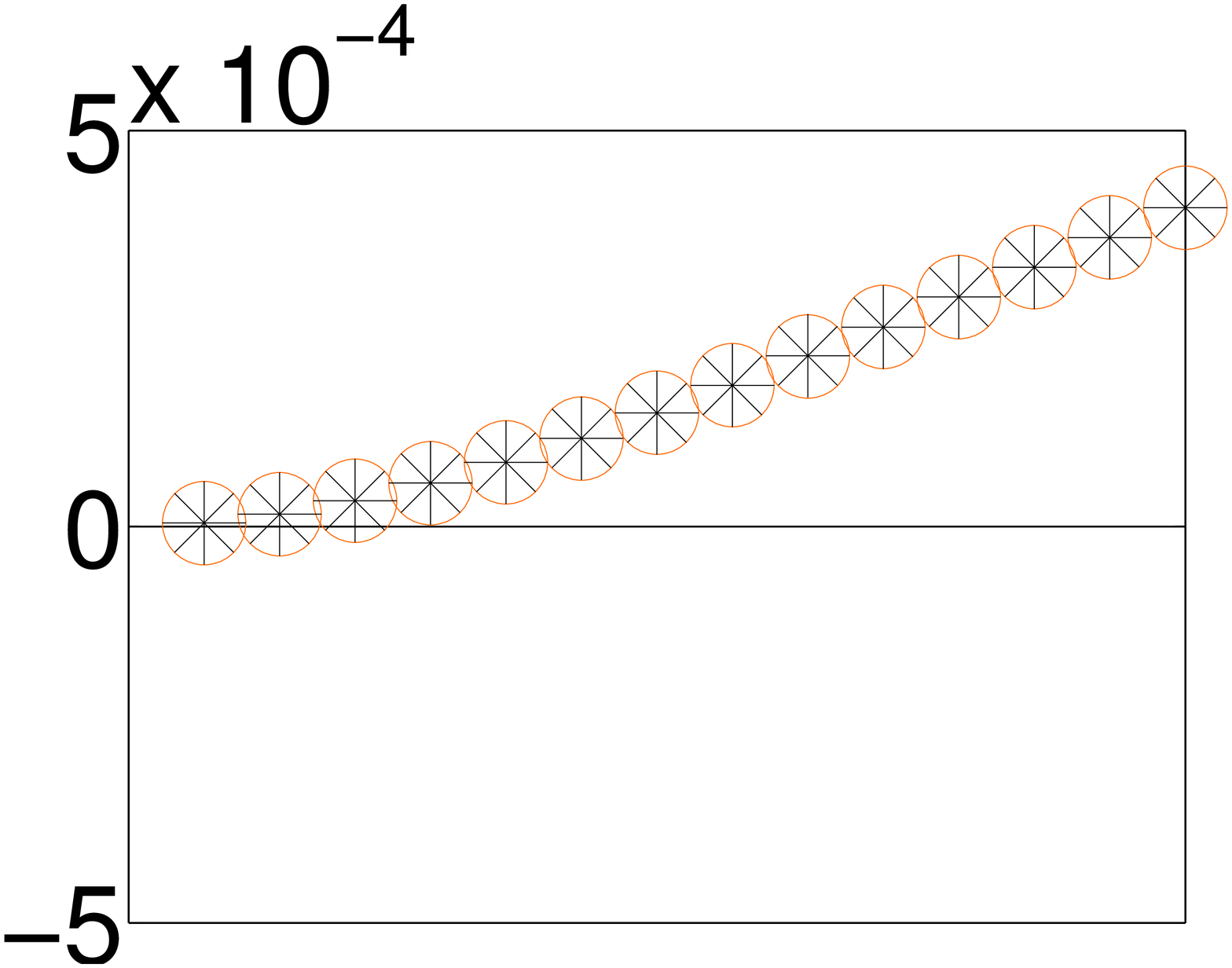}}
		\put(337,363){\includegraphics[width=25mm]{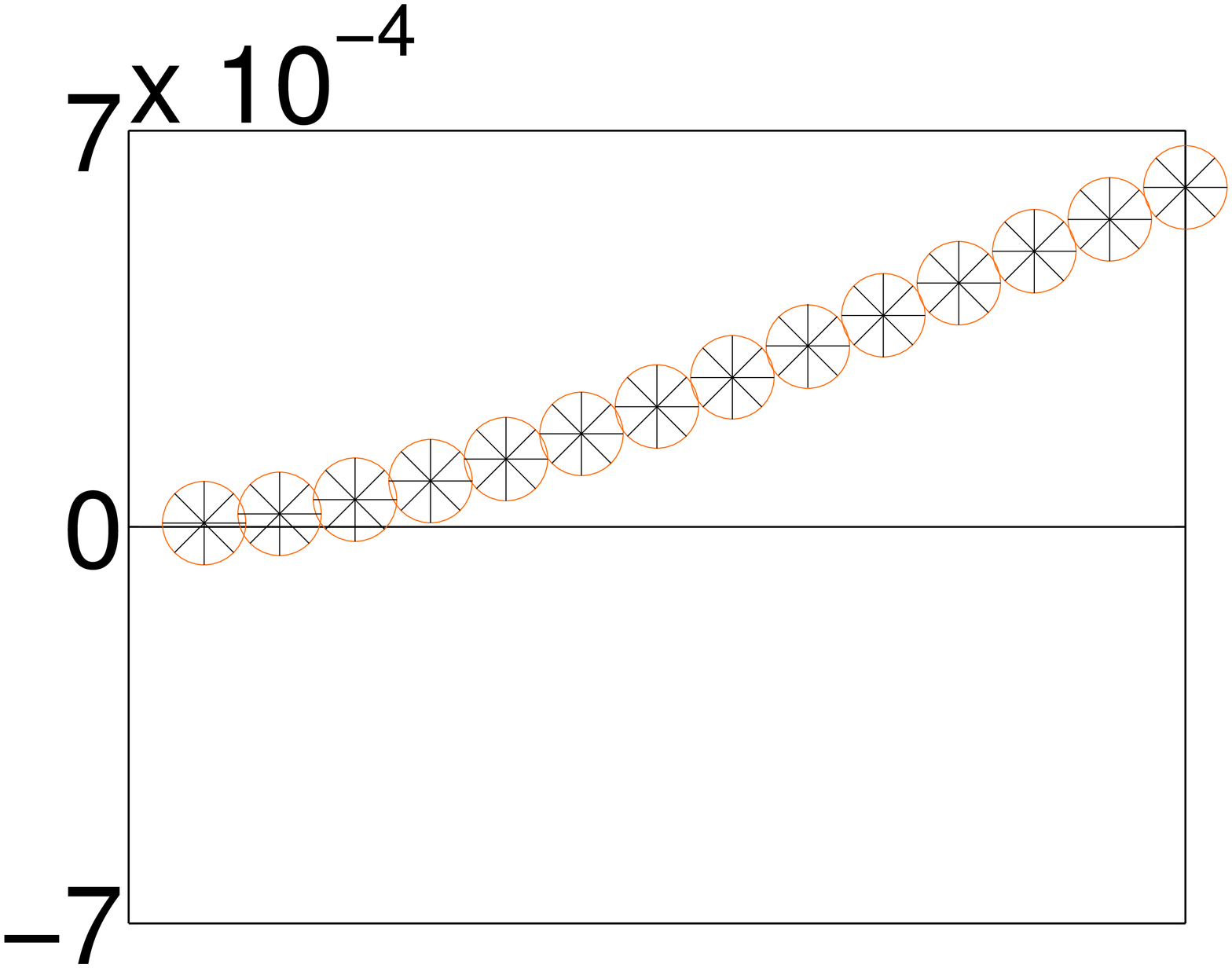}}
		\put(543,363){\includegraphics[width=25mm]{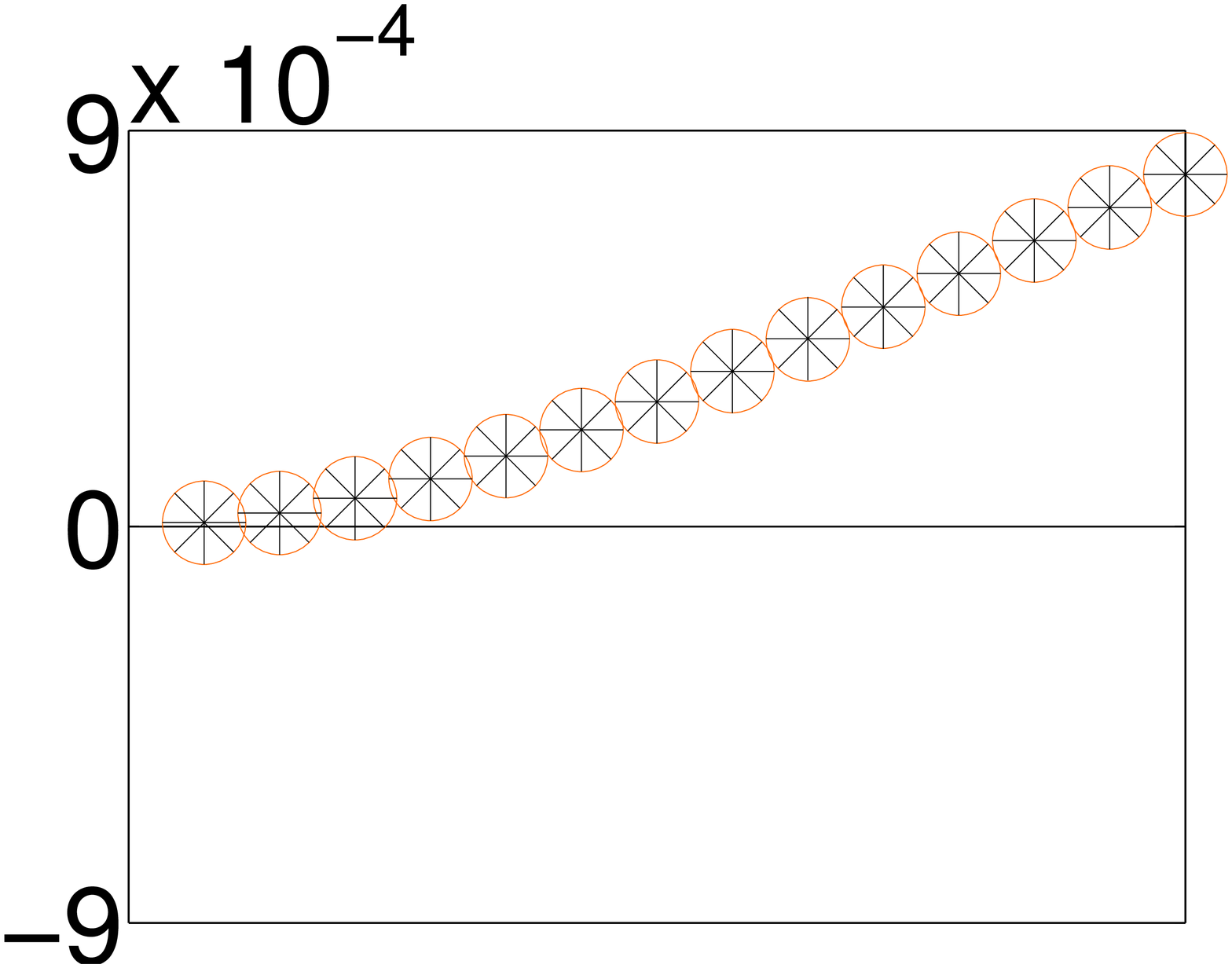}}
		\put(160,535){\footnotesize (a)}
		\put(357,535){\footnotesize (b)}
		\put(357,338){\footnotesize (c)}
		\put(564,338){\footnotesize (d)}
\end{overpic}
\caption{Comparison between the theoretical (in black) and identified (in orange) frequency-amplitude evolution of the first NNM of the nonlinear beam. The NNM shapes (displacement amplitudes of the main beam) at four amplitude levels, namely 0.2, 0.4, 0.6 and 0.8 $mm$, are inset in (a -- d).}
\label{Fig:NLBeam_NNM1_NPS}
\end{center}
\end{figure}

\begin{figure}[p]
\begin{center}
\begin{overpic}[width=140mm]
		{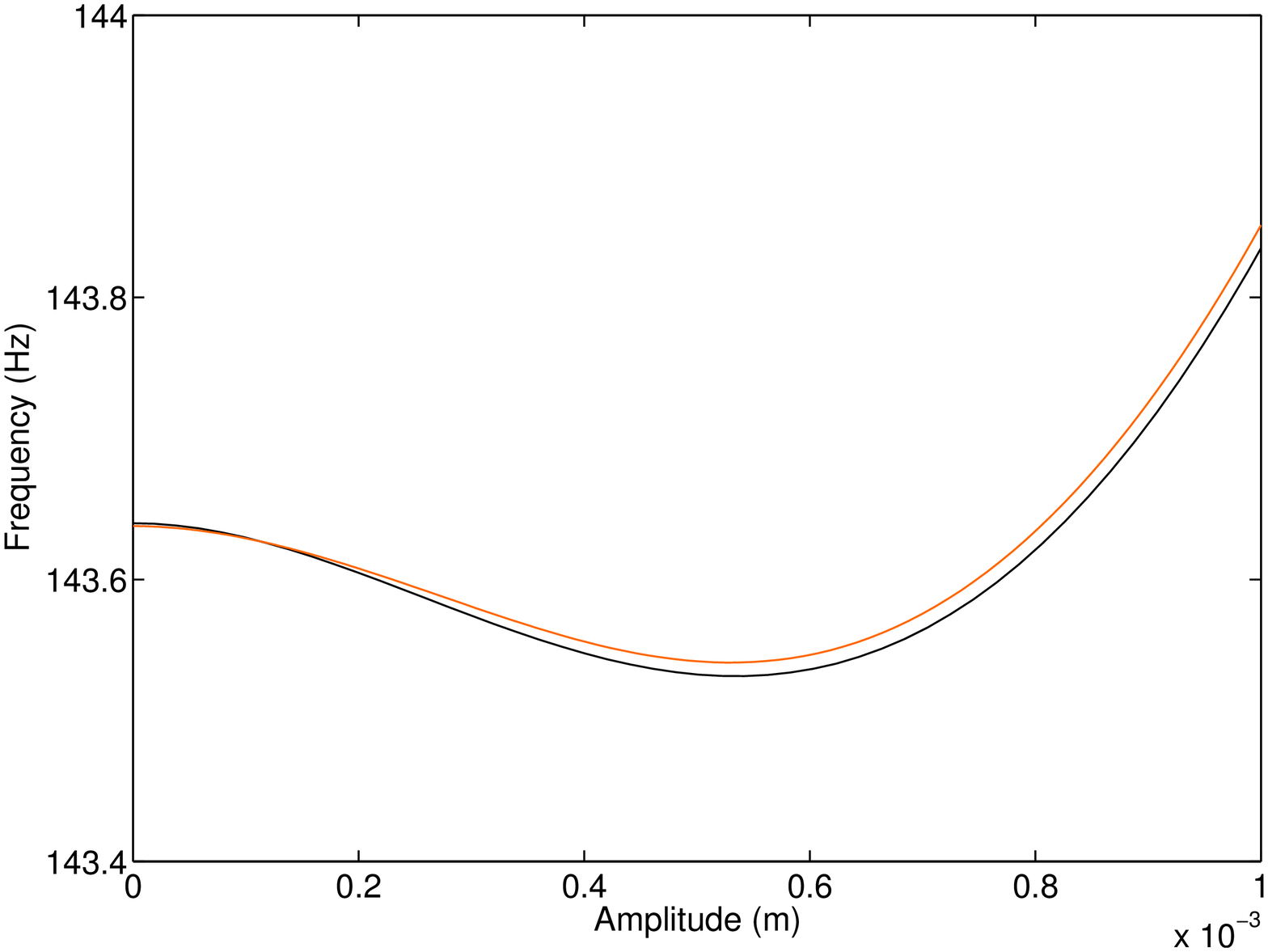}
		\put(140,560){\includegraphics[width=25mm]{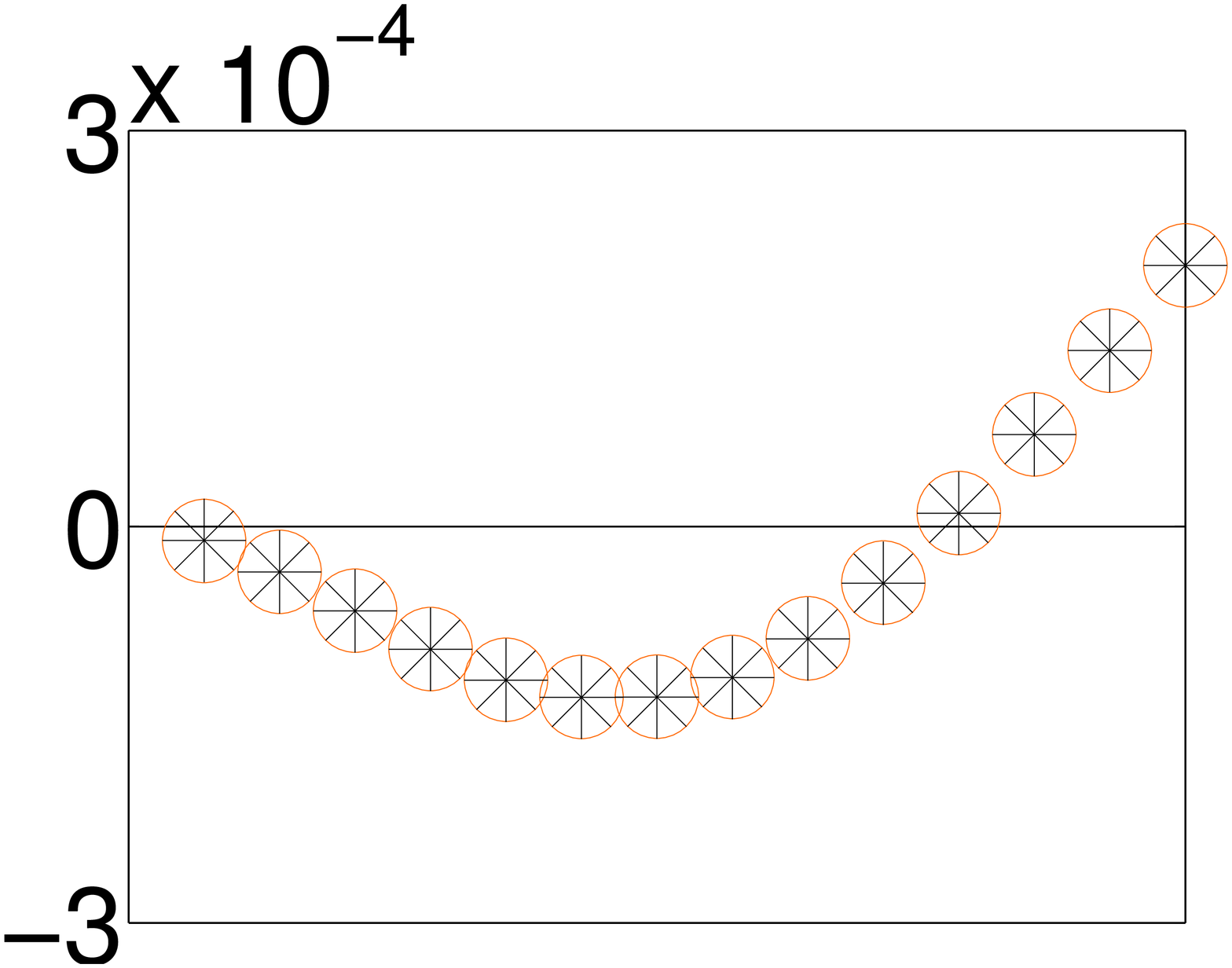}}
		\put(337,560){\includegraphics[width=25mm]{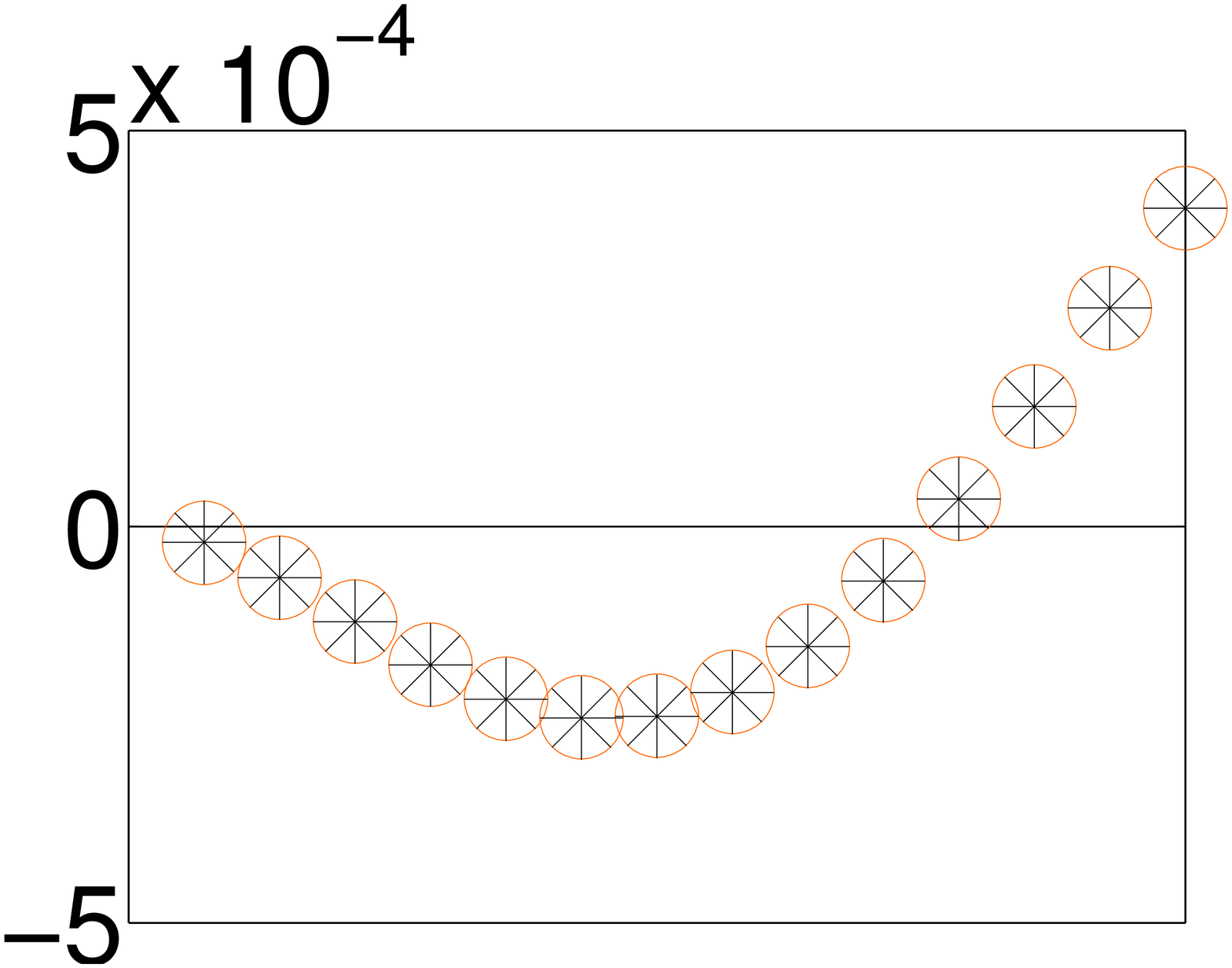}}
		\put(337,363){\includegraphics[width=25mm]{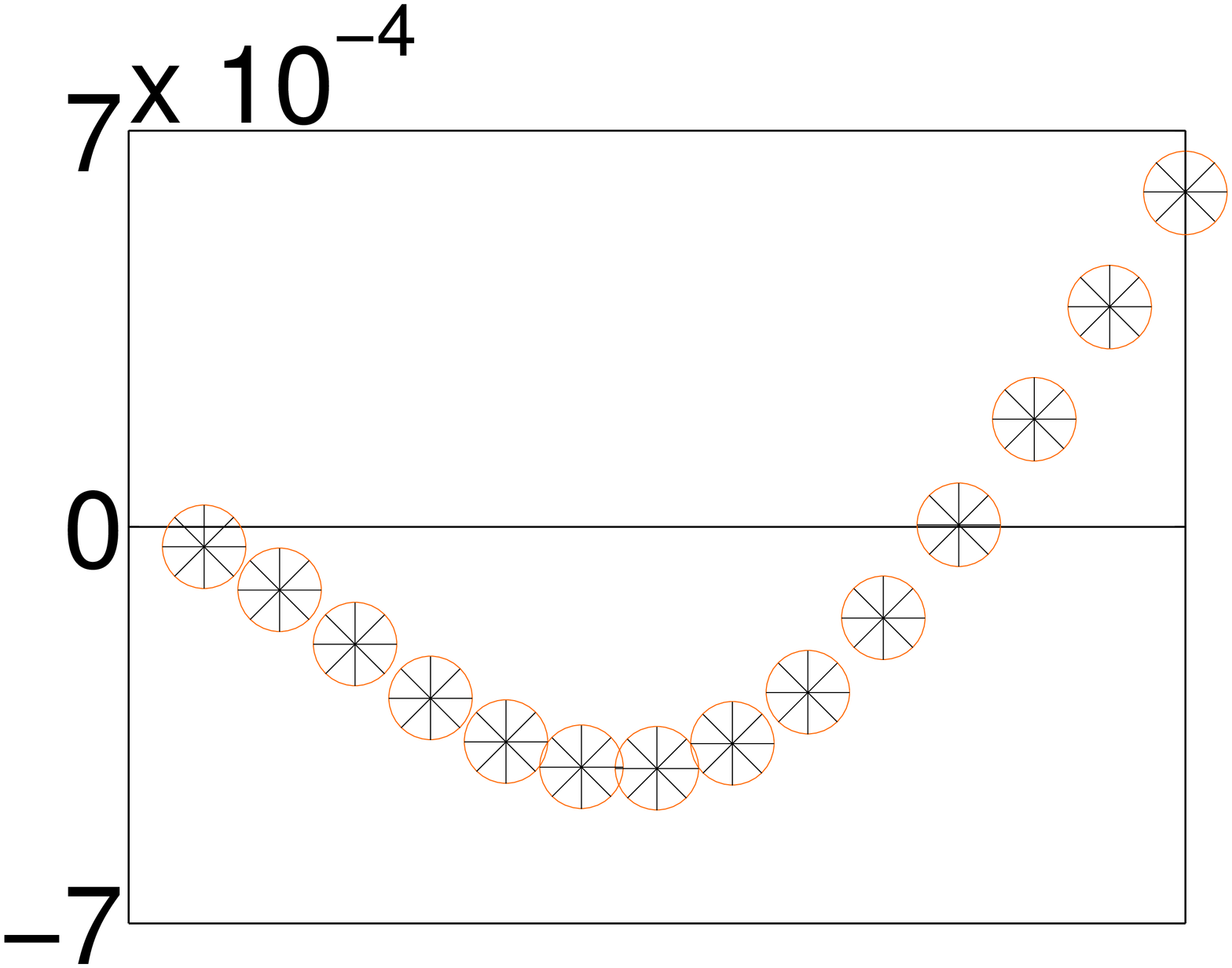}}
		\put(543,363){\includegraphics[width=25mm]{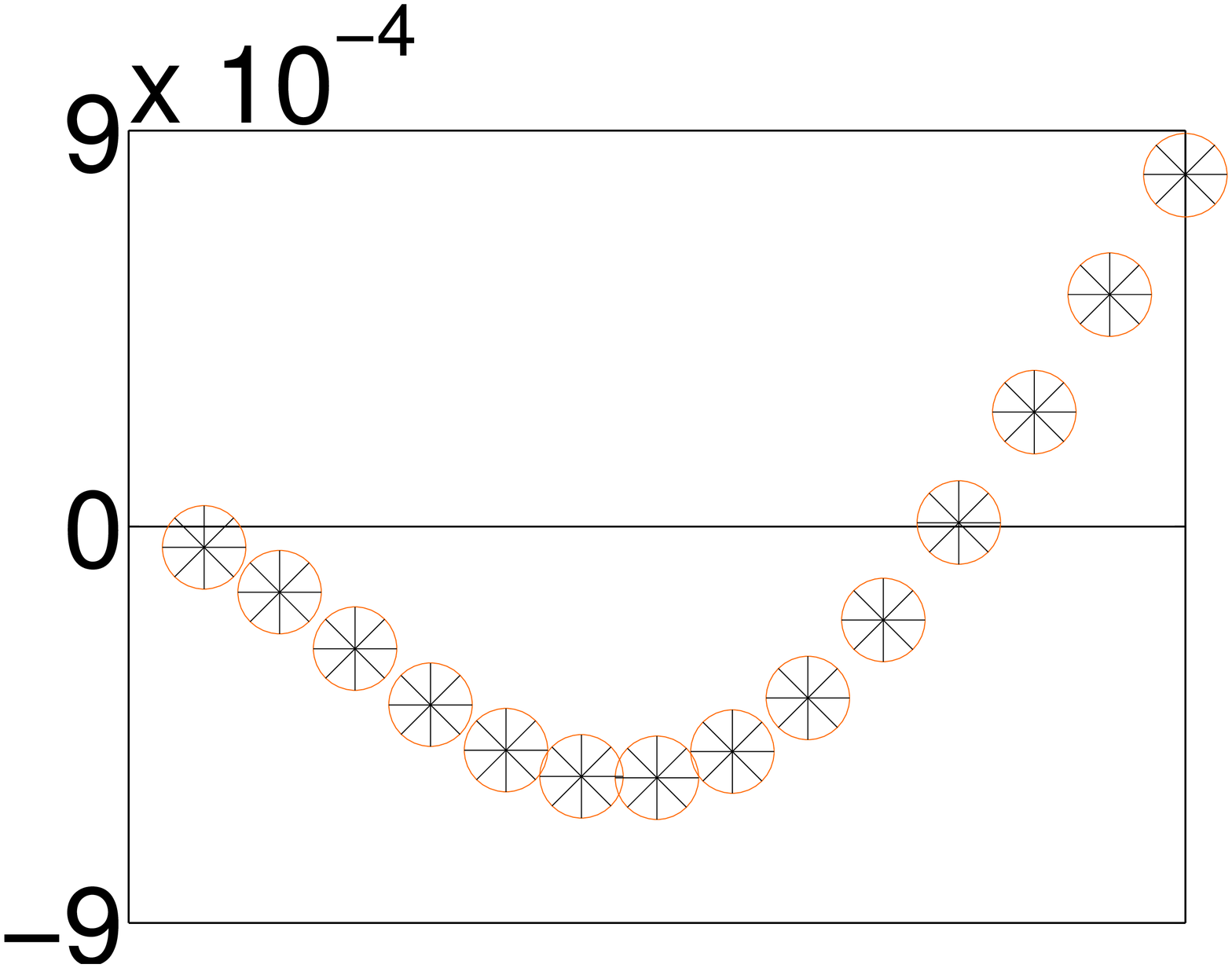}}
		\put(160,535){\footnotesize (a)}
		\put(357,535){\footnotesize (b)}
		\put(357,338){\footnotesize (c)}
		\put(564,338){\footnotesize (d)}
\end{overpic}
\caption{Comparison between the theoretical (in black) and identified (in orange) frequency-amplitude evolution of the second NNM of the nonlinear beam. The NNM shapes (displacement amplitudes of the main beam) at four amplitude levels, namely 0.2, 0.4, 0.6 and 0.8 $mm$, are inset in (a -- d).}
\label{Fig:NLBeam_NNM2_NPS}
\end{center}
\end{figure}

Fig.~\ref{Fig:NLBeam_NNM1_NPS} shows the evolution of the frequency of the first bending mode of the structure as a function of the amplitude of the motion evaluated at the main beam tip. The use of a displacement amplitude as horizontal axis in this plot is justified by the absence of direct access to the total conserved energy associated with the considered mode in experimental conditions~\cite{PRM_Part2}. The identified frequency-amplitude curve depicted in orange in Fig.~\ref{Fig:NLBeam_NNM1_NPS} is seen to closely match the exact NNM presented in black (with a maximum relative error of 0.25 $\%$). The identification accuracy is confirmed through the comparison between the exact and identified modal shapes inserted in Fig.~\ref{Fig:NLBeam_NNM1_NPS} at four specific amplitude levels, namely 0.2, 0.4, 0.6 and 0.8 $mm$. Similar conclusions are drawn from the quality of the identification of the second mode of the nonlinear beam plotted in Fig.~\ref{Fig:NLBeam_NNM2_NPS}. It is important to remark that the introduced methodology makes no approximation throughout, and retrieves exactly the NNMs of the system when the SNR is infinite and the exact nonlinear basis functions are utilized.

\begin{figure}[p]
\begin{center}
\begin{tabular}{c c}
\subfloat[]{\label{NLBeam_NNM1_ConfSpace1}\includegraphics[width=75mm]{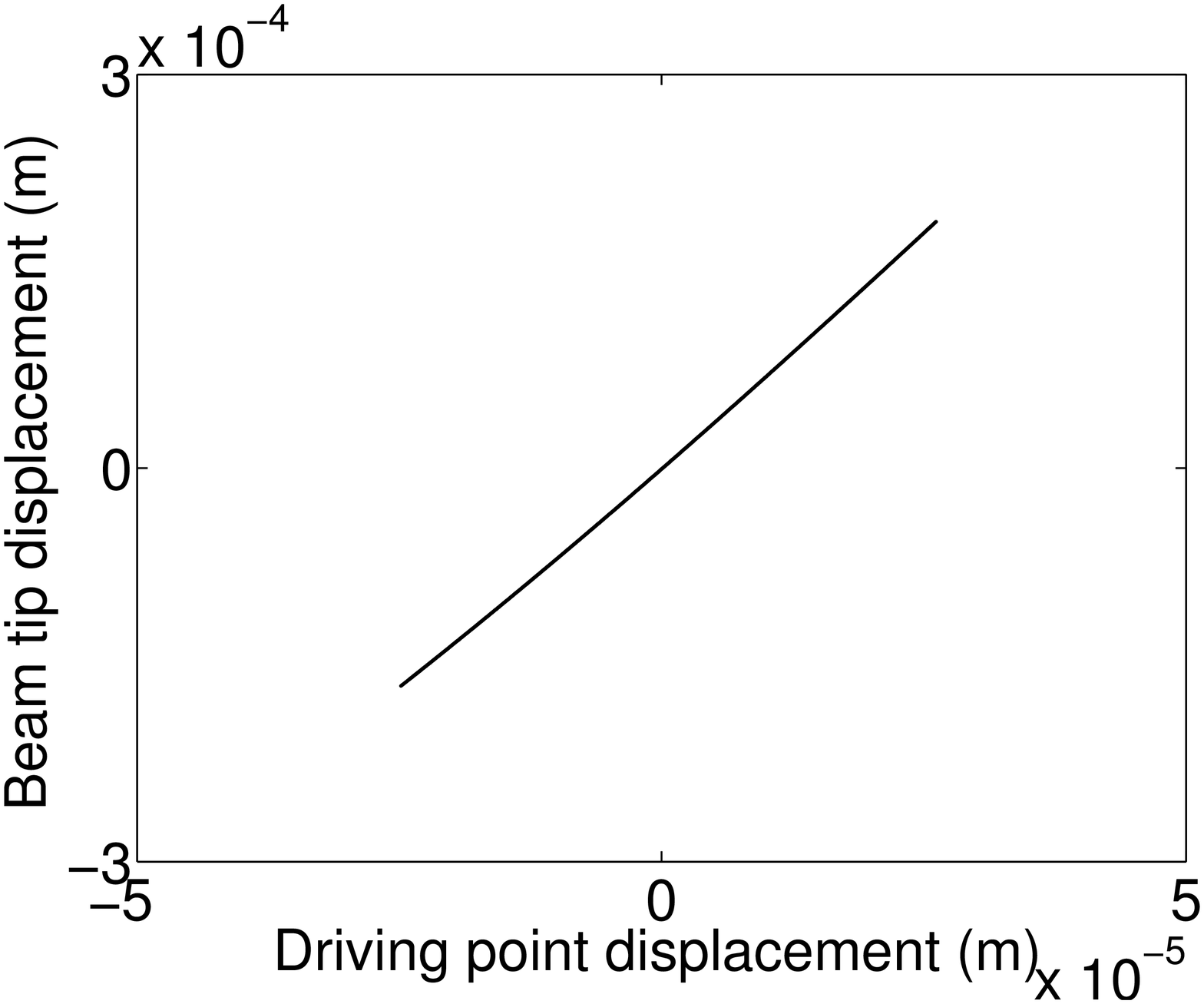}} &
\subfloat[]{\label{NLBeam_NNM1_ConfSpace2}\includegraphics[width=75mm]{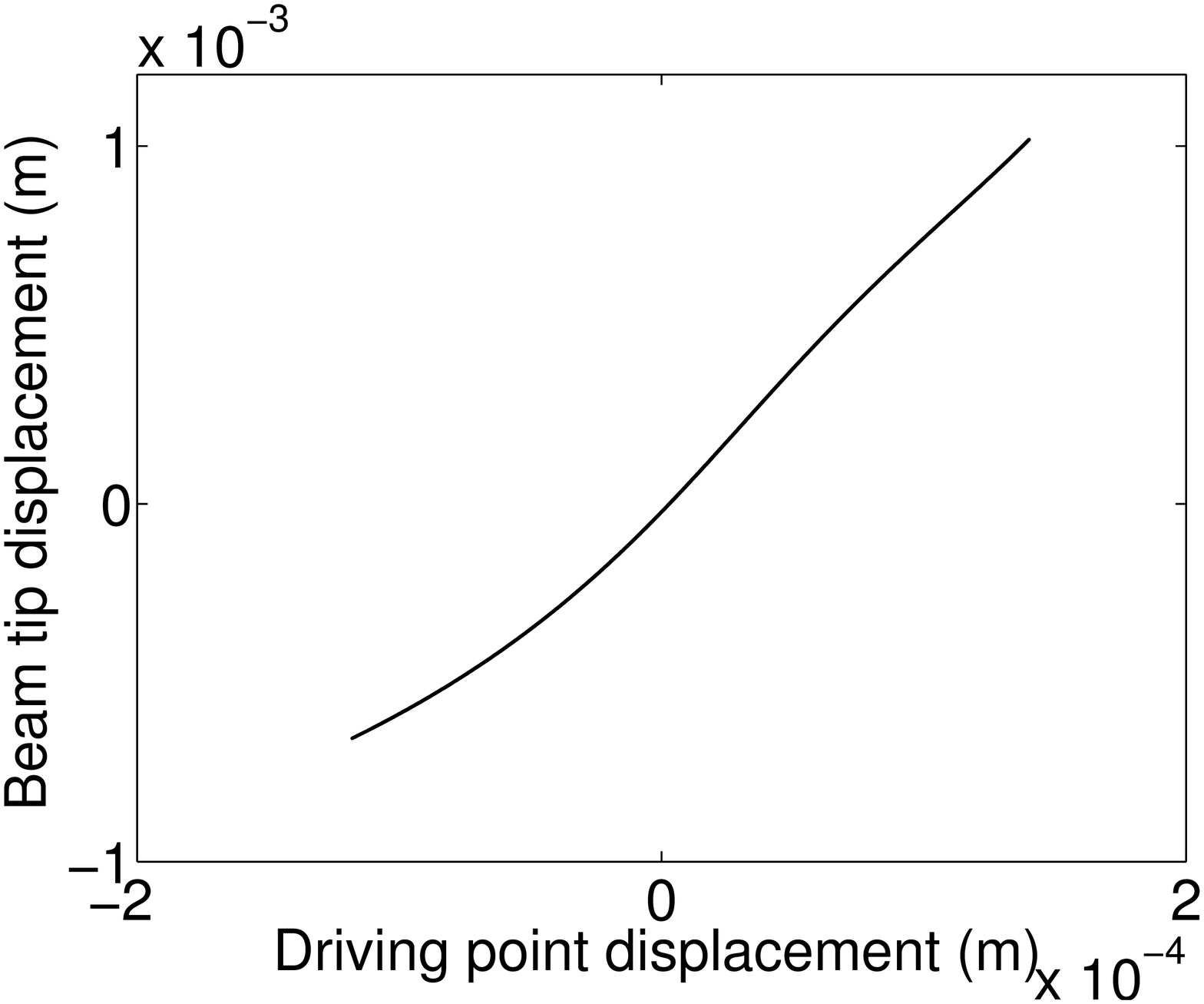}} \\
\end{tabular}
\caption{First NNM of the nonlinear beam represented in configuration space at two amplitude levels. The configuration space is constructed using the displacements measured at the driving point and the main beam tip. (a) 0.2 $mm$ ; (b) 1 $mm$.}
\label{Fig:NLBeam_NNM1_ConfSpace}
\end{center}
\end{figure}

\begin{figure}[p]
\begin{center}
\begin{tabular}{c c}
\subfloat[]{\label{NLBeam_NNM2_ConfSpace1}\includegraphics[width=75mm]{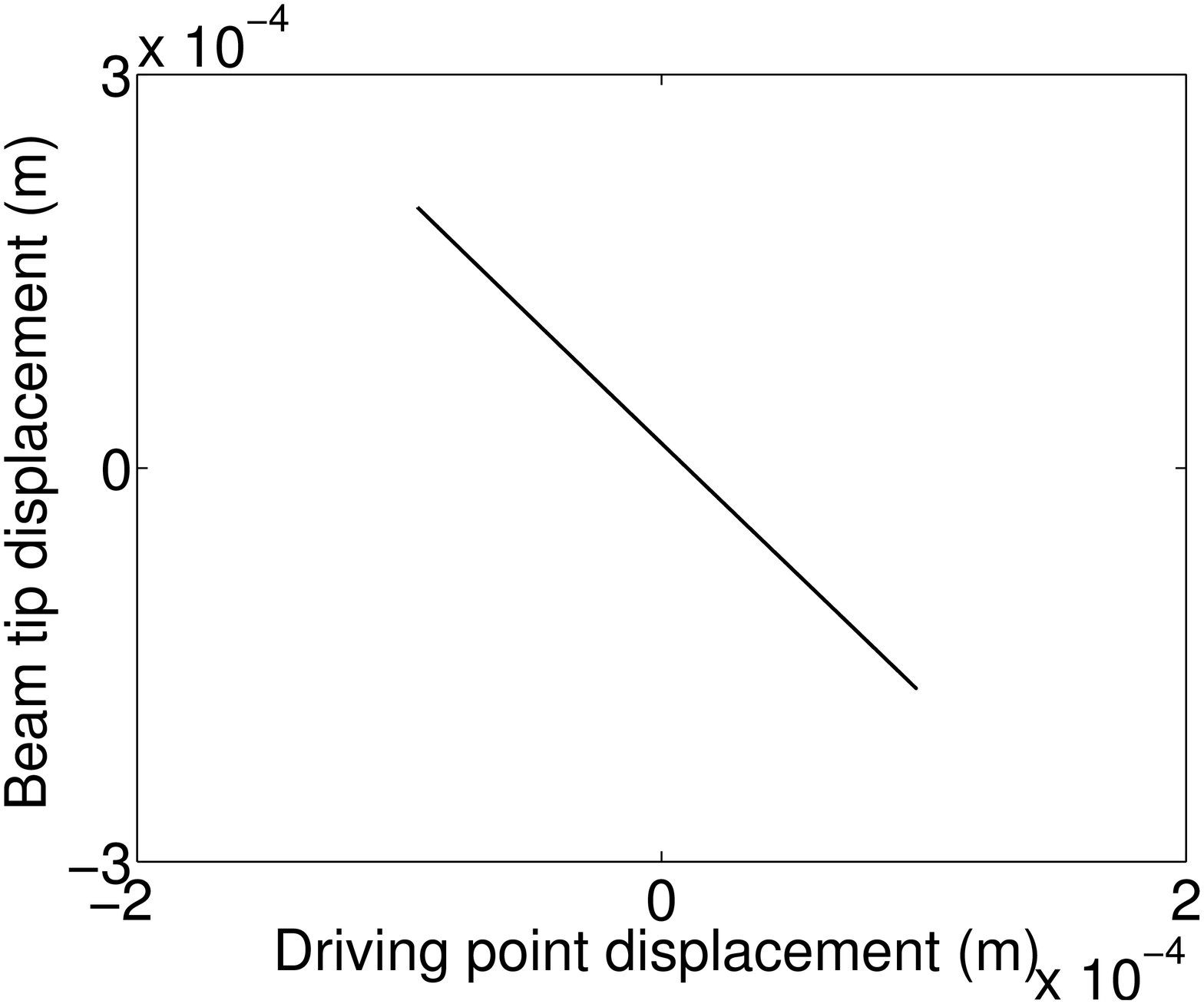}} &
\subfloat[]{\label{NLBeam_NNM2_ConfSpace2}\includegraphics[width=75mm]{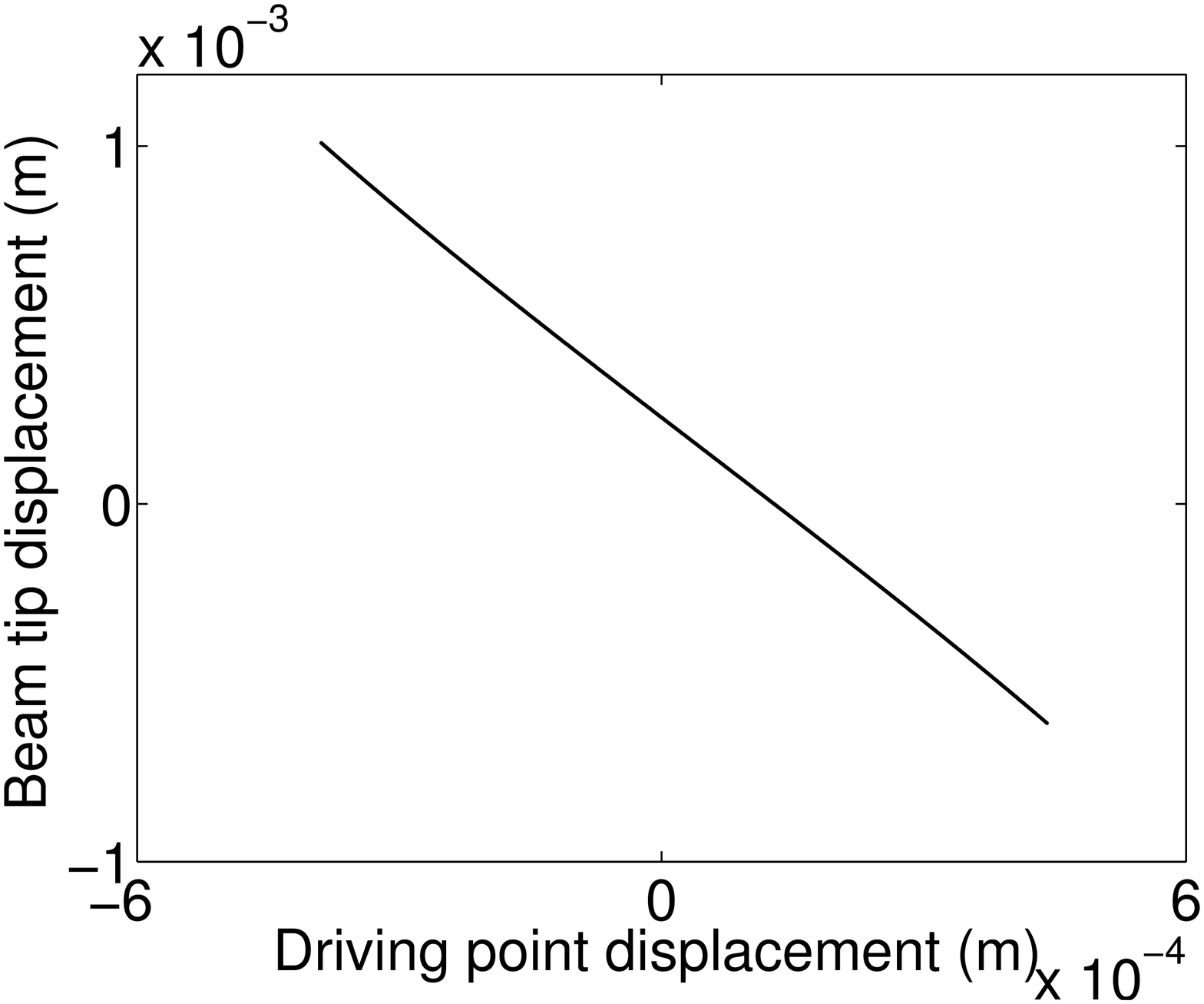}} \\
\end{tabular}
\caption{Second NNM of the nonlinear beam represented in configuration space at two amplitude levels. The configuration space is constructed using the displacements measured at the driving point and the main beam tip. (a) 0.2 $mm$ ; (b) 1 $mm$.}
\label{Fig:NLBeam_NNM2_ConfSpace}
\end{center}
\end{figure}

The results in Fig.~\ref{Fig:NLBeam_NNM1_NPS} prove the validity of the NNM identification methodology for strongly nonlinear regimes of motion. Specifically, positive and negative variations of the natural frequency of approximately 4 and 1 $\%$, respectively, are found to be accurately captured in this plot. These variations correspond to an amplitude of motion of 1 $mm$ at the main beam tip, which is twice the thickness of the thin beam. The importance of nonlinearity in the beam dynamics is well evidenced in Fig.~\ref{Fig:NLBeam_NNM1_ConfSpace}~(a -- b), where the first NNM of the structure is represented in the configuration space for amplitudes of 0.2 and 1 $mm$, respectively. The configuration space is spanned in this figure by the displacements measured at nodes 4 and 14, \textit{i.e.} at the driving point and the main beam tip, respectively. One observes that, at low amplitude level in Fig.~\ref{Fig:NLBeam_NNM1_ConfSpace}~(a), the NNM is a straight line, whereas it corresponds to a curved line for high amplitudes in Fig.~\ref{Fig:NLBeam_NNM1_ConfSpace}~(b), revealing the appearance of harmonics in the time series. The asymmetry in the nonlinear restoring force in the system is also clearly visible. A similar analysis is achieved in Fig.~\ref{Fig:NLBeam_NNM2_ConfSpace}~(a -- b) for the second NNM of the beam. These two graphs show that, owing to the displacement nature of the involved nonlinearity, higher-frequency modes are less impacted by harmonic distortions, and translate into straight lines in configuration space even for large amplitudes of motion. In this study, the amplitude interval over which the continuation was performed was merely selected by observing that the maximum amplitude of displacement recorded in Section~\ref{Sec:Demo_FNSI} at the main beam tip under multisine forcing was of the order of 1 $mm$. However, a rigorous evaluation of the validity ranges of identified frequency-amplitude plots deserves more investigation.

\section{Comparison with NNMs identified using nonlinear phase resonance}\label{Sec:Comparison}

\begin{figure}[p]
\begin{center}
\begin{tabular}{c c}
\subfloat[]{\label{NLBeam_NNM1_Appropriation_MIF}\includegraphics[width=75mm]{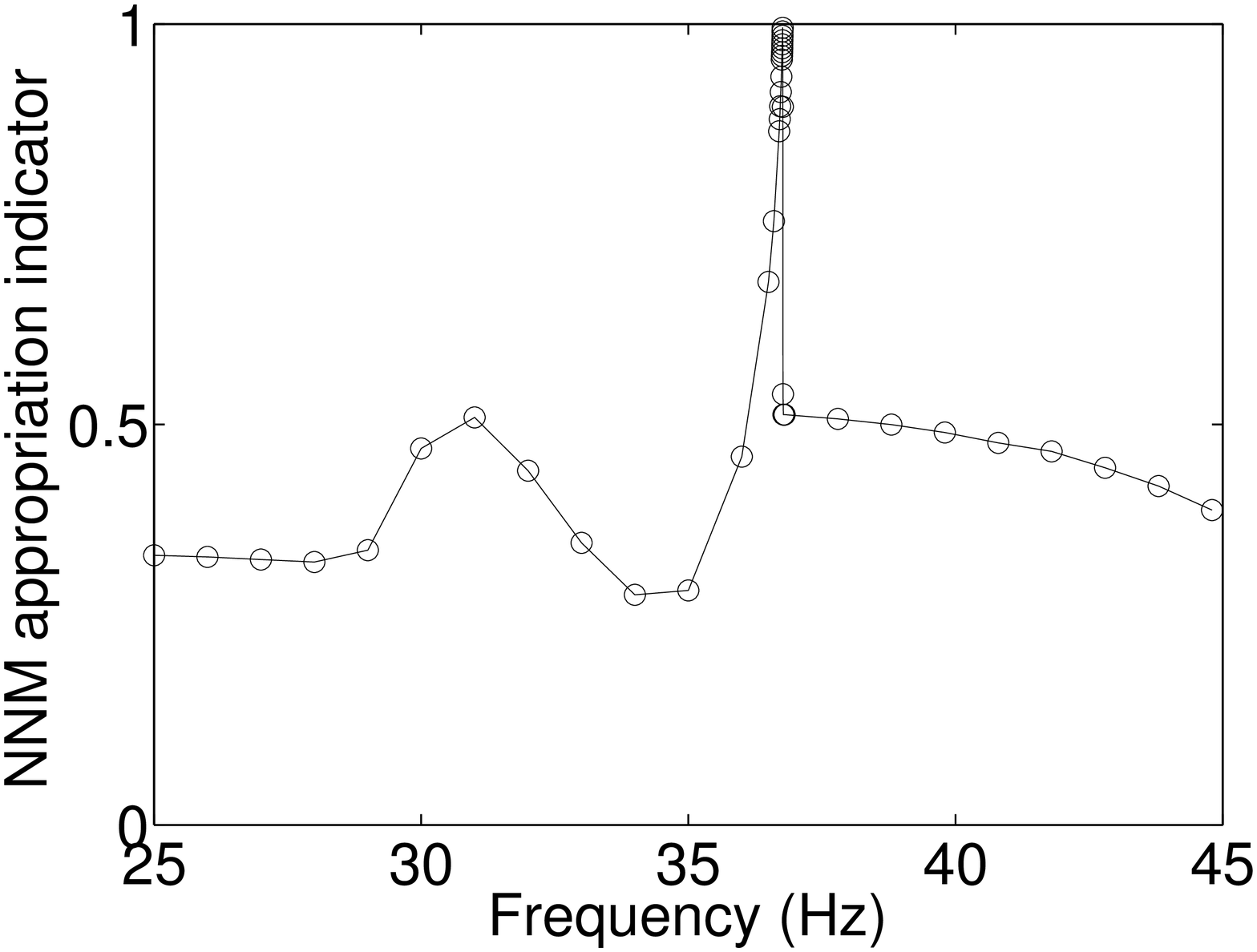}} &
\subfloat[]{\label{NLBeam_NNM1_Appropriation_Response}\includegraphics[width=75mm]{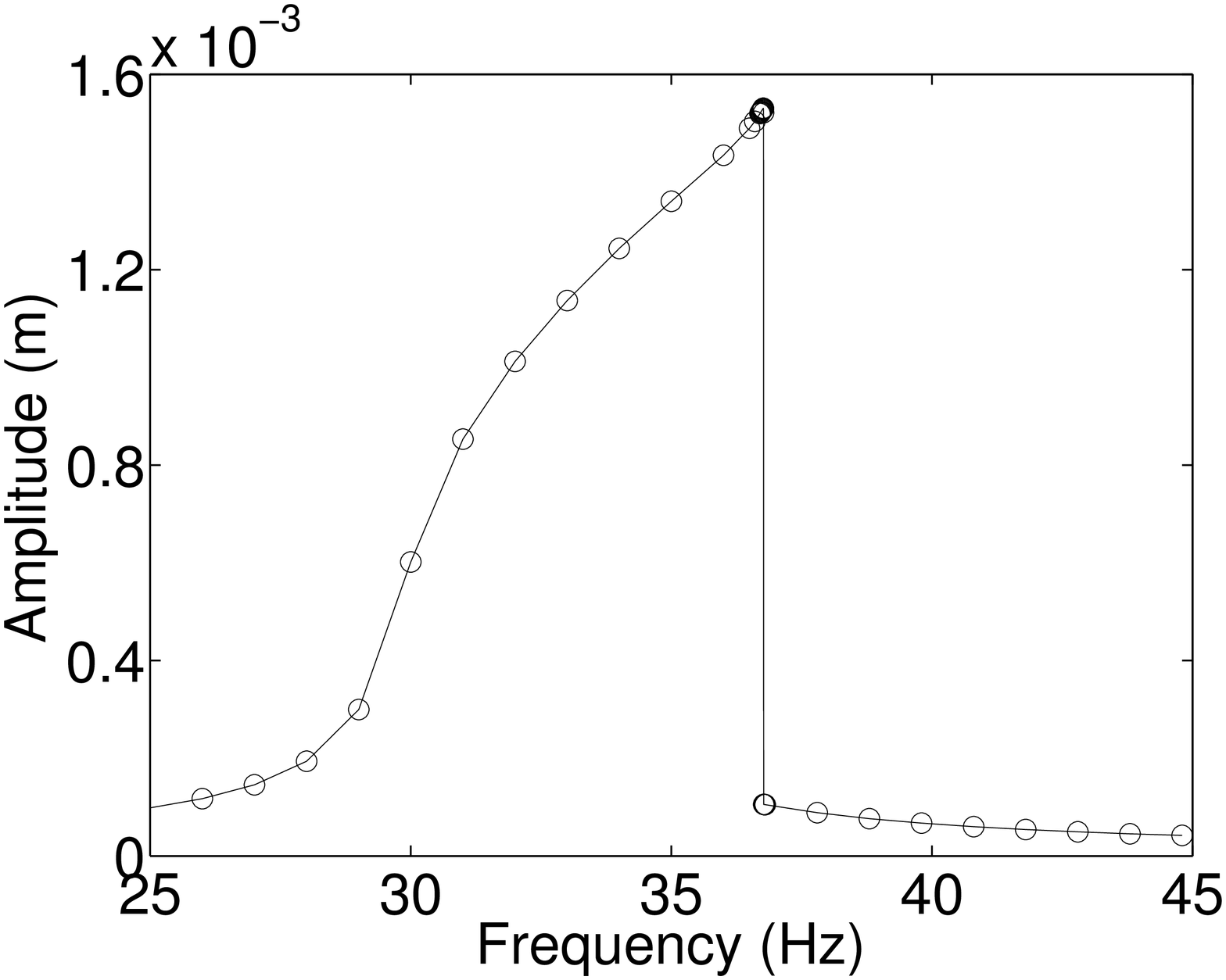}} \\
\end{tabular}
\caption{Appropriation of the first NNM of the beam structure in the nonlinear phase resonance method. (a) NNM appropriation indicator; (b) amplitude of the response at the main beam tip.}
\label{Fig:NLBeam_NNM1_Appropriation}
\end{center}
\end{figure}

As depicted in Fig.~\ref{Fig:NPR_Methodology}, the first step of the nonlinear phase resonance testing procedure is the isolation of the NNM of interest. To this end, a 3 $N$ sine signal is applied to node 4 of the structure. The frequency of this stepped-sine excitation is tuned until the force appropriation indicator derived from the generalized phase lag quadrature criterion is equal to 1~\cite{PRM_Part2}. Fig.~\ref{Fig:NLBeam_NNM1_Appropriation}~(a) shows that NNM appropriation is achieved at 36.8 $Hz$ for the first beam mode. The corresponding amplitude of the forced response at beam tip in Fig.~\ref{Fig:NLBeam_NNM1_Appropriation}~(b) depicts the distorted frequency response of the mode and the sudden jump occurring as soon as resonance is passed.

\begin{figure}[p]
\begin{center}
\begin{overpic}[width=140mm]
		{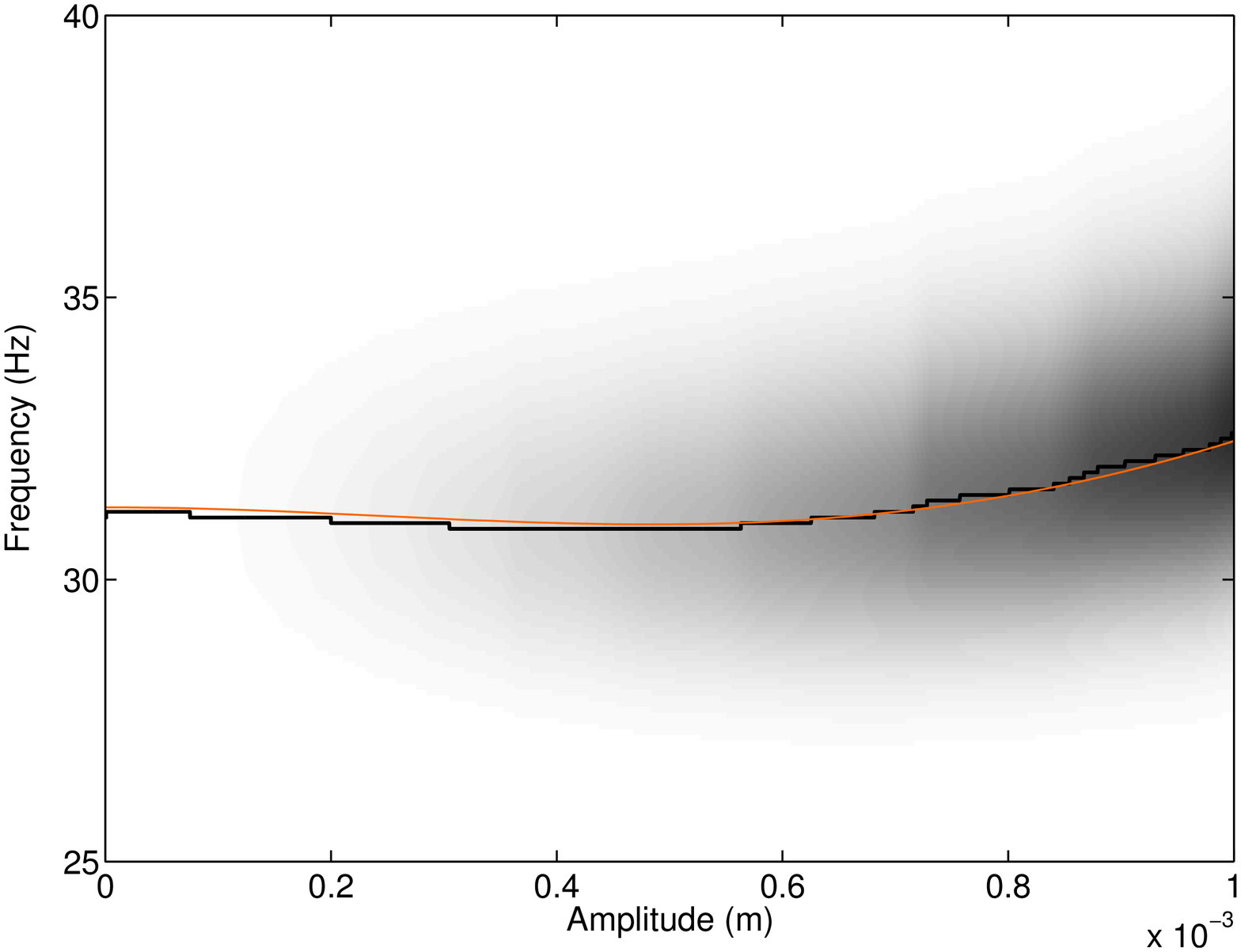}
		\put(110,580){\includegraphics[width=25mm]{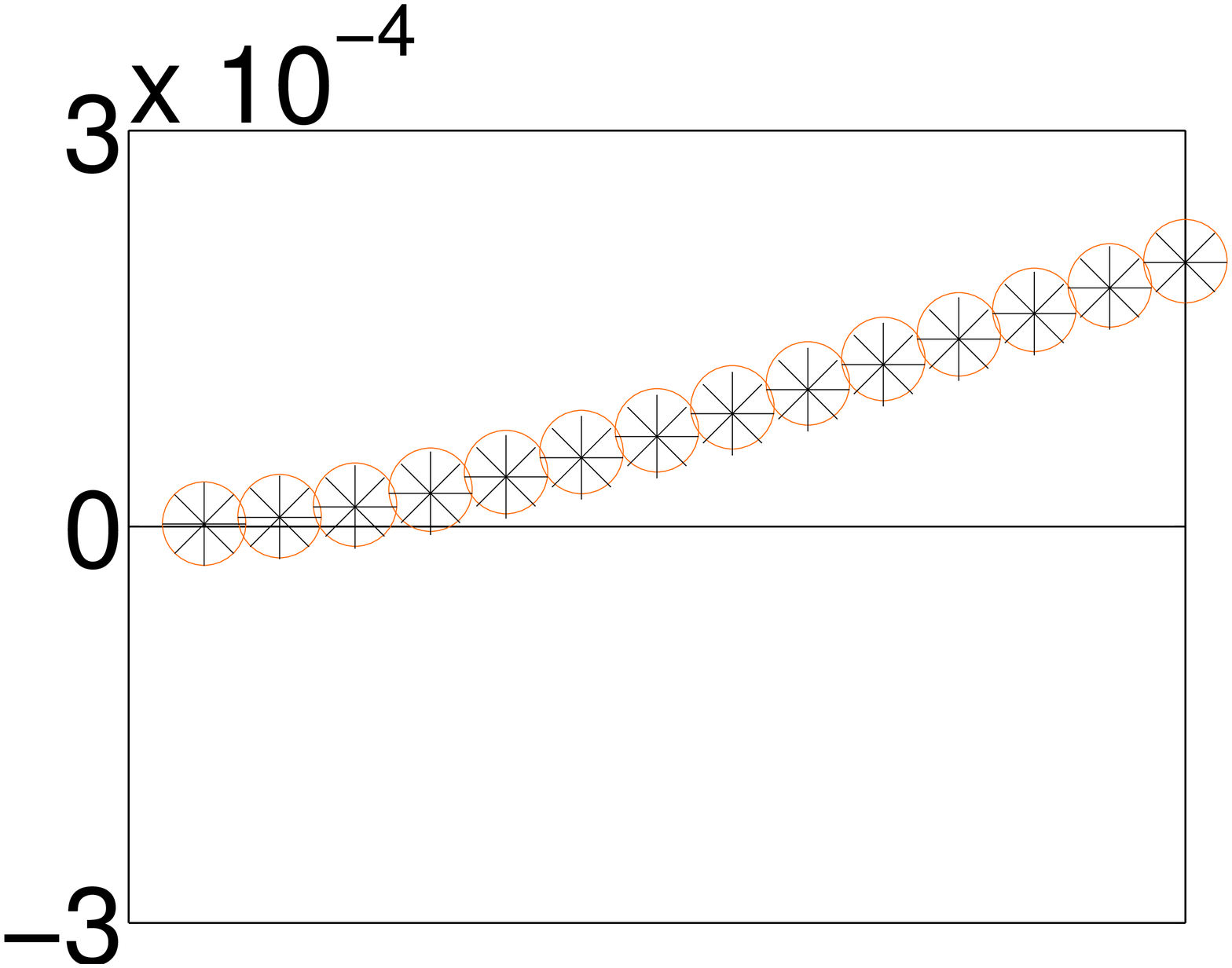}}
		\put(315,580){\includegraphics[width=25mm]{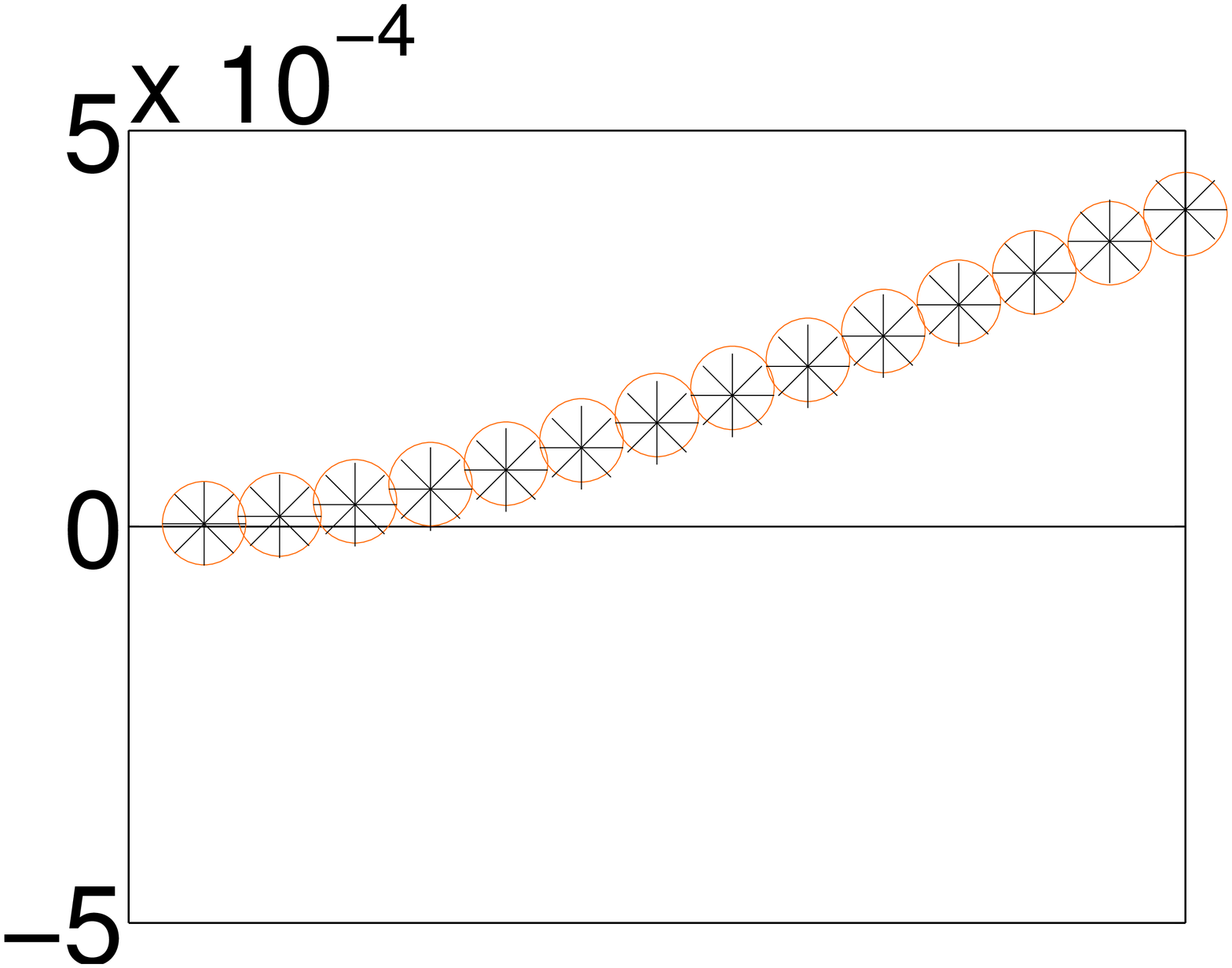}}
		\put(520,580){\includegraphics[width=25mm]{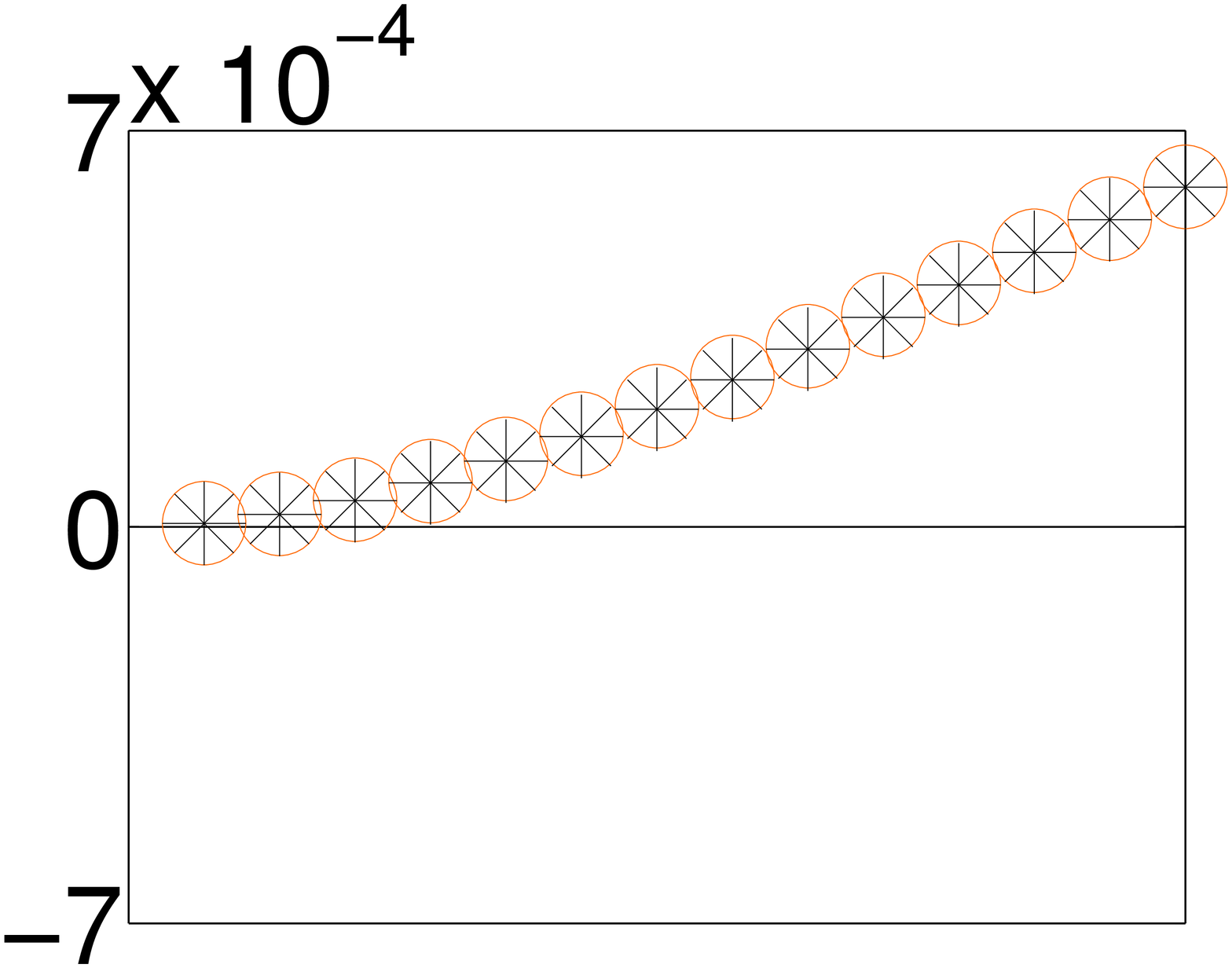}}
		\put(725,580){\includegraphics[width=25mm]{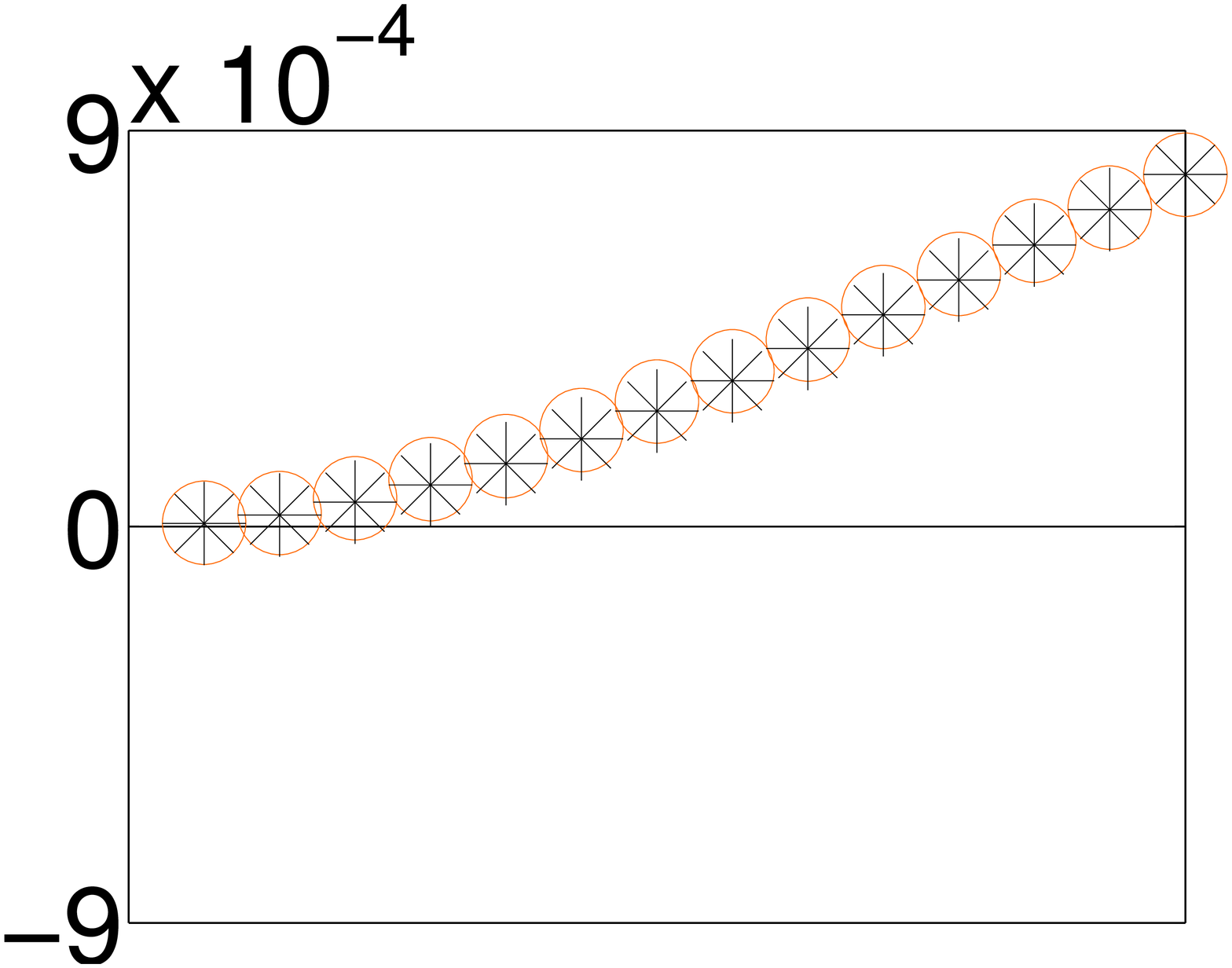}}
		\put(129,555){\footnotesize (a)}
		\put(334,555){\footnotesize (b)}
		\put(539,555){\footnotesize (c)}
		\put(744,555){\footnotesize (d)}
\end{overpic}
\caption{Decay along the first NNM branch calculated using time-frequency analysis (in black) and corresponding NNM obtained in Section~\ref{Sec:Demo_Continuation} using the proposed identification methodology (in orange). The NNM shapes (displacement amplitudes of the main beam) at four amplitude levels, namely 0.2, 0.4, 0.6 and 0.8 $mm$, are inset in (a -- d).}
\label{Fig:NLBeam_NNM1_NPRS}
\end{center}
\end{figure}

\begin{figure}[p]
\begin{center}
\begin{tabular}{c c}
\subfloat[]{\label{NLBeam_NNM2_Appropriation_MIF}\includegraphics[width=75mm]{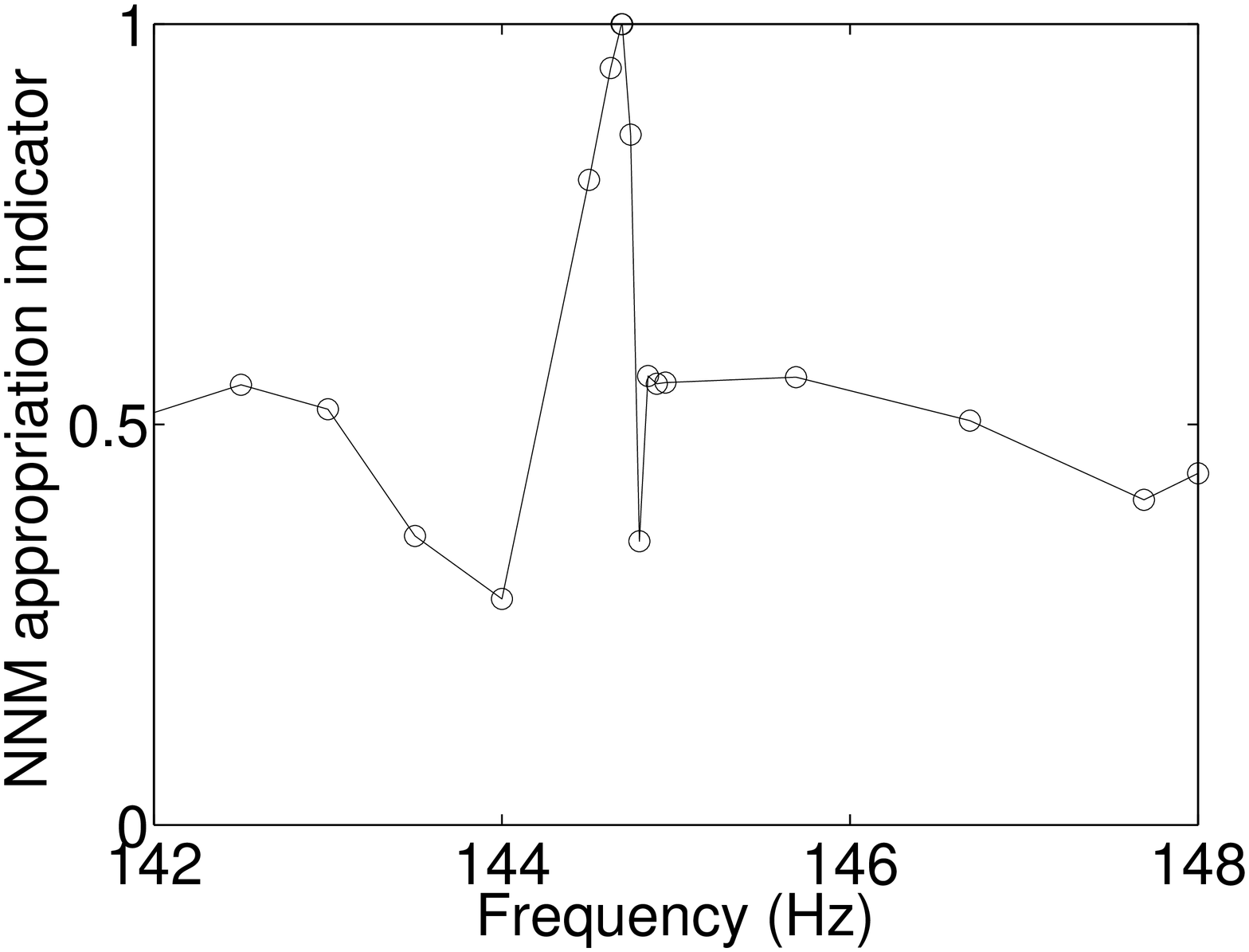}} &
\subfloat[]{\label{NLBeam_NNM2_Appropriation_Response}\includegraphics[width=75mm]{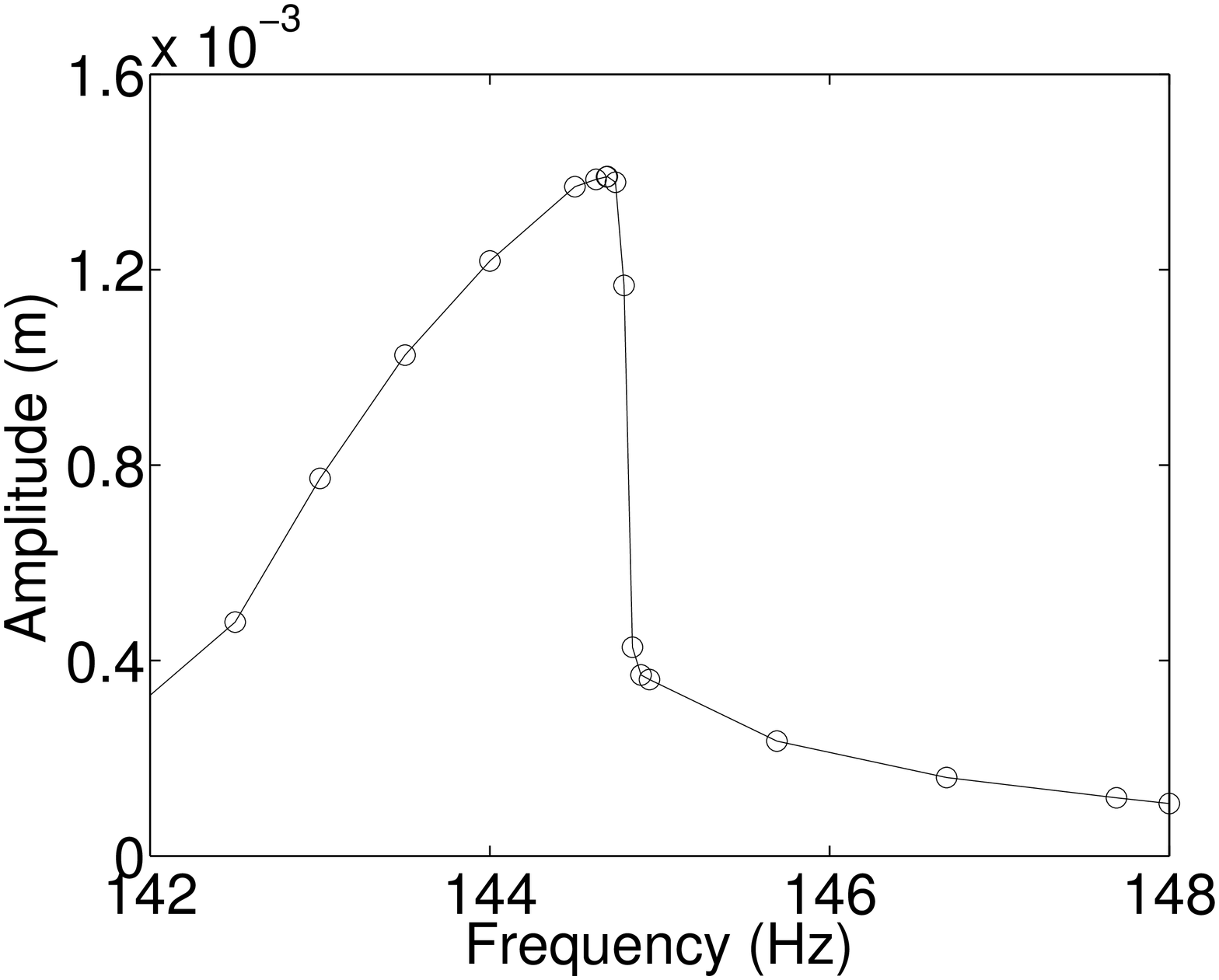}} \\
\end{tabular}
\caption{Appropriation of the second NNM of the beam structure in the nonlinear phase resonance method. (a) NNM appropriation indicator; (b) amplitude of the response at the main beam tip.}
\label{Fig:NLBeam_NNM2_Appropriation}
\end{center}
\end{figure}

\begin{figure}[p]
\begin{center}
\begin{overpic}[width=140mm]
		{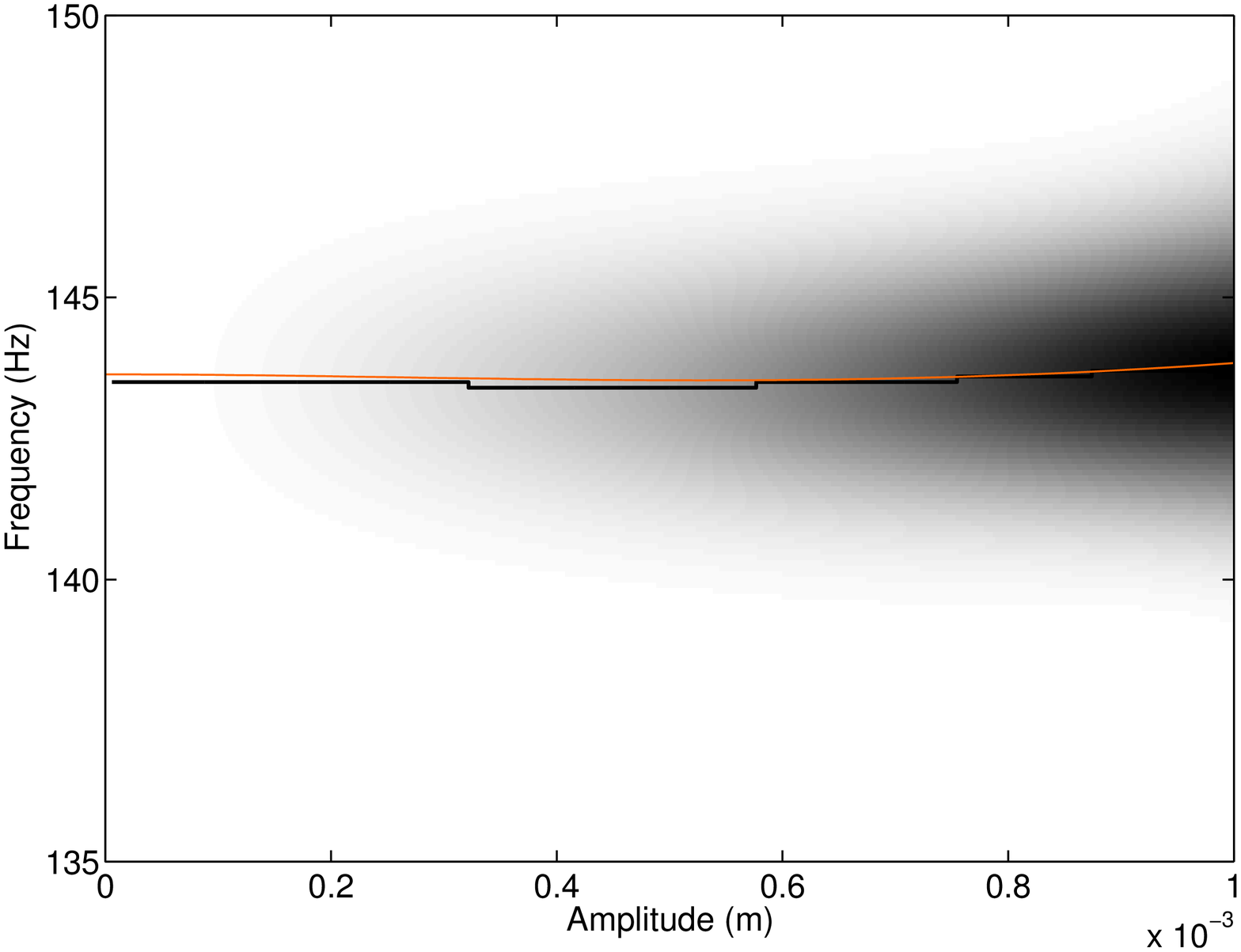}
		\put(110,130){\includegraphics[width=25mm]{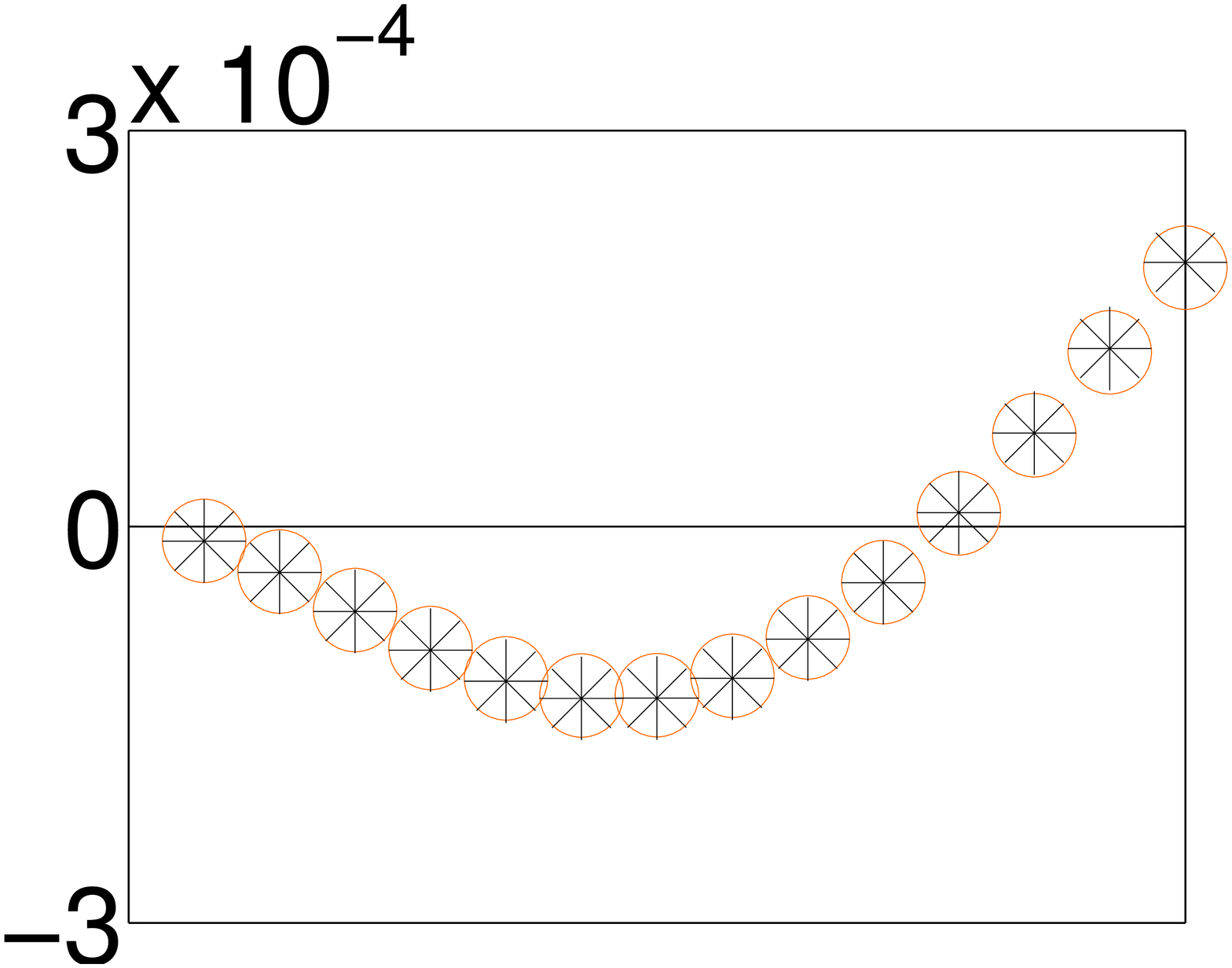}}
		\put(315,130){\includegraphics[width=25mm]{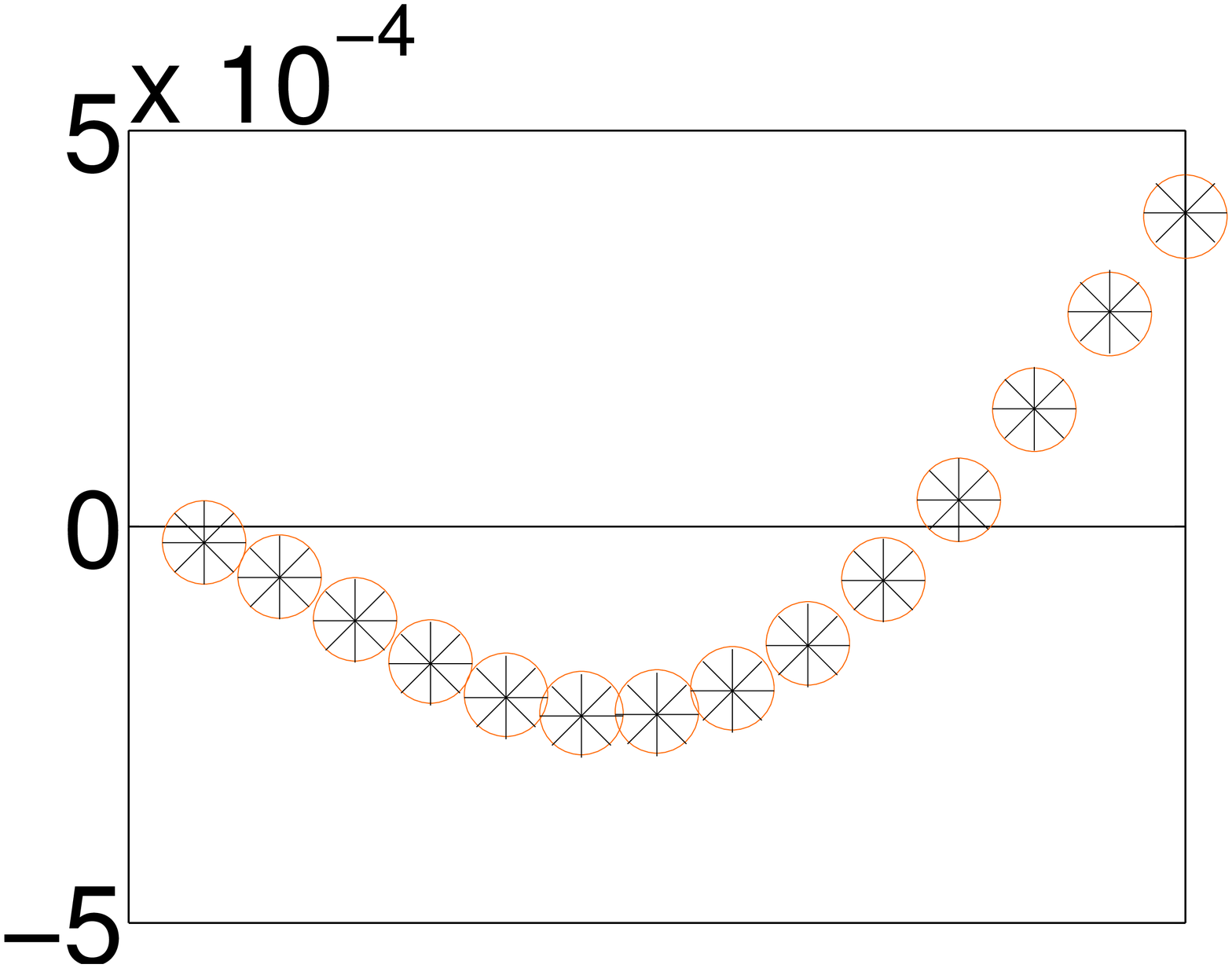}}
		\put(520,130){\includegraphics[width=25mm]{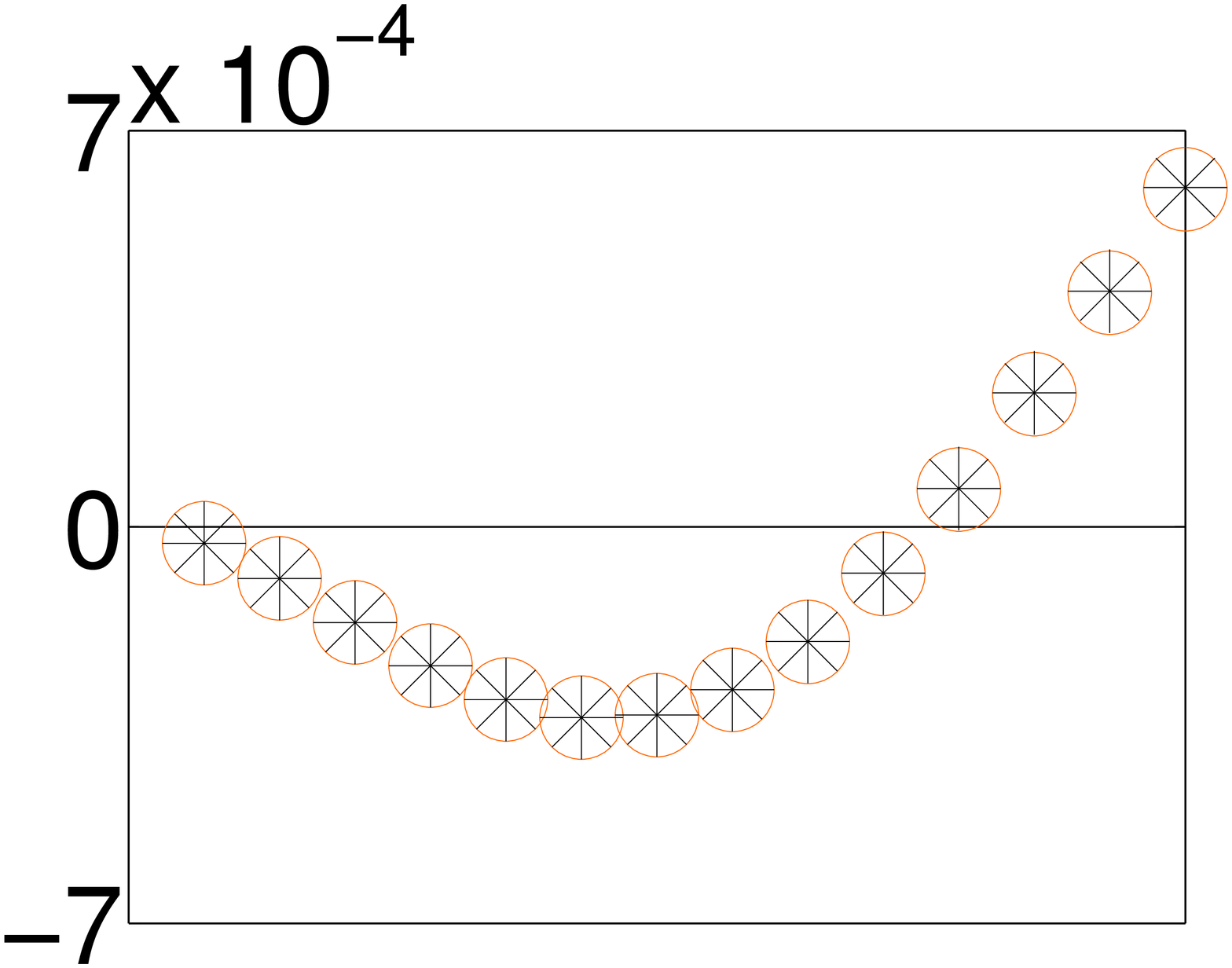}}
		\put(725,130){\includegraphics[width=25mm]{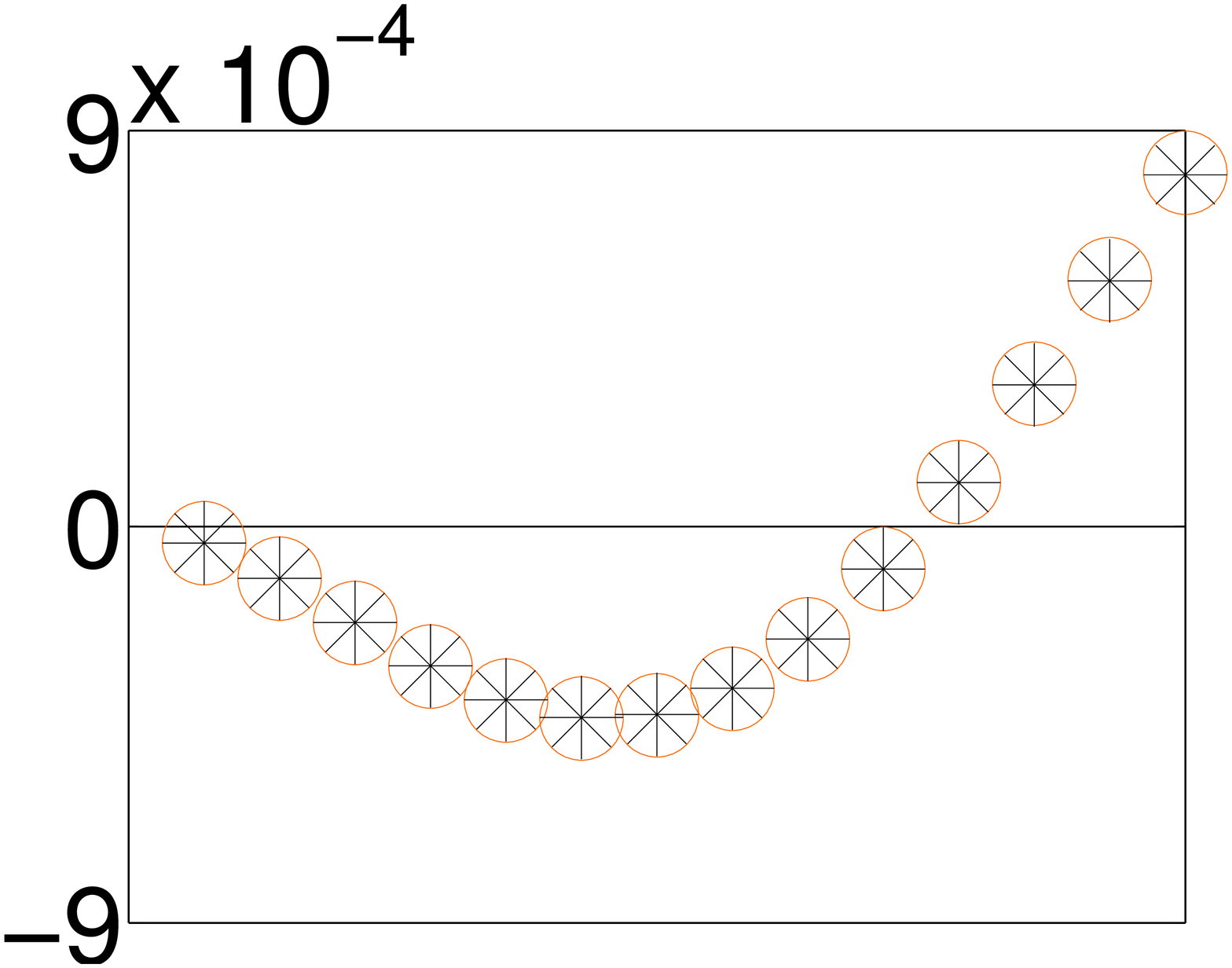}}
		\put(129,105){\footnotesize (a)}
		\put(334,105){\footnotesize (b)}
		\put(539,105){\footnotesize (c)}
		\put(744,105){\footnotesize (d)}
\end{overpic}
\caption{Decay along the second NNM branch calculated using time-frequency analysis (in black) and corresponding NNM obtained in Section~\ref{Sec:Demo_Continuation} using the proposed identification methodology (in orange). The NNM shapes (displacement amplitudes of the main beam) at four amplitude levels, namely 0.2, 0.4, 0.6 and 0.8 $mm$, are inset in (a -- d).}
\label{Fig:NLBeam_NNM2_NPRS}
\end{center}
\end{figure}

When the considered NNM is appropriated, the second step of the procedure turns off the excitation in order to observe the free decay of the system along the NNM branch. A time-frequency analysis of the decaying time response is then carried out to extract the frequency-energy dependence of the mode. This is achieved in Fig.~\ref{Fig:NLBeam_NNM1_NPRS} where the wavelet transform of the displacement measured at beam tip is represented. The ridge of the wavelet, \textit{i.e.} the locus of maximum amplitude with respect to frequency, is presented as a black line, and is seen to closely coincide with the NNM identified in the previous section and plotted in orange. The comparable accuracy of the phase separation and phase resonance approaches is confirmed by the modal shapes superposed at four amplitude levels in Fig.~\ref{Fig:NLBeam_NNM1_NPRS}. The results in this figure, together with the analysis of the second NNM appropriation in Figs.~\ref{Fig:NLBeam_NNM2_Appropriation} and~\ref{Fig:NLBeam_NNM2_NPRS}, clearly confirm the accuracy of this new nonlinear phase separation technique. In summary, Table~\ref{Table:Comp} lists the strengths and limitations of the two methodologies.

\begin{table}[ht]
\vspace*{1cm}
\begin{center}
\begin{tabular*}{1.00\textwidth}{@{\extracolsep{\fill}} c c}
\hline
\textbf{Nonlinear phase}  & \textbf{Nonlinear phase} \\
\textbf{separation method} & \textbf{resonance method} \\
& \\
Fast & Time-consuming \\
(multiple NNMs & (one NNM at a time) \\
\vspace{3mm}
simultaneously) &  \\
Need of an & Model-free \\
\vspace{3mm}
experimental model & \\
Classical random & Harmonic \\ 
\vspace{3mm}
excitation can be utilized & forcing must be tuned \\
Nonlinear components must & Shaker must \\
\vspace{3mm}
be instrumented on both sides & be turned off \\
Nonlinearity characterization & Limited information needed \\
\vspace{3mm}
is required & about the nonlinearities \\
 & \\
\hline
\end{tabular*}
\caption{Comparison of the strengths and limitations of the two methodologies.} 
\label{Table:Comp}
\end{center}
\end{table}

\newpage
\section{Conclusion}\label{Sec:Conclusion}

The present paper extracted the nonlinear normal modes (NNMs) of vibrating systems from measurements collected under broadband forcing. A key feature of the proposed method is that it makes no assumption as to the strength of the nonlinearities and the modal couplings. Together with the previously-developed NNM identification method based on stepped-sine forcing, they provide a rigorous generalization of modal testing to nonlinear systems.
 
The accuracy of the method was demonstrated numerically considering important noise perturbations and no prior knowledge about the nonlinearities in the system, which paves the way for a future experimental validation of the method. In its current state, the method can only handle stiffness nonlinearities. Further research will address this limitation by computing damped NNMs from the experimentally-derived state-space model using, \textit{e.g.}, the computational technique proposed in Ref.~\cite{Renson_NNM_FEM}.

\section*{Acknowledgments}

The author J.P. Noël is a Postdoctoral Researcher of the \textit{Fonds de la Recherche Scientifique -- FNRS} which is gratefully acknowledged. The author L. Renson is a Marie-Curie COFUND Postdoctoral Fellow of the University of Liège, co-funded by the European Union. The authors C. Grappasonni and G. Kerschen would finally like to acknowledge the financial support of the European Union (ERC Starting Grant NoVib 307265).


\end{document}